\documentclass
{article}
\usepackage{amsmath,amssymb}
\usepackage{graphicx,epsf}
\title
{Hilbert functions of points on\\
	Schubert varieties in Orthogonal Grassmannians
	}                    
\author{K.~N.~Raghavan\\
	Institute of Mathematical Sciences\\
        C.~I.~T.~Campus,  Taramani\\
	Chennai 600\,113, INDIA\\
	{\sffamily email\/}: {\ttfamily knr@imsc.res.in}
        \and
        Shyamashree Upadhyay\\
	Chennai Mathematical Institute\\
	Plot No. H1, SIPCOT IT Park\\
	Padur Post, Siruseri 603\,103,	Tamilnadu, INDIA\\
	{\sffamily email\/}: {\ttfamily shyama@cmi.ac.in}}
\date{04 April 2007}

\usepackage{ifthen}     

\newcommand\finalized{yes}

\newcommand\ignore[1]{}
\usepackage{color}
\providecommand\wantcolor{yes}   %
\ifthenelse{\equal{\wantcolor}{no}}{\renewcommand\color[1]{}}{}
\definecolor{backgroundyellow}{cmyk}{.2,.1,.8,.2}
\definecolor{backgroundblue}{rgb}{0,0,1}
\definecolor{backgroundred}{rgb}{1,0,0}
\definecolor{backgroundmagenta}{cmyk}{0,1,0,0}
\definecolor{GoodForInverseVideo}{rgb}{.6,.8,1}
\definecolor{myyellow}{rgb}{1,1,0}

\newcommand\mypart[1]{\vfill\eject\part{#1}}
\newcommand\mysection{\section}
\newcommand\mysectionstar[1]{\section*{#1}}
\newcommand\mysubsection{\subsection}
\newcommand\mysubsectionstar{\subsection*}
\newcommand\mysubsubsection[1]{%
		\subsubsection{\sffamily\upshape\mdseries #1}}
\newcommand\mysss{\mysubsubsection}

\newtheorem{annotation}{Annotation}[subsection]
\newtheorem{theorem}[annotation]{
		Theorem}
\newtheorem{lemma}[annotation]{
		Lemma}

\newtheorem{corollary}[annotation]{
		Corollary}
\newtheorem{proposition}[annotation]{
		Proposition}
\newtheorem{example}[annotation]{
		Example}
\newcommand\bexample{\begin{example}\begin{rm}}
\newcommand\eexample{\end{rm}\hfill$\Box$\end{example}}
\newtheorem{notation}[annotation]{
		Notation}
\newcommand\bnotation{\begin{notation}\begin{rm}}
\newcommand\enotation{\end{rm}\end{notation}}

\newtheorem{remark}[annotation]{
		Remark}
\newcommand\bremark{\begin{remark}\begin{sffamily}\begin{upshape}}
\newcommand\eremark{\end{upshape}\end{sffamily}\end{remark}}
\newenvironment{myproof}{%
\par\noindent{\scshape Proof:}\begin{rm}}{\hfill$\Box$\end{rm}\newline}
\newenvironment{myproofnobox}{%
\par\noindent{\scshape Proof:}\begin{rm}}{\end{rm}\hfill\newline}

\numberwithin{equation}{subsection}
\numberwithin{figure}{subsection}

\providecommand\finalized{no}
\ifthenelse{\equal{\finalized}{yes}}{}{\marginparwidth=2cm}
\ifthenelse{\equal{\finalized}{yes}}%
{\newcommand\mylabel[1]{\label{#1}}}%
{\newcommand\mylabel[1]{\label{#1}\marginpar{[{\ttfamily\upshape\tiny #1}]}}}
\usepackage{ragged2e}
\ifthenelse{\equal{\finalized}{yes}}%
{\newcommand\checked[1]{}}%
{\newcommand\checked[1]{\marginpar{[{\ttfamily\upshape\tiny CHECKED: #1}]}}}
\ifthenelse{\equal{\finalized}{yes}}%
{\newcommand\spellchecked[1]{}}%
{\newcommand\spellchecked[1]{\marginpar{[{\ttfamily\upshape\tiny SPELLCHECKED: #1}]}}}

\providecommand\version{public}   
\ifthenelse{\equal{\version}{public}}%
{\newcommand\mcomment[1]{}}%
{\newcommand\mcomment[1]{\marginpar{{\raggedright\sffamily\upshape\tiny #1}}}}
\ifthenelse{\equal{\version}{public}}%
{\newcommand\fcomment[1]{}}%
{\newcommand\fcomment[1]{\footnote{#1}}}

\newcommand\field{\mathfrak{k}}

\newcommand\idd{I(d,2d)}
\newcommand\id{I(d)}
\newcommand\isubn{I_n}
\newcommand\modV{(V)}
\newcommand\misd{\mathfrak{M}_d\modV}
\newcommand\miso{\misd}
\newcommand\sov{{\rm SO}\modV}

\newcommand\lift[1]{\widetilde{#1}}
\newcommand\liftn{\lift{n}}
\newcommand\liftV{\lift{V}}
\newcommand\lifte{\lift{e}}
\newcommand\lifti{\lift{i}}

\newcommand\misoliftV{\mathfrak{M}_{d+1}(\liftV)}
\newcommand\soliftV{{\rm SO}(\liftV)}
\newcommand\liftT{\lift{T}}
\newcommand\liftB{\lift{B}}

\newcommand\pointe{\mathfrak{e}}

\newcommand\yvw{Y(w)}
\newcommand\affinev{\mathbb{A}\modv}
\newcommand\kyvw{\field[\yvw]}

\newcommand\ortho{\ensuremath{\mathfrak{O}}}
\newcommand\ophi{\ensuremath{\ortho\phi}}
\newcommand\opi{\ensuremath{\ortho\pi}}

\newcommand\pos{\mathfrak{N}}
\newcommand\roots{\mathfrak{R}}

\newcommand\modv{}

\newcommand\modm{(m)}
\newcommand\modvstar{(v^*)}

\newcommand\sv{S\modv}
\newcommand\svw{S^w\modv}
\newcommand\svwm{\svw\modm}

\newcommand\sm{SM}
\newcommand\smv{\sm\modv}
\newcommand\smvadv{\sm_v\modv}
\newcommand\smvw{\sm^w\modv}
\newcommand\smvv{\sm^v\modv}
\newcommand\smvwm{\smvw\modm}
\newcommand\smvvm{\smvv\modm}
\newcommand\smvadvwm{\sm^w_v\modv\modm}

\newcommand\uv{U\modv}
\newcommand\uvm{\uv\modm}
\newcommand\tv{T\modv}
\newcommand\tvm{\tv\modm}
\newcommand\tvwm{T^w\modv\modm}

\newcommand\pass[1]{\widetilde{#1}}

\newcommand\upass{\pass{u}}

\newcommand\idstar{I(d)^*}
\newcommand\vsp{\pass{v^*}}
\newcommand\modvsp{(\vsp)}
\newcommand\smvspvsp{\sm_{\vsp}\modvsp}
\newcommand\modddag{(\ddag)}
\newcommand\smvspvspddag{\smvspvsp\modddag}
\newcommand\smvspvspmodddag{\smvspvsp\modddag}
\newcommand\tvspddag{T\modvsp\modddag}

\newcommand\st{\,|\,}

\newcommand\vchain{$v$-chain}
\newcommand\vchains{$v$-chains}

\newcommand\posv{\ortho\pos\modv}
\newcommand\oroots{\ortho\roots}
\newcommand\rootsv{\ortho\roots\modv}
\newcommand\andposv{\pos\modv}
\newcommand\androotsv{\roots\modv}

\newcommand\diag{\mathfrak{d}}
\newcommand\diagv{\diag\modv}

\newcommand\vdeg[1]{\textup{$v$-degree}(#1)}
\newcommand\degree[1]{\textup{degree}(#1)}

\newcommand\hash{\#}

\newcommand{\mon}{\mathfrak{S}}
\newcommand\monw{\mon_w}
\newcommand{\mont}{\mathfrak{T}}
\newcommand{\monwsupj}{\monw^j}
\newcommand{\monwsubjjpone}{\mon_{w,j,j+1}}
\newcommand{\wsupj}{w^j}
\newcommand{\wsupjptwo}{w^{j+2}}
\newcommand{\wsubjjpone}{w_{j,j+1}}
\newcommand{\monu}{\mathfrak{U}}

\newcommand{\montwjjpone}{\mathfrak{T}_{w,j,j+1}}
\newcommand\montwjjponehashstar{(\montwjjpone\cup\montwjjpone^\hash)^\star}
\newcommand\montwjjponehash{\montwjjpone\cup\montwjjpone^\hash}
\newcommand\montwjjponestar{\montwjjpone^\star}
\newcommand\montwstar{\mont_w^\star}

\newcommand\monj{\mon_j}
\newcommand\monjprime{\mon_j^\prime}
\newcommand\monjponeprime{\mon_{j+1}^\prime}
\newcommand\monjext{\monj\textup{(ext)}}
\newcommand\monjpone{\mon_{j+1}}
\newcommand\monjponeext{\monjpone\textup{(ext)}}
\newcommand\monjjpone{\mon_{j,j+1}}
\newcommand\monjjponeext{\monjjpone\textup{(ext)}}
\newcommand\monjprimeup{\mon_{j}'(\textup{up})}
\newcommand\monjponeprimeup{\mon_{j+1}'(\textup{up})}
\newcommand\block{\mathfrak{B}}
\newcommand\blockb{\mathfrak{B}}
\newcommand\blockc{\mathfrak{C}}

\newcommand\monjpr{\mon_{j,j+1}^{\textup{pr}}}
\newcommand\monjjponepr{\mon_{j,j+1}^{\textup{pr}}}
\newcommand\monprimejjpone{\mon'_{j,j+1}}

\newcommand\piece{\mathfrak{P}}

\newcommand\pbeta{\piece_\beta}

\newcommand\pbetap{\piece_{\beta'}}
\newcommand\pbstar{\piece_\beta^*}
\newcommand\pbpstar{\piece_{\beta'}^*}

\newcommand\opbstar{\ortho\piece_\beta^*}
\newcommand\opbpstar{\ortho\piece_{\beta'}^*}

\newcommand\monk{\mon_{k}}

\newcommand\monkpr{\mon_k^\textup{pr}}

\newcommand\monkprime{\mon_{k}'}

\newcommand\monkprimeup{\mon_{k}'(\textup{up})}

\newcommand\monkext{\mon_{k}(\textup{ext})}

\newcommand\monjonlypr{\mon_j^\textup{pr}}
\newcommand\monjponepr{\mon_{j+1}^\textup{pr}}

\newcommand\deltaj{\delta_j}

\newcommand\odepth{$\ortho$-depth}
\newcommand\odepths{$\ortho$-depths}
\newcommand\odepthinC[1]{\ensuremath{\ortho\textup{-depth}_C(#1)}}
\newcommand\odepthin[2]{\ensuremath{\ortho\textup{-depth}_{#1}(#2)}}
\newcommand\odepthinmon[1]{\ensuremath{\ortho\textup{-depth}_\mon(#1)}}

\newcommand\depthin[2]{\ensuremath{\textup{depth}_{#1}#2}}

\newcommand\depthinmonC[1]{\ensuremath{\textup{depth}_{\mon_C}(#1)}}

\newcommand\even{\textup{(even)}}
\newcommand\odd{\textup{(odd)}}
\newcommand\up{\textup{(up)}}

\newcommand\down{\textup{(down)}}
\newcommand\ext{\textup{(ext)}}

\newcommand\myrho{q}

\newcommand\xjjpone{x_{j,j+1}}
\newcommand\pr{\textup{pr}}

\newcommand\monsupjjpone{\mon^{j,j+1}}
\newcommand\monsupjptwojpthree{\mon^{j+2,j+3}}
\newcommand\monsupthreefour{\mon^{3,4}}

\newcommand\start{\textup{(start)}}
\newcommand\finish{\textup{(finish)}}
\newcommand\bstart{\beta\start}

\newcommand\bfinish{\beta\finish}

\newcommand\monomvw{\mathfrak{Mon}^w\modv}
\newcommand\pathsvw{\mathfrak{Paths}^w\modv}
\newcommand\path{\Lambda}

\usepackage{makeidx}
\ifthenelse{\equal{\finalized}{no}}{\makeindex}{}

\ifthenelse{\equal{\finalized}{no}}{
\usepackage[editors,firstabbrev,miniindex]{authorindex}
\usepackage{multicol}
\renewcommand{\cite}{\aicite}    
\aisee{, see }
}{}
\begin{document}
\maketitle
\begin{abstract} 
A solution is given to the following problem:   
how to compute the multiplicity, or more generally the Hilbert function,
at a point on a Schubert variety in an orthogonal Grassmannian.
Standard monomial theory is applied to translate the problem from
geometry to combinatorics.  The solution of the resulting combinatorial 
problem forms the bulk of the paper.    This approach has been followed
earlier to solve the same problem for the Grassmannian and the symplectic
Grassmannian.   

As an application, we present\mcomment{Change two to one
if Gr\"obner basis is not being included here and other appropriate changes in
this paragraph.   

Added 16 March: Done. }
an interpretation of the multiplicity as the number of
non-intersecting lattice paths of a certain kind.

Taking the Schubert variety to be of a special kind and the point to
be the ``identity coset,'' our problem specializes to a problem about
Pfaffian ideals treatments of which by different methods exist in the
literature.  Also available in the literature is a geometric solution
when the point is a ``generic singularity.''
\end{abstract}


\vfill\eject
\tableofcontents
\vfill\eject
\mysectionstar{Introduction}\mylabel{s.introduction}
\addcontentsline {toc}{section}{Introduction}
\checked{Jan 18, 2007}\spellchecked{Jan 18, 2007}
In this paper the following problem is solved:   given
a Schubert variety in an {\em orthogonal Grassmannian\/} (by which is
meant the variety of isotropic subspaces of maximum possible
dimension of a finite dimensional vector space with a symmetric 
non-degenerate form---see~\S\ref{s.setup} for precise definitions)
and an arbitrary point on the Schubert variety,
how to compute the multiplicity, or more generally the Hilbert function,
of the local ring of germs of functions at that point.
In a sense,  our solution is 
but a translation of the problem:   we do not give closed form formulas
but {\em alternative\/} combinatorial descriptions.    The meaning of
``alternative'' will presently become clear.

The same problem for the Grassmannian was treated
in~\cite{kl1,kratII,kr,k,kl2} and for the symplectic Grassmannian in~\cite{gr}.
The present paper is a sequel to~\cite{kl1,kr,k,kl2,gr}
and toes the same line as them.
In particular, its strategy is borrowed from them and runs as
follows:  first translate the problem from geometry to combinatorics, or,
more precisely, apply {\em standard monomial theory\/} to obtain an
{\em initial\/} combinatorial description of the Hilbert function (the
earliest version of the theory capable of handling Schubert varieties
in an orthogonal Grassmannian is to be found in~\cite{s}); 
then transform the initial combinatorial description to obtain the desired
alternative description.    But that is easier said than done.

While the problem makes sense for Schubert varieties of any kind
and standard monomial theory itself is available 
in great generality~\cite{ls,li},  the translation of the problem from
geometry to combinatorics has been
made---in~\cite{lw}---
only for ``minuscule\footnote{Symplectic Grassmannians are not
minuscule but can be treated as if they were.} generalized Grassmannians.''
Orthogonal Grassmannians being minuscule,  this
translation is available to us and we have an initial combinatorial
description of the Hilbert function.    As to the passage from the initial
to the alternative description---and this is where the content of the
present paper lies---neither the end nor the means is clear at the outset.

The first problem then is to {\em find\/} a good alternative
description.     But how to measure the worth of an alternative description?
The interpretation of multiplicity as the number of
certain non-intersecting lattice paths
(deduced in~\S\ref{s.application} 
from our alternative description)  seems to testify to the 
correctness of our alternative description,  but we are not sure if there
are others that are equally or more correct.

The proof of the equivalence of the initial and alternative combinatorial
descriptions is, unfortunately, a little technically involved.
It builds on the details of the proofs of the corresponding equivalences
in the cases of the Grassmannian and the symplectic
Grassmannian.     
In~\cite{k2} it is shown that the equivalence in
the case of the Grassmannian is a kind of KRS correspondence,
called ``bounded KRS.''  The proof there is short and elegant and it
would be nice to realise the main result of the present paper too in
a similar spirit as a kind of KRS correspondence.

The initial description is in terms of ``standard monomials'' and the
alternative description in terms of ``monomials in roots.''   The
equivalence of the two descriptions thus 
gives a bijective correspondence between standard  monomials and monomials
in roots.  Roughly---but not actually---the correspondence
maps each standard monomial
to its initial term (with respect to a certain monomial order).   Thus
it is natural to wonder whether we can compute the initial ideal of 
the ideal of the tangent cone to the Schubert variety
at the given point.   
\mcomment{Paragraph changed because we are not doing Groebner basis}
We believe that this can be done but that it is far
more involved and difficult than the corresponding computation for
Grassmannians and symplectic Grassmannians (the natural set of
generators of the ideal of the tangent cone do not form a
Gr\"obner basis unlike in those cases).   If all goes well,
the computation will soon appear~\cite{ru}.


Taking the Schubert variety to be of a special kind and the point to be
the ``identity coset,''  our problem specializes to a problem about 
Pfaffian ideals considered in~\cite{ht,dn}.
On the other side of the spectrum from the identity coset, so to speak,
lie the ``generic singularities,'' points that are generic in the
complement of the open orbit of the stabiliser of the Schubert variety.
For these,  a geometric solution to the problem appears in~\cite{bp}. 

Given that our solution of the problem is but a translation,  it makes sense
to ask if one can extract more tangible information---closed form formulas
for example---from our alternative description.    See the papers quoted
in the previous paragraph and also~\cite{gk} 
for some answers in the special cases they consider.
\mysubsectionstar{Organization of the paper}  
\addcontentsline{toc}{subsection}{Organization of the paper}
The table of contents indicates how the paper is organized.
\mcomment{Addition}
There is a brief description at the beginning of every subdivision of the
contents therein.    An index of definitions
and notation is included,  for it would otherwise be difficult to find
the meanings of certain words and symbols.
\mysubsectionstar{Important note added}  
\addcontentsline{toc}{subsection}{Important note added}
The recent article~\cite{in} treats some of the questions addressed here and
some that could be addressed by using the main result proved here.
It includes:
\begin{itemize}
\item  an interpretation of the multiplicity similar to ours.
\item  a closed formula for the multiplicity (as a specialization of
a factorial Schur function),   thereby answering the question we raised above.
\item a formula for the restriction to the torus fixed point of the
equivariant cohomology class of a Schubert variety.
\end{itemize}
The approach in~\cite{in} is quite different from ours.
In fact, it is the
opposite of ours in that
it circumvents the lack of results about
initial ideals of tangent cones,  while our prime motivation
is to remedy the lack.    
The starting points in the two approaches are also different:
\cite{in}~takes off from certain results of
Kostant-Kumar and Arabia on equivariant cohomology, while our launchpad
is standard monomial theory.

The appearance of~\cite{in} notwithstanding,  our approach is worthwhile,
for, quite apart from the difference in starting points,   
there is no way,
as far as we can tell, to the Hilbert function via the approach of~\cite{in},
nor to the initial ideal,  both of which are interesting in their own right.




\mypart{The theorem}\mylabel{p.theorem}
Definitions are recalled, the problem formulated, and the theorem stated.
\mysection{The set up}\mylabel{s.setup}
\checked{\S\ref{s.setup}; 21 Nov, 1000~hrs}
\spellchecked{21 Nov, 1000~hrs}
In this section,  we state the problem to be addressed after 
recalling the necessary basic definitions,
make some choices that are convenient for studying the problem,
and see why it is enough to focus
on a particular case of the problem.


\mysubsection{The statement of the problem}\mylabel{ss.problem}%
\mylabel{ss.xwinmiso}
Fix an algebraically closed field of characteristic not
equal to~$2$%
\index{k@$\field$, base field, ($\textup{characteristic}\neq2$)}.
Fix a vector space~$V$\index{V@$V$, vector space of dimension $n$} of finite 
dimension~$n$\index{n@$n:=\dim V$, (even from \S\ref{ss.basicnotn} on)} 
over this field
and a non-degenerate symmetric bilinear form~$\langle\ ,\ \rangle\index{form@$\langle\ ,\ \rangle$, bilinear form on $V$}$ on~$V$.
Let~$d$\index{d@$d$, integral part of $n/2$ (unfortunately also used otherwise)} be the integer such that either~$n=2d$ or~$n=2d+1$.
A linear subspace of $V$ is said to be {\em isotropic\/}%
\index{isotropic subspace}
if the form $\langle\ ,\ \rangle$ vanishes identically
on it. It is elementary to see that an isotropic subspace of~$V$
has dimension at most~$d$ and that every isotropic subspace is
contained in one of dimension~$d$.   Denote by~
$\miso'$\index{mdvprime@$\miso'$} the closed sub-variety of the Grassmannian of $d$-dimensional
subspaces consisting of the points corresponding
to isotropic subspaces.

The orthogonal group ${\rm O}(V)$\index{O(V)@${\rm O}(V)$} of linear automorphisms of $V$
preserving $\langle\ ,\ \rangle$ acts transitively on
$\miso'$, for by Witt's theorem an isometry
between subspaces can be lifted to one of the whole vector space.
If $n$ is odd the special orthogonal
group~${\rm SO}(V)$\index{SO(V)@${\rm SO}(V)$}
(consisting of form preserving linear automorphisms with trivial determinant)
itself acts transitively\fcomment{Proof: Choose basis such that the matrix
of the form has $1$s on the anti-diagonal and $0$s elsewhere.
The transformation whose matrix with respect to this basis 
has $1$s on the diagonal except on row $d+1$ where it is $-1$
and $0$s elsewhere is in the isotropy at $\langle e_1,\ldots,e_d\rangle$.
Thus there is a transformation of determinant~$-1$ in the isotropy at
$\langle e_1,\ldots,e_d\rangle$.} on $\miso'$.   If $n$ is even
the special orthogonal group ${\rm SO}(V)$ does not\fcomment{
Proof:  We will show that the isotropy at a point is contained in $SO$.
Choose basis so that the matrix of the form is 
{\tiny
$\left(\begin{array}{cc}
0 & I \\
I & 0 \end{array}\right)$}.
Given a transformation in the isotropy at $\langle e_1,\ldots,e_d\rangle$ multiply it by 
{\tiny
$\left(\begin{array}{cc}
g & 0 \\
0 & {g^{t}}^{-1} \end{array}\right)$}
to assume it has form 
{\tiny
$\left(\begin{array}{cc}
I & A \\
0 & B \end{array}\right)$}.   The condition that this last matrix is
in the isotropy translates easily to $B=I$ and $A+A^t=0$
(which means that the matrix has determinant~$1$).}
act transitively on $\miso'$, and $\miso'$ has two connected components.
We define the {\em orthogonal \index{orthogonal Grassmannian ($\miso$)}Grassmannian\/}~$\miso$\index{mdv@$\miso$, orthogonal Grassmannian} to be $\miso'$ if
$n$ is odd and to be one of the two components of $\miso'$ if $n$ is even.

The {\em Schubert\index{Schubert varieties} varieties\/} of~$\miso$ are defined to be the $B$-orbit
closures in~$\miso$ (with canonical reduced scheme structure), where
$B$ is a Borel subgroup of~$\sov$.   The choice of~$B$ is immaterial, for
any two of them are conjugate.
The question that is tackled in this paper is this:  given
a point on a Schubert variety in~$\miso$, how to compute the multiplicity
(and more generally, the Hilbert function) of the Schubert variety at
the given point?     The answers are contained in Theorem~\ref{t.main}
and Corollary~\ref{c.t.main}.     But in order to make sense of those
statements,   we need some preparation.



\mysubsection{Some convenient choices}\mylabel{ss.choices}
We now make some choices that are convenient for the study of 
Schubert varieties.    For $k$ an integer such that $1\leq k\leq n$,  
set $k^*:=n+1-k$\index{kstar@$k^*(:=n+1-k)$}.   Fix a basis 
$e_1,~\ldots,~e_{n}\index{eone@$e_1,\ldots,e_n$, a specific basis of $V$}$ 
of $V$ such that
\[
\langle e_i, e_k \rangle = \left\{ \begin{array}{rl}
    1 & \mbox{ if $i=k^*$ }\\ 
    0 & \mbox{ otherwise}\\ \end{array}\right.
 \]
The advantage of this choice is:
the elements of $\sov$ for which each $e_k$ is an eigenvector
form a maximal torus,  and the elements that
are upper triangular with respect to this basis form a Borel subgroup
(a linear transformation is {\em upper triangular\/} if
for each $k$, $1\leq k\leq n$,  the image of $e_k$  under
the transformation is a linear combination of $e_1,\ldots,e_k$).
We denote this maximal torus and this Borel subgroup by $T$\index{T@$T$, a specific maximal torus} and $B$\index{B@$B$, a specific Borel subgroup}
respectively. Our Schubert varieties\index{Schubert varieties} will be orbit closures of this particular 
Borel subgroup~$B$.

The $B$-orbits of $\miso'$ are naturally indexed by its $T$-fixed points:
each orbit contains one and only one such point.     The $T$-fixed points
are evidently of the form $\langle e_{i_1},\ldots,e_{i_d}\rangle$, where
$1\leq i_1<\ldots<i_d\leq n$ and for each $k$, $1\leq k\leq d$, there
does not exist~$j$, $1\leq j\leq d$, such that $i_k^*=i_j$---in other words,
for each $\ell$, $1\leq\ell\leq n$,  such that $\ell\neq\ell^*$,  exactly one
of $\ell$ and $\ell^*$ appears in $\{i_1,\ldots,i_d\}$; in addition,
if~$n$ is odd,  then $d+1$ does not appear in~$\{i_1,\ldots,i_d\}$. Denote the
set of such $d$-element subsets~$\{i_1<\ldots<i_d\}$ 
by $\isubn'$\index{Isubn'@$\isubn'$}.    
We thus have a bijective correspondence between~$\isubn'$ and
the $B$-orbits of $\miso'$.
Each $B$-orbit being irreducible and open in its closure,  it follows
that $B$-orbit closures are indexed by the $B$-orbits.    Thus $\isubn'$
is an indexing set for $B$-orbit closures in $\miso'$.

Suppose that $n$ is even---it will be shown presently that it is enough
to consider this case.
As already observed, $\miso'$ has two connected
components on each of which $\sov$ acts transitively.   The $B$-orbits
belong to one or the other component accordingly as the parity of the
cardinality of the number of entries bigger than~$d$ in the corresponding
element of $\isubn'$.
We take~$\miso$\index{mdv@$\miso$, orthogonal Grassmannian} 
to be the component in which
these cardinalities are even.     We let~$\isubn$\index{Isubn@$\isubn$} 
denote the subset of 
$\isubn'$ consisting of elements for which this cardinality is even.
Schubert varieties in $\miso$ are thus indexed by elements of $\isubn$.

\mysubsection{Reduction to the case $n$ even}\mylabel{ss.neven}
We now argue that it is enough to consider the case~$n$ even.
Suppose that $n$ is odd.   Let $\liftn:=n+1$ and $\liftV$
be a vector space of dimension~$\liftn$ with a non-degenerate symmetric
form.   Let $\lifte_1,\ldots,\lifte_{\liftn}$ be a basis of $\liftV$ as
in~\ref{ss.choices}.    Put $e:=\lifte_{d+1}$ and $f:=\lifte_{d+2}$.
Take $\lambda$ to be an element of the field such that $\lambda^2=1/2$.
We can take $V$ to be the subspace of $\liftV$ spanned by the vectors
$\lifte_1,\ldots,\lifte_{d},\lambda e+\lambda f,\lifte_{d+3},\ldots,
\lifte_{\liftn}$,  and a basis of $V$ to be these vectors in that order.

There is a natural map from $\misoliftV'$ to $\miso$:  intersecting with~$V$ an
isotropic subspace of $\liftV$ of dimension $d+1$ gives an isotropic subspace
of~$V$ of dimension~$d$.  This map is onto, for every isotropic
subspace of $\liftV$ (and hence of $V$) is contained in an isotropic
subspace of $\liftV$ of dimension~$d+1$.   It is also elementary to
see that the map is two-to-one (essentially because in a two-dimensional
space with a non-degenerate symmetric form there are two isotropic lines),
and that the two points in any fiber lie one in each component (there
is clearly an element in ${\rm O}(\liftV)\setminus {\rm SO}(\liftV)$
that moves one element of the fiber to the other,  and so if there
was an element of ${\rm SO}(\liftV)$ that also moved one point to the
other, the isotropy at the point would not be contained in ${\rm SO}(\liftV)$,
a contradiction).

We therefore get a natural isomorphism between
$\misoliftV$ and $\miso$.     We will now show that the $\liftB$-orbits
in $\misoliftV$ correspond under the isomorphism to $B$-orbits of
$\miso$ (we denote by $\liftT$ and $\liftB$ the maximal torus and Borel
subgroups of $\soliftV$ as in \S\ref{ss.choices}).    It will then
follow that Schubert varieties in $\misoliftV$ are isomorphic to those 
in~$\miso$ and the purpose of this subsection will be achieved.

The group $\sov$ can be realized as the subgroup of $\soliftV$ consisting
of the elements that fix~$e-f$.  The isomorphism $\misoliftV\cong\miso$ above is
equivariant for~$\sov$,  and we have  
$\liftT\cap\sov=T$ and $\liftB\cap\sov=B$.
It should now be clear that the preimages in $\misoliftV$ of two elements in
the same $B$-orbit of $\miso$ are in the same $\liftB$-orbit: an element
of $B$ that moves one to the other considered as an element of $\liftB$
moves also the preimage of the one to that of the other.  

On the other hand,  the preimages 
of distinct $T$-fixed points are distinct $\liftT$-fixed points,
the corresponding map from $\isubn'$ to $I_{\liftn}$ being given as follows:
\[ \underline{i}=
\{i_1<\ldots<i_d\}\mapsto\left\{
\begin{array}{ll}
\{\lifti_1,\ldots,\lifti_d,d+1\} & \textup{if $\underline{i}\in\isubn$}\\
\{\lifti_1,\ldots,\lifti_d,d+2\} &
                   \textup{if $\underline{i}\in\isubn'\setminus\isubn$}
\end{array}\right.   \]
where 
\[
\lifti_k=\left\{\begin{array}{ll}
i_k & \textup{if $i_k\leq d$}\\
i_k+1 & \textup{if $i_k\geq d+2$}
\end{array}\right.         \]
(Note that $d+1$ never occurs as an entry in any element of $\isubn'$
and that the elements $\lifti_1$, \ldots, $\lifti_d$, $d+1$
(respectively $\lifti_1$, \ldots, $\lifti_d$, $d+2$) are not in increasing
order except in the trivial case~$\underline{i}=\{1<\ldots<d\}$.)
Given that each $B$-orbit has a $T$-fixed point and that distinct 
$\liftT$-fixed points belong to distinct $\liftB$-orbits,  this implies that
the preimages of two elements in distinct $B$-orbits belong to distinct
$\liftB$-orbits, and the proof is over.~\hfill$\Box$

\mysection{The theorem}\mylabel{s.theorem}
\checked{\S\ref{s.theorem}, 21 Nov, 1100hrs}
\checked{Added definition of degree; 18 Jan 2007}
\spellchecked{\S\ref{s.theorem}, 21 Nov, 1100hrs}
The purpose of this section is to state the main theorem and its corollary.
We first set down some basic notation and two fundamental 
definitions needed in order to
state the theorem.
\subsection{Basic notation}\mylabel{ss.basicnotn}
We keep the terminology and notations of \S\ref{ss.problem},~\ref{ss.choices}.
As observed in~\S\ref{ss.neven},  it is enough to consider the case $n$
\index{n@$n:=\dim V$, (even from \S\ref{ss.basicnotn} on)} even.
So from now on let $n=2d$\index{d@$d$, integral part of $n/2$ (unfortunately also used otherwise)}.
Recall that, for an integer $k$, $1\leq k\leq 2d$, $k^*:=2d+1-k$\index{kstar@$k^*(:=n+1-k)$}.
As observed in~\S\ref{ss.choices}, Schubert
varieties in $\miso$ are indexed by $\isubn$.   

Since $d$ now determines $n$,  we will henceforth write $\id$\index{I(d)@$\id$}
instead of $\isubn$.    In other words, $\id$ is the set of $d$-element of 
subsets of $\{1,\ldots,2d\}$ such that
\begin{itemize}
\item for each $k$, $1\leq k\leq 2d$, the subset contains exactly one of $k$, $k^*$,
and
\item the number of elements in the subset that exceed~$d$ is even.
\end{itemize}
We write $I(d,2d)$\index{I(d,2d)@$I(d,2d)$} for the set of all $d$-element subsets of $\{1,\ldots,2d\}$.
There is a natural partial order on $I(d,2d)$ and so
also on $\id$:  $v=(v_1<\ldots<v_d)\leq%
\index{\leq@$\leq$, partial order on $I(d,2d)$} %
w=(w_1<\ldots<w_d)$ if and only
if $v_1\leq w_1$, \ldots, $v_d\leq w_d$.

Given $v\in\id$,  the corresponding $T$-fixed point in $\miso$ (namely,
the span of $e_{v_1}$, \ldots, $e_{v_d}$) is denoted $\pointe^v$%
\index{ev@$\pointe^v$, $T$-fixed point}.
Given $w\in\id$,   the corresponding Schubert variety
in $\miso$ (which, by definition, is the closure of the $B$-orbit of
the $T$-fixed point $\pointe^w$ with canonical reduced scheme structure)
is denoted $X(w)$\index{X(w)@$X(w)$, Schubert variety}.   
The point $\pointe^v$ belongs to $X(w)$ if and only if $v\leq w$ in the
partial order just defined.  Since, under the natural action of $B$ on $X(w)$,
each point of $X(w)$ is in the $B$-orbit of a $T$-fixed point $\pointe^v$ for 
some $v$ such that $v\leq w$,  it is enough to focus attention on such
$T$-fixed points.

For the rest of this section an element $v$\index{v@$v$, fixed element of $I(d)$} of $\id$ will remain fixed.

We will be dealing extensively with ordered pairs $(r,c)$,
$1\leq r,c\leq 2d$,  such that $r$ is not and $c$ is an entry of~$v$.
Let $\androotsv$\index{rootsv@$\androotsv$} denote the set of all such ordered pairs, and set
\\[1mm]
\begin{minipage}{6cm}
\begin{align*}
  \andposv &:= \index{N@$\protect\andposv$}%
  \left\{(r,c)\in\androotsv\st r>c\right\}\\
  \rootsv &:= \index{orootsv@$\protect\rootsv$}%
  \left\{(r,c)\in\androotsv\st r<c^*\right\}\\
  \posv &:= \index{on@$\protect\posv$}%
  \left\{(r,c)\in\androotsv\st r>c, r<c^*\right\}\\
  &=\rootsv\cap\andposv\\
  \diagv &:= \index{diagonal, $\protect\diagv$}%
  \left\{(r,c)\in\androotsv\st r=c^*\right\}\\
\end{align*}\vfill
\end{minipage}
\hfill \begin{minipage}{6cm}
\setlength{\unitlength}{.34cm}
\begin{picture}(12,12)(-3,0)
\multiput(0,0)(12,0){2}{\line(0,1){12}}
\multiput(0,0)(0,12){2}{\line(1,0){12}}
\linethickness{.05mm}
\multiput(0,0)(1,0){13}{\line(0,1){12}}
\multiput(0,0)(0,1){13}{\line(1,0){12}}
\thicklines
\put(0,0){\line(1,1){12}}
\put(3.25,2.75){diagonal}
\thicklines
\put(5.25,10.25){boundary}\put(5.25,9.25){of $\andposv$}
\put(0,12){\line(1,0){3}}
\put(3,12){\line(0,-1){1}}
\put(3,11){\line(1,0){2}}
\put(5,11){\line(0,-1){2}}
\put(5,9){\line(1,0){2}}
\put(7,9){\line(0,-1){1}}
\put(7,8){\line(1,0){1}}
\put(8,8){\line(0,-1){1}}
\put(8,7){\line(1,0){1}}
\put(9,7){\line(0,-1){2}}
\put(9,5){\line(1,0){2}}
\put(11,5){\line(0,-1){2}}
\put(11,3){\line(1,0){1}}
\put(12,3){\line(0,-1){3}}
\put(1,5){\circle{.4}}
\put(.25,5.5){$(r,c)$}
\put(1.25,.75){$(c^*,c)$}
\put(5.25,4.75){$(r,r^*)$}
\put(1.25,3){leg}
\put(2.75,5.25){leg}
\multiput(1,5)(0,-.5){8}{\line(0,-1){.3}}
\multiput(1,5)(.5,0){8}{\line(1,0){.3}}
\put(0,1){\circle*{.4}}
\put(1,7){\circle*{.4}}
\put(6,8){\circle*{.4}}
\put(7,9){\circle*{.4}}
\end{picture}
\end{minipage}
\\

The picture shows a drawing of $\androotsv$.   We think of $r$ and $c$
in $(r,c)$ as row index and column index respectively.    The columns
are indexed from left to right 
by the entries of~$v$ in ascending order,  the rows from top to bottom by
the entries of $\{1,\ldots,2d\}\setminus v$ in ascending order.
The points of $\diagv$ are those on the diagonal, 
the points of $\rootsv$ are those that are (strictly) above the diagonal, and
the points of $\andposv$ are those that are to the South-West of the
poly-line captioned ``boundary of $\andposv$''---we draw the
boundary so that points on the boundary belong to $\andposv$.
The reader can readily verify that  $d=13$ and $v=(1,2,3,4,6,7,10,11,13,15,18,19,22)$ for the particular picture drawn.   The points of $\posv$ indicated
by solid circles form a $v$-chain (see~\S\ref{sss.vchain} below).
 
We will be considering {\em monomials}, also called 
{\em multisets\index{multiset@$\textup{multiset}:=\textup{monomial}$}}, in some of these sets.
A {\em monomial}\index{monomial},  as usual,  is a subset with each member being allowed
a multiplicity (taking values in the non-negative integers).  The 
{\em degree\/}\index{degree, of a monomial}
of a monomial has also the usual sense: 
it is the sum of the multiplicities in the
monomial over all elements of the set.  The {\em intersection\/}%
\index{intersection (of a monomial in a set with a subset)}
of a monomial in a set with a subset of the set has also the natural meaning:
it is a monomial in the subset,  the multiplicities being those in the
original monomial.

We will refer to $\diagv\index{diagonal, $\diagv $}$ as the 
{\em diagonal}\index{diagonal, $\diagv $}.
\mysubsection{Two fundamental definitions}\mylabel{ss.2funddef}%
\mylabel{ss.vchaindom}
\mysubsubsection{Definition of\/ $v$-chain}\mylabel{sss.vchain}\index{v-chain@$v$-chain}
Given two elements $(R,C)$ and $(r,c)$ in $\posv$,  we write
$(R,C)>(r,c)$\index{>@$>$, relation on $\posv$}
if $R>r$ and $C<c$ (note that these are strict
inequalities).
An ordered sequence $\alpha$, $\beta$, \ldots\ 
of elements of $\posv$ is called a {\em $v$-chain\/} if
$\alpha>\beta>\ldots\ $.  
A $v$-chain $\alpha_1>\ldots>\alpha_\ell$
has {\em head%
\index{head, of a $v$-chain}%
\/}~$\alpha_1$, {\em tail%
\index{tail, of a $v$-chain}%
\/}~$\alpha_\ell$,
and {\em length%
\index{length, of a $v$-chain}%
\/}~$\ell$.

\mysubsubsection{Definition of \ortho-domination}\mylabel{sss.domination}
To a $v$-chain $C:\alpha_1>\alpha_2>\ldots$ in $\posv$ there corresponds,
as described in~\S\ref{sss.montoC}, 
a subset $\mon_C$ of $\andposv$ which, as observed
in Proposition~\ref{p.montoC}, is ``distinguished'' in the
sense of~\S\ref{sss.ddist}.  
To a distinguished subset of~$\andposv$  there corresponds,
as described below in~\S\ref{sss.montow}, 
an element of~$\idd$.    Following these correspondences through,
we get an element of~$\idd$ attached to the $v$-chain~$C$. 
Let~$w(C)$\index{w(C)@$w(C)$ (or $w_C$), where $C$ is $v$-chain} 
denote this element---sometimes we write $w_C$.  (All this makes sense even
when $C$ is empty---$w(C)$ will turn out to be~$v$ itself in that case.)

Furthermore, as will be obvious from its definition,
the monomial~$\mon_C$ is ``symmetric'' in the sense
of~\S\ref{sss.hashonposv} and contains evenly many elements of
the diagonal~$\diagv$.    Thus, by Proposition~\ref{p.distid},
the element $w(C)$ of $\idd$ belongs to~$\id$.   

An element $w$ of $\id$
is said to {\em \ortho-dominate\/}\index{ortho-domination@$\ortho$-domination} $C$ if $w\geq w(C)$, or, equivalently---and this is important for the proofs---if
$w$ dominates in the sense of~\cite{kr} the monomial~$\mon_C$ 
(for the proof of the equivalence, see~\cite[Lemma~5.5]{kr}).
An element~$w$ of $\id$ 
{\em \ortho-dominates\/}\index{ortho-domination@$\ortho$-domination} 
a monomial $\mon$ 
of $\posv$ (repsectively of $\rootsv$)
if it \ortho-dominates
every $v$-chain in $\mon$ (respectively in $\mon\cap\posv$).
\subsection{The main theorem and its corollary}\mylabel{ss.mainthm}
\begin{theorem}\mylabel{t.main}
Fix a positive integer $d$ and elements $v\leq w$ of $\id$.   Let
$V$ be a vector space of dimension $2d$ with a symmetric non-degenerate
bilinear form (over a field of characteristic not~$2$).
Let $X(w)$ be the Schubert variety corresponding to~$w$ in
the orthogonal Grassmannian $\misd$, and~$\pointe^v$ the torus fixed
point of~$X(w)$ corresponding to~$v$.  
Let $R_v^w$ denote the associated graded ring
with respect to the unique maximal ideal of the local ring of germs
at~$\pointe^v$ of functions on $X(w)$.   Then, for any non-negative integer $m$,
the dimension as a vector space of the homogeneous piece of $R_v^w$
of degree $m$ equals the cardinality of the set $S^w(v)(m)$%
\index{Svwm@$S^w(v)(m)$}  
of monomials of degree $m$ of\/ $\rootsv$ that are \ortho-dominated by~$w$.
\end{theorem}
\noindent
The proof of this theorem occupies us for most of this paper.
It is reduced in~\S\ref{s.reduction},
by an application of standard monomial theory, 
to combinatorics.    The resulting combinatorial problem
is solved in~\S\ref{s.furtherred}--\ref{s.proof}.
For now, let us note the following immediate consequence:
\begin{corollary}\mylabel{c.t.main}\mylabel{c.main}
The multiplicity at the point~$\pointe^v$ of the Schubert variety $X(w)$ equals
the number of monomials in $\posv$ of maximal cardinality that are
square-free and \ortho-dominated by~$w$.
\end{corollary}
\begin{myproof} The proof of Corollary~2.2 of \cite{kr} holds verbatim here too.
\end{myproof}

\mypart{From geometry to combinatorics}\mylabel{p.gtoc}
The problem is translated from geometry to combinatorics.
The main combinatorial results are formulated.
\mysection{Reduction to combinatorics}\mylabel{s.reduction}
\checked{\S\ref{s.reduction}, 21 Nov 1230hrs}
\spellchecked{\S\ref{s.reduction}, 21 Nov 1230hrs}
In this section we translate the problem from geometry to combinatorics.
In \S\ref{ss.smt} we recall from~\cite{s} the theorem that enables the
translation.    The translation itself is done in~\ref{ss.translate}
and follows~\cite{lw}.
\mysubsection{Homogeneous co-ordinate ring of the Schubert variety~$X(w)$}%
\mylabel{ss.smt}
\mysubsubsection{The line bundle~$L$ on~$\miso$}\mylabel{sss.linebundle}
Let $\miso\subseteq G_{d}(V)\hookrightarrow \mathbb{P}(\wedge^d
V)$ be the Pl\"ucker embedding (where $G_d(V)$ denotes the
Grassmannian of all $d$-dimensional subspaces of $V$).
The pull-back to $\miso$ of the line bundle
$\mathcal{O}(1)$ on $\mathbb{P}(\wedge^d V)$ is the square of the
ample generator of the Picard group of $\miso$.
Letting $L$\index{L@$L$, line bundle} denote the ample generator, we observe that it is very ample and
want to describe the homogeneous coordinate rings of~$\miso$
and its Schubert subvarieties in the embedding defined by $L$.
\mysubsubsection{The section~$q_\theta$ of $L$}\mylabel{sss.qtheta}
\newcommand{\grass}{G_d(V)}
\newcommand{\matb}{M}
For $\theta$ in $I(d,2d)$,  let $p_\theta$\index{ptheta@$p_\theta$, Pl\"ucker coordinate} denote the corresponding Pl\"ucker
coordinate. Consider the affine patch $\mathbb{A}$ of $\mathbb{P}(\wedge^d V)$
given by $p_{\epsilon}=1$, where $\epsilon:=(1,\ldots,d)$.
The intersection $\mathbb{A}\cap\grass$ of
this patch with the Grassmannian is an affine space.     Indeed the $d$-plane
corresponding to an arbitrary point $z$ of $\mathbb{A}\cap\grass$ has a basis consisting
of column vectors of a matrix of the form
\[ \matb= \left(\begin{array}{c}I\\A\\ \end{array}\right)\]
where $I$ is the identity matrix and $A$
an arbitrary matrix both of size $d\times d$.
The association $z\mapsto A$ is bijective.
The restriction
of a Pl\"ucker coordinate $p_\theta$ to $\mathbb{A}\cap\grass$ is given by
the determinant of a submatrix of size $d\times d$ of $\matb$,  the entries of
$\theta$ determining the rows to be chosen from $\matb$ to form the submatrix.

As can be readily verified, a point $z$ of $\mathbb{A}\cap\grass$ represents
an isotropic subspace if and only if the corresponding matrix 
$A=(a_{ij})$ is {\em skew-symmetric
with respect to the anti-diagonal\/}:   $a_{ij}+a_{j^*i^*}=0$,  where the
columns and rows of $A$ are numbered $1,\ldots,d$ and $d+1,\ldots,2d$
respectively.
For example, if $d=4$, then a matrix that is skew-symmetric with respect to the
anti-diagonal looks like this:
\[
\left(\begin{array}{rrrr}
-d & -c & -b & 0\\
-g & -f & 0 & b\\
-i & 0 & f & c\\
0 & i & g & d\\ \end{array}\right)  \]
Since the set of these matrices is connected and contains the point that
is spanned by $e_1,\ldots,e_d$,   it follows that $\mathbb{A}\cap\grass$
does not intersect the other component of $\miso'$. In other words, $p_\epsilon$
vanishes everywhere on $\miso'\setminus\miso$.

Now suppose that $\theta$ belongs to $\id$.
Computing $p_\theta/p_\epsilon$ as a function on the affine patch $p_\epsilon\neq0$,
we see that it is the determinant of a skew-symmetric matrix of even size, and
therefore a square.    The square root, which is determined up to sign, is
called the {\em Pfaffian}.    
This suggests that $p_\theta$ itself is a square:  more
precisely that there exists a section $q_\theta$\index{qtheta@$q_\theta$, Pfaffian} of the line bundle $L$ on
$\miso$ such that $q_\theta^2=p_\theta$.     A weight calculation confirms this
to be the case.   The $q_\theta$ are also called {\em Pfaffians\index{Pfaffian $q_\theta$}}.
\mysubsubsection{Standard monomial theory for $\miso$}\mylabel{sss.smt}
A {\em standard monomial\/}\index{standard monomial}
in $\id$ is a totally ordered sequence $\theta_1\geq\ldots\geq \theta_t$
(with repetitions allowed) of elements of $\id$.    Such a
standard monomial is said to be {\em $w$-dominated\/}\index{standard monomial!w-dominated@$w$-dominated} for
$w\in\id$ if $w\geq \theta_1$. To a standard monomial
$\theta_1\geq\ldots\geq \theta_t$ in $\id$ we associate the product
$q_{\theta_1}\cdots q_{\theta_t}$,  where the $q_\theta$ are the sections
defined above of the line bundle~$L$. Such a product is also called a {\em
standard monomial\/}\index{standard monomial} and it is said to 
be {\em dominated by $w$\/}\index{standard monomial!w-dominated@$w$-dominated}
for $w\in\id$ if the underlying monomial in $\id$ is dominated
by $w$.  {\em Standard monomial theory\/} for $\miso$  says:
\begin{theorem}\mylabel{t.smt}
{\rm (Seshadri~\cite{s})} 
Standard mo\-nomials $q_{\theta_1}\cdots q_{\theta_r}$ of degree $r$ 
form a basis for the space of forms of
degree $r$ in the homogeneous coordinate ring of\/ $\miso$ in the
embedding defined by the ample generator $L$ of the Picard group.
More generally,   for $w \in \id$,  the $w$-dominated standard
monomials of degree $r$ form a basis for the space of forms of
degree $r$ in the homogeneous coordinate ring of the Schubert
subvariety $X(w)$ of\/ $\miso$.
\end{theorem}
\mysubsection{Co-ordinate rings of affine patches and tangent cones of $X(w)$}\mylabel{ss.translate}
{From}~Theorem~\ref{t.smt} one can deduce rather easily, as we now show,
bases for co-ordinate rings of affine patches of the form $q_v\neq0$
and of tangent cones of Schubert varieties.
An element $v$ of $\id$ will remain
fixed for the rest of this section.   To simplify notation we will
suppress explicit reference to~$v$.

\mysubsubsection{Standard monomial theory for affine patches}%
\mylabel{sss.smtforpatch}
Let $\affinev$\index{A@$\affinev$, affine patch $q_v\neq 0$ of $\miso$}
denote the affine patch of $\mathbb{P}(H^0(\miso,L)^*)$ 
given by $q_v\neq 0$.   The origin of the affine space~$\affinev$ is
identified as the $T$-fixed point $\pointe^v$.
The functions $f_\theta:=q_\theta/q_v$\index{ftheta@$f_\theta:=q_\theta/q_v$}, $v\neq\theta\in\id$, provide a 
set of coordinate functions on~$\affinev$.   Monomials in these $f_\theta$
form a $\field$-basis for the polynomial ring $\field[\affinev]$ of
functions on~$\affinev$,  where~$\field$%
\index{k@$\field$, base field, ($\textup{characteristic}\neq2$)}
denotes the underlying field.

Fix $w\geq v$ in $\id$, so that the point $\pointe^v$ belongs to 
the Schubert variety $X(w)$, and 
let~$\yvw$\index{Y(w)@$Y(w)(:=X(w)\cap\affinev)$}
be the affine patch of $X(w)$ defined thus:
\[\yvw:=X(w)\cap\affinev.\]
The coordinate ring $\field[\yvw]$
of $\yvw$ is a quotient of the polynomial ring $\field[\affinev]$, and
the proposition that follows identifies a subset of the monomials
in $f_\theta$ which forms a $\field$-basis for~$\field[\yvw]$.

We say that a standard monomial $\theta_1\geq\ldots\geq\theta_t$
in $\id$ is {\em $v$-compatible\/}\index{standard monomial!v-compatible@$v$-compatible} if for each $k$, $1\leq k\leq t$,
either $\theta_k\gneq v$ or $v\gneq\theta_k$.   
Given $w$ in $I(d)$, we denote by $\smvw$\index{SMvw@$\smvw$} the set
of $w$-dominated $v$-compatible standard monomials. 

\begin{proposition}\mylabel{p.smtforpatch}
As $\theta_1\geq\ldots\geq\theta_t$ runs over the set $\smvw$ of
$w$-dominated $v$-compatible standard monomials, the elements
$f_{\theta_1}\cdots f_{\theta_t}$  form a basis for the coordinate
ring $\field[\yvw]$ of the affine patch $\yvw=X(w)\cap\affinev$ of the
Schubert variety~$X(w)$.
\end{proposition}
\begin{myproof}
The proof is similar to the proof of Proposition~3.1 of \cite{kr}.
First consider a linear dependence relation among the
$f_{\theta_1}\cdots f_{\theta_t}$.  Replacing $f_\theta$ by $q_\theta$ 
and ``homogenizing'' by $q_v$ yields a linear dependence relation among the
$w$-dominated standard monomials $q_{\theta_1}\cdots q_{\theta_s}$ 
restricted to $X(w)$, and so the original relation must only have been 
the trivial one, for by Theorem~\ref{t.smt} the $q_{\theta_1}\cdots q_{\theta_s}$
are linearly independent on~$X(w)$.

To prove that $f_{\theta_{1}}\cdots f_{\theta_{t}}$ generate $\kyvw$ as 
a vector space, 
we make the following claim: if $q_{\mu_{1}}\cdots q_{\mu_{r}}$ be any monomial
in the Pfaffians $q_\theta$, and $q_{\tau_{1}}\cdots q_{\tau_{s}}$ 
a standard monomial that occurs 
with non-zero co-efficient in the expression for (the restriction 
to~$X(w)$ of) $q_{\mu_{1}}\cdots q_{\mu_{r}}$ 
as a linear combination of $w$-dominated standard monomials, then 
$\tau_{1}\cup\cdots\cup\tau_{s}=\mu_{1}\cup\cdots\cup\mu_{r}$ 
as multisets of $\{1,\ldots,2d\}$.   To prove the claim, consider the
maximal torus $T$ of $\sov$ as in~\S\ref{ss.choices}.  The affine
patch~$\affinev$ is $T$-stable and there is an action of~$T$ on~$\kyvw$.
The sections~$q_\theta$ are eigenvectors for $T$ with corresponding
characters $\epsilon_{\theta_1}+\cdots+\epsilon_{\theta_d}$, where $\epsilon_k$
denotes the character of $T$ given by the projection to the diagonal entry
on row~$k$.    The claim now follows since eigenvectors corresponding to
different characters are linearly independent.

Let $f_{\mu_1}\cdots f_{\mu_r}$ be an arbitrary monomial in the $f_\theta$.
Fix an integer $h$ such that $h>r(d-1)$ and 
consider the expression for (the restriction to $X(w)$ of)
$q_{\mu_{1}}\cdots q_{\mu_{r}}\cdot q_{v}^{h}$
as a linear combination of $w$-dominated standard monomials.   
We claim that $q_{v}$ occurs in every standard monomial
$q_{\tau_{1}}\cdots q_{\tau_{r+h}}$ in this expression (from which it
will follow that the~$\tau_j$ are all comparable to~$v$).
Suppose that none of $\tau_{1},\ldots,\tau_{r+h}$ equals $v$. For each 
$\tau_{j}$ there is at least one entry of $v$ that does not occur in it.
The number of occurrences of entries of $v$ in $\tau_{1}\cup\cdots\cup\tau_{r+h}$
is thus at most $(r+h)(d-1)$.   But these entries occur at least $hd$ times 
in $\mu_{1}\cup\cdots\cup\mu_{r}\cup v\cup\cdots\cup v$ 
(where $v$ is repeated $h$ times), a contradiction to the claim proved in
the previous paragraph.  Hence our claim is proved.
Dividing by $q_{v}^{r+h}$ the expression for 
$q_{\mu_{1}}\cdots q_{\mu_{r}}.q_{v}^{h}$ as a linear combination of 
$w$-dominated standard monomials provides an expression for 
$f_{\mu_{1}}\cdots f_{\mu_{r}}$ as a linear combination of 
$f_{\theta_{1}}\cdots f_{\theta_{t}}$, as 
$\theta_{1}\geq\ldots\geq\theta_{t}$ varies over $\smvw$.  
\end{myproof}
\mysubsubsection{Standard monomial theory for tangent cones}%
\mylabel{sss.smttgtcone}
The affine patch $\miso\cap\affinev$ of the orthogonal Grassmannian
$\miso$ is an affine space 
whose coordinate ring can be taken to be the polynomial ring in 
variables of the form $X_{(r,c)}$\index{Xrc@$X_{r,c}$, variable}
with $(r,c)\in\rootsv$, where 
(as in~\S\ref{ss.basicnotn})
\[ \rootsv=\{(r,c)\st 1\leq r,c \leq 2d, r\not\in v, c\in v, r<c^*\}\]
Taking $d=5$ and $v=(1,3,4,6,9)$ for example,   a general element
of $\miso\cap\affinev$ has a basis consisting of column vectors
of a matrix of the following form:
\[\left(\begin{array}{ccccc}
1 & 0 & 0 & 0 & 0 \\
X_{21} & X_{23} & X_{24} & X_{26} & 0 \\
0 & 1 & 0 & 0 & 0 \\
0 & 0 & 1 & 0 & 0 \\
X_{51} & X_{53} & X_{54} & 0 & -X_{26}\\
0 & 0 & 0 & 1 & 0\\
X_{71} & X_{73} & 0 & -X_{54} & -X_{24}\\
X_{81} & 0 & -X_{73} & -X_{53} & -X_{23}\\
0 & 0 & 0 & 0 & 1\\
0 & -X_{81} & -X_{71} & -X_{51} & -X_{21}\\
\end{array}\right)\]
The expression for $f_\theta=q_\theta/q_v$ in terms of the $X_{(r,c)}$
is a square root of the determinant of the submatrix of a matrix 
like the one above obtained by choosing the rows 
given by the entries of $\theta$.
Thus $f_\theta$ is a homogeneous polynomial of degree the 
$v$-degree of $\theta$, where the {\em $v$-degree}\index{v-degree@$v$-degree}
of $\theta$ is defined as
one half of the cardinality of $v\setminus\theta$.

Since the ideal of the Schubert variety $X(w)$ in the homogeneous  
coordinate ring of $\miso$ is generated\footnote{
This is a consequence of Theorem~\ref{t.smt}.  
It is easy to see that the $q_\tau$ 
such that $\tau\not\leq w$  vanish on $X(w)$.    
Since all standard monomials form a basis for the homogeneous  
coordinate ring of $\miso$ in $\mathbb{P}(H^0(\miso,L)^*)$, it follows that
$w$-dominated standard monomials span
the quotient ring by the ideal generated by such $q_\tau$.
Since such monomials are linearly independent in the homogeneous
 coordinate ring of $X(w)$,    the desired result follows.} 
by the $q_\tau$, $\tau\in\id$ such that $\tau\not\leq w$,
it follows that the ideal of $\yvw:=X(w)\cap\affinev$ in 
$\miso\cap\affinev$ is generated by the the $f_\tau$, $\tau\in\id$
such that $\tau\not\leq w$. 
We are interested in the tangent cone to $X(w)$ at~$\pointe^v$ 
(or, what is the same, the tangent cone
to $\yvw$ at the origin),   and since $\kyvw$ is graded,
its associated graded ring with respect to
the maximal ideal corresponding to the origin
is $\kyvw$ itself.   

Proposition~\ref{p.smtforpatch} says that the graded piece 
of $\kyvw$ of degree~$m$
is generated as a $\field$-vector space by elements of degree $m$ of the
set $\smvw$ of $w$-dominated $v$-compatible standard monomials,
where the {\em degree\/}\index{degree, of a standard monomial} of a standard monomial 
$\theta_1\geq\ldots\geq\theta_t$ is
defined to be the sum of the $v$-degrees of $\theta_1,\ldots,\theta_t$.   
To prove Theorem~\ref{t.main} it therefore suffices to prove the following:
\begin{theorem}\mylabel{t.main.1} 
The set $\smvwm$ of standard monomials in~$\id$ of degree~$m$ that are
$w$-dominated and $v$-compatible is in bijection with
the set $S^w(v)(m)$ of monomials in $\oroots$ of degree~$m$ 
that are $\ortho$-dominated by~$w$.
\end{theorem}

\mysection{Further reductions}\mylabel{s.furtherred}
In the last section, we reduced the proof of our main theorem 
(Theorem~\ref{t.main}) to that of Theorem~\ref{t.main.1}.     
We now reduce the proof of Theorem~\ref{t.main.1} to that of
Propositions~\ref{p.4.1.kr},~\ref{p.4.2.kr} and~\ref{p.new} below.
These propositions will eventually be proved in~\S\ref{s.proof}.     

\mysubsection{The main propositions}\mylabel{ss.p.main}
Fix once and for all an element~$v$ of $I(d)$.   The bijection
stated in Theorem~\ref{t.main.1} will be described by means of
two maps~$\opi$\index{opi@$\opi$} and $\ophi$\index{ophi@$\ophi$} whose definitions will be given
in \S\ref{s.opi} and \S\ref{s.ophi} below.  We will now state
some properties of these maps. In~\S\ref{ss.mpropstomt} we will see how
Theorem~\ref{t.main.1} follows once these properties are established.

The map $\opi$ associates to a monomial~$\mon$ in~$\posv$ a pair
$(w,\mon')$   consisting of an element $w$ of $\id$ and a ``smaller''
monomial $\mon'$ in $\posv$.   This map enjoys the following good
properties:
\begin{proposition}\mylabel{p.4.1.kr}
\begin{enumerate}
\item\mylabel{i.1.p.4.1.kr}
$w\geq v$.
\item\mylabel{i.2.p.4.1.kr}
$\vdeg{w}+\degree{\mon'}=\degree{\mon}$.
\item\mylabel{i.3.p.4.1.kr}
$w$ \ortho-dominates $\mon'$.
\item\mylabel{i.4.p.4.1.kr}
$w$ is the least element of $\id$ that \ortho-dominates $\mon$.
\end{enumerate}
\end{proposition}
The map $\ophi$, on the other hand,   associates a monomial in
$\posv$ to a pair $(w,\mont)$ consisting of an element $w$ of $\id$ with
$w\geq v$  and a monomial $\mont$ in $\posv$ that is $\ortho$-dominated by~$w$.
\begin{proposition}\label{p.4.2.kr}
The maps $\opi$ and $\ophi$ are inverses of each other.
\end{proposition}

For an integer $f$, $1\leq f\leq 2d$,   consider the following conditions, the
first on a monomial $\mon$ in $\posv$, the second on an element $w$ of~$\id$:
\begin{quote}
(\ddag)~$f$ is not the row index of any element of $\mon$ and
$f^\star$ is not the column index of any element of $\mon$.\\
  (\ddag) $f$ is not an entry of $w$.
\end{quote}
(It is convenient to the use the same notation~(\ddag) for both conditions.)
\begin{proposition}\mylabel{p.new}\mylabel{p.4.3.kr}
Assume that $v$ satisfies (\ddag)---all references to (\ddag) in this
proposition are with respect to a fixed $f$, $1\leq f\leq 2d$.
  \begin{enumerate} 
  \item Let $w$ be an element of $\id$ with $w\geq v$ and $\mont$ a monomial
in $\posv$ that is \ortho-dominated by~$w$.   If $w$ and $\mont$ both
satisfy (\ddag),  then so does $\ophi(w,\mont)$.
\item If a monomial $\mon$ in $\posv$ satisfies (\ddag),  then so do the
``components'' $w$ and $\mon'$ of its image under $\opi$.
  \end{enumerate}
\end{proposition}
\mysubsection{From the main propositions to the main theorem}%
\mylabel{ss.mpropstomt}
Let us now see how Theorem~\ref{t.main.1} follows from the propositions
of~\S\ref{ss.p.main}.
Most of the following argument runs parallel to its counterparts in the
case of the Grassmannian and symplectic Grassmannian
(Propositions~\ref{p.4.1.kr} and~\ref{p.4.2.kr} have their counterparts
in~\cite{kr,gr}), but, in
the case that $d$ is odd, the part involving the ``mirror image''
requires additional work.   This is where Proposition~\ref{p.new} comes in.

Let $\sv$%
\index{S@$\sv$, set of monomials in $\rootsv$}%
\index{S@$\sv$, set of monomials in $\rootsv$!modifications|see{Notation~\ref{n.explain}}},
$\tv$%
\index{T@$\tv$, set of monomials in $\posv$}%
\index{T@$\tv$, set of monomials in $\posv$!modifications|see{Notation~\ref{n.explain}}},
and $\uv$%
\index{U@$\uv$, set of monomials in $\rootsv\setminus\posv$}%
\index{U@$\uv$, set of monomials in $\rootsv\setminus\posv$!modifications|see{Notation~\ref{n.explain}}}, denote 
respectively the sets of monomials 
in $\rootsv$, $\posv$, and $\rootsv\setminus\posv$.    
Let $\smvadv$\index{smsubv@$\smvadv$} 
denote the set of $v$-compatible standard monomials that are ``anti-dominated''
by~$v$:   a standard monomial $\theta_1\geq \ldots\geq \theta_t$ is
{\em anti-dominated}\index{anti-domination} 
by $v$ if $\theta_t\geq v$   (we can also write
$\theta_t>v$ since $\theta_t\neq v$ by $v$-compatibility).

Define the {\em domination map}\index{domination map} from $\tv$ to $\id$ by
sending a monomial in $\posv$ to the least element that \ortho-dominates it.
Define the {\em domination map}\index{domination map} from $\smvadv$
to $\id$ by sending $\theta_1\geq\ldots\geq\theta_t$ to
$\theta_1$. 
Both these maps take, by definition, the value~$v$ on the empty monomial.

\bnotation\mylabel{n.explain}
In the following, we use subscripts, superscripts, suffixes, and combinations
thereof to modify the meanings of $\sv$, $\tv$, $\uv$, $\smv$\index{SMmodifications@$\smv$}%
\index{SMmodifications@$\smv$!modifications|see{Notation~\ref{n.explain}}}, and $\smvadv$.
\begin{itemize}
\item superscript: this will be an element $w$ of $\id$; when used on
$T$ it denotes $\ortho$-domination (more precisely, 
$T^w$ denotes the subset of $T$ consisting of those elements that are
$\ortho$-dominated by~$w$);  when used on $\smv$ or $\smvadv$
it denotes domination by $w$.
\item subscript: denotes anti-domination (applied only to 
standard monomials).
\item suffix ``$(m)$'':  indicates degree (for example, $\smvadvwm$ denotes
the set of $v$-compatible standard monomials that are anti-dominated by $v$,
dominated by $w$, and of degree $m$).
\end{itemize}
\enotation

Repeated application of $\opi$ gives a map from $\tv$ to $\smvadv$
that commutes with domination (as just defined) and preserves degree.   
Repeated application
of $\ophi$ gives a map from $\smvadv$ to $\tv$.   These two maps being
inverses of each other (Proposition~\ref{p.4.2.kr}) 
and so we have a bijection between
$\smvadv$ and $\tv$.   In fact, since domination and degree are
respected (Proposition~\ref{p.4.1.kr}),   
we get a bijection $\smvadvwm\cong\tvwm$. 

As explained below,  the ``mirror image'' of the bijection
$\smvadv\modm\cong\tvm$ gives a bijection $\smvvm\cong\uvm$.
Putting these bijections together, we get the desired result:
\begin{eqnarray*}
\smvwm &=& 
	\bigcup\limits_{k=0}^m \sm^w_v\modv(k) \times
	\sm^v\modv(m-k)\\
&\cong&\bigcup\limits_{k=0}^m T^w\modv(k) \times
				\uv(m-k)
				= \svwm.
\end{eqnarray*}

We now explain how to realize the bijection $\smvvm\cong\uvm$ as
the ``mirror image'' of the bijection $\smvadv\modm\cong\tvm$.
For an element $u$ of $\id$, define $u^*:=(u_d^*,\ldots,u_1^*)$%
\index{ustar@$u^*$, for $u\in\id$}.
In the case $d$ is even,  the association $u\mapsto u^*$ is an
order reversing involution, and the argument in~\cite{gr} for
the symplectic Grassmannian holds here too.
In the case $d$ is odd, $u^*$ is not an element of $\id$, and so
some additional work is required.

Recall that a ``base element''~$v$ of~$\id$ has been 
fixed and that our notation does not
explicitly indicate this dependence upon~$v$:
for example,  $\oroots$ is dependent upon~$v$.   
For a brief while now (until the end of this section)
we need to simultaneously handle several base elements
of~$\id$.    We will use the following convention:  when the
base element of~$\id$ is not~$v$,  we will explicitly indicate it by means
of a suffix.   For instance, $\sm(v^*)$ denotes the set of
$v^*$-compatible standard monomials in~$\id$.

Let us first do the case when $d$ is even.   We get a bijection
$\smvv\cong\sm_{v^*}\modvstar$ by associating to
$\theta_1\geq\ldots\geq\theta_t$ the element $\theta_t^*
\geq\ldots\geq\theta^*_1$.   The sum of the $v$-degrees of 
$\theta_1,\ldots,\theta_t$ equals the sum of the $v^*$-degrees of
$\theta_t^*,\ldots,\theta_1^*$,  so that we get a bijection
$\smvv\modm\cong\sm_{v^*}\modvstar\modm$.

For an element $(r,c)$ of $\ortho\pos\modvstar$, consider its flip $(c,r)$.
Since $v$ belongs to~$\id$, the complement of $v^*$ in~$\{1,\ldots,2d\}$
is $v$,  and it follows that $(c,r)$ belongs to
$\ortho\roots\modv\setminus\ortho\pos\modv$.    This induces
a degree preserving bijection $T\modvstar\cong U\modv$.
Putting this together with the bijection of the previous paragraph
and the one deduced earlier in this section (using $\opi$ and
$\ophi$),  we get what we want:
\[
\smvv\modm\cong\sm_{v^*}\modvstar\modm\cong T\modvstar\modm\cong \uvm.
								\]

Now suppose that $d$ is odd.  Then the map $x\mapsto x^*$ does not 
map $\id$ to $\id$ but to $\idstar$ (defined as the set consisting
of those elements $u$ of $I(d,2d)$ such that, for each $k$, $1\leq k\leq 2d$,
exactly one of $k$, $k^*$ belongs to $u$, and the number of entries
of $u$ greater than $d$ is odd).   We define a map $u\mapsto\upass$ from
$\idstar$ to $I(d+1)$ as follows: 
  $\upass:=\{\pass{u_1},\ldots,\pass{u_d},d+2\}$
(the elements are not in increasing order except in the trivial case
$u=(1,\ldots,d)$),
where,  for an integer $e$, $1\leq e\leq 2d$, we set 
\[
\pass{e}:=\left\{\begin{array}{ll}
e & \textup{if $1\leq e\leq d$}\\
e+2 & \textup{if $d+1\leq e\leq 2d$}\\
\end{array}
\right. \]
This map $u\mapsto\upass$ is an order preserving injection.

Consider the composition $x\mapsto x^*\mapsto \pass{x^*}$ from
$\id$ to $I(d+1)$.   This is an order reversing injection.
The induced map on standard monomials is an 
injection from $\smvv$ to $\smvspvsp$.  It is readily seen that the image
under this map is the subset $\smvspvspddag$ consisting of those
standard monomials all of whose elements satisfy (\ddag) with
$f=d+1$.     We have already established (using the maps $\opi$
and $\ophi$) a bijection $\smvspvsp\cong T\modvsp$.   It follows
from Proposition~\ref{p.new} that under this bijection the subset
$\smvspvspmodddag$ maps to $\tvspddag$ (defined as the set of those
monomials in $\ortho\pos\modvsp$ satisfying (\ddag) with $f=d+1$).

Now $\tvspddag$ is in degree preserving bijection with
$\uv$:  every element of degree~$1$ of $\tvspddag$ is uniquely of the form
$(\pass{c},\pass{r})$ for $(r,c)$ in $\ortho\roots\modv\setminus\ortho
\pos\modv$,  and the desired bijection is induced from this.
Putting all of these together, we finally have
\[
\smvv\cong\smvspvspmodddag\cong\tvspddag\cong\uv.     \]

Thus, in order to prove our main theorem (Theorem~\ref{t.main}),  it
suffices to describe the maps $\opi$ and $\ophi$ and to prove 
Propositions~\ref{p.4.1.kr}--\ref{p.new}.

\mypart{The proof}\mylabel{p.proof}
The main combinatorial results formulated in~\S\ref{ss.p.main} are proved.   
An attempt is made to maintain parallelism with the proofs in~\cite{kr}.    
\mysection{Terminology and notation}\mylabel{s.notation}
\mysubsection{Distinguished subsets}\mylabel{ss.distinguished}
\mysubsubsection{Distinguished subsets of\/ $\andposv$}%
\mylabel{sss.ddist}
Following~\cite[\S4]{kr}, we define a multiset~$\mon$ of~$\andposv$
to be {\em distinguished\/},%
\index{distinguished (a subset of $\andposv$)}  
if, first of all, it is a subset in
the usual sense (in other words, it is ``multiplicity free''), and if,
for any two distinct elements~$(R,C)$ and~$(r,c)$ of~$\mon$, 
the following conditions are satisfied:
\begin{enumerate}
\item[A.] 
$R\neq r$ and $C\neq c$.
\item[B.] 
If $R>r$,  then either $r<C$ or $C<c$. 
\end{enumerate}
In terms of pictures,  
condition~A says that $(r,c)$ cannot lie exactly due North or East of $(R,C)$
(or the other way around);  so we can assume, interchanging the two points 
if necessary,  that $(r,c)$ lies strictly 
to the Northeast or Northwest of $(R,C)$;
condition~B now says that, if $(r,c)$ lies to the Northwest of $(R,C)$,
then the point that is simultaneously due North of $(R,C)$ and
due East of $(r,c)$ (namely $(r,C)$) does not belong to~$\andposv$.
\mysss{Attaching elements of $\idd$ to distinguished subsets of $\andposv$}%
\mylabel{sss.montow}
To a distinguished subset $\mon$ of $\andposv$ there is naturally associated
an element $w$ of $\idd$ as follows:  start with $v$,  remove all members of
$v$ which appear as column indices of elements of $\mon$,  and 
add row indices of all elements of $\mon$.    As observed in~\cite[%
Proposition~4.3]{kr},  this association gives a bijection between
distinguished subsets of $\andposv$ and elements $w\geq v$ of $\idd$.
The unique distinguished subset of $\andposv$ corresponding to an element
$w\geq v$ of $\idd$ is denoted $\mon_w$.%
\index{Sw@$\monw$, $w$ in $\idd$ or $I(d,n)$}
\mysubsection{The involution~$\hash$}\mylabel{ss.hash}
\mysss{The involution $\hash$ on $\idd$}\mylabel{sss.hashonidd}
There are two natural order reversing involutions on $\idd$.
First there is $w\mapsto w^*$%
\index{wstar@$w^*$, for $w$ an element of $\idd$}
 induced by the natural order reversing
involution $j\mapsto j^*$ on $\{1,\ldots,2d\}$: 
here $w^*$ has the obvious meaning,
namely,  it consists of all $j^*$ such that $j$ belongs to $w$.   
Then there is the map taking $w$ to its complement $\{1,\ldots,2d\}\setminus w$.
These two involutions commute.   Composing the two we get an order
preserving involution on $\idd$ which we denote by $w\mapsto w^\hash$\index{whash@$w^\hash$, for $w$ an element of $\idd$}.
The elements of the subset~$\id$ are fixed points under
this involution (there are points not in $\id$ that are also fixed).
\mysss{The involution $\hash$ on $\andposv$ and $\androotsv$}
\mylabel{sss.hashonposv} 
For $\alpha=(r,c)$ in $\andposv$, or more generally in $\androotsv$,
define $\alpha^\hash=(c^*,r^*)$.%
\index{alphahash@$\alpha^\hash$ for $\alpha$ in $\andposv$}      
The involution $\alpha\mapsto
\alpha^\hash$ is just the reflection with respect to the diagonal~$\diagv$.
For a subset or even multiset $\mon$ of $\andposv$ (or $\androotsv$),  
the symbol $\mon^\hash$%
\index{Shash@$\mon^\hash$, for monomial $\mon$ in $\andposv$ or $\androotsv$} 
has the obvious meaning.    We call $\mon$
{\em symmetric\/}\index{symmetric (monomial of $\andposv$)}
if $\mon=\mon^\hash$.  
\begin{proposition}\mylabel{p.distid}
An element $w\geq v$ of $\idd$ belongs to $\id$ if and only if the
distinguished subset~$\mon_w$ of~$\andposv$ corresponding to it
as described in~\S\ref{sss.montow}
is symmetric and has evenly many diagonal elements.
\end{proposition}
\begin{myproof}
That the symmetry of $\mon_w$ is equivalent to the condition that 
$w=w^\hash$ is proved in~\cite[Proposition~5.7]{gr}.     Now suppose
that $\mon_w$ is symmetric.  
We claim that for an element $(r,c)$ of $\mon_w$
that is not on the diagonal, either both $r$ and~$c$ are bigger
than $d$ or both are less than $d+1$.    It is enough to prove the claim,
for~$w$ is obtained from~$v$ by
removing the column indices and adding the row indices of elements of~$\mon_w$,
and it would follow that the number of entries in~$w$ that
are bigger than~$d$ equals the number of such entries in~$v$ plus the number
of diagonal elements in~$\mon_w$.

We now prove the claim.   
Since $\mon_w$ is symmetric, it follows 
that~$(c^*,r^*)$ also belongs to $\mon_w$.   Since $\mon_w$ is distinguished,  
it follows that 
in case $r<c^*$ (that is, if $(r,c)$ lies above the diagonal),  we have
$r<r^*$,  and so $c<r<r^*$;  and in case $r>c^*$,  we have $c^*<c$, and
so $c^*<c<r$.    Thus the claim is proved.
\end{myproof}
\mysubsection{The subset $\mon_C$ attached to a $v$-chain $C$}%
\mylabel{ss.montoC}
\mysss{Vertical and horizontal projections of an element of $\posv$}%
\mylabel{sss.projections}   For $\alpha=(r,c)$ in $\posv$ (or more
generally in $\rootsv$),  the 
elements $p_v(\alpha):=(c^*,c)$ and $p_h(\alpha):=(r,r^*)$ of the
diagonal~$\diagv$ are called respectively the {\em vertical and horizontal
projections of $\alpha$\/}.%
\index{vertical projection $p_v(\alpha)$}%
\index{horizontal projection $p_h(\alpha)$}%
\index{pv@$p_v(\alpha)$, vertical projection}%
\index{ph@$p_h(\alpha)$, horizontal projection}
In terms of pictures,  the vertical
projection is the element of the diagonal due South of $\alpha$;
the horizontal projection is the element of the diagonal due East of $\alpha$.
The vertical line joining $\alpha$ to its vertical projection $p_v(\alpha)$
and the horizontal line joining $\alpha$ to its horizontal projection 
$p_h(\alpha)$ are called the {\em legs\/}%
\index{legs of $\alpha$, for $\alpha\in\posv$} 
of $\alpha$.
\mysss{The ``connection'' relation%
\index{connectedness of two succcessive elements in a \vchain} on elements of a $v$-chain}%
\mylabel{sss.connection}
Let $C: \alpha_1=(r_1,c_1)>\alpha_2=(r_2,c_2)>\cdots$ be a $v$-chain
in~$\posv$.   Two consecutive elements $\alpha_j$ and $\alpha_{j+1}$ of $C$ 
are said to be {\em connected\/} if the following conditions are both satisfied:
\begin{itemize}
\item their legs are ``intertwined''%
\index{legs,  intertwining of};
equivalently and more precisely, 
this means that $r_j^*\geq c_{j+1}$,
or, what amounts to the same, $r_j\leq c_{j+1}^*$.
\item the point $(r_{j+1}, r_j^*)$ belongs to $\andposv$; this just means
that $r_{j+1}> r_j^*$. 
\end{itemize}
Consider the coarsest equivalence relation on the elements of $C$ generated
by the above relation.   The equivalence classes of $C$ with respect to this
equivalence relation are called the {\em connected components\/}%
\index{connected components of a \vchain} of the
$v$-chain $C$.

This definition has its quirks:\mcomment{picture added}
\\[1mm]
\begin{minipage}{6cm}
{\small
The $v$-chain~$C:
\alpha>\beta>\gamma$ in the picture has  
$\{\alpha,\beta\}$ and~$\{\gamma\}$ as its connected components;
but the ``sub''~$v$-chain $\alpha>\gamma$ of $C$ is connected
(as a $v$-chain in its own right).  }
\vfill
\end{minipage}
\hfill \begin{minipage}{6cm}
\setlength{\unitlength}{.34cm}
\begin{picture}(12,12)(-3,0)
\multiput(0,0)(12,0){2}{\line(0,1){12}}
\multiput(0,0)(0,12){2}{\line(1,0){12}}
\linethickness{.05mm}
\multiput(0,0)(1,0){13}{\line(0,1){12}}
\multiput(0,0)(0,1){13}{\line(1,0){12}}
\thicklines
\put(0,0){\line(1,1){12}}
\put(3.25,2.75){diagonal}
\thicklines
\put(9.25,6.25){boundary}\put(9.25,5.25){of $\andposv$}
\put(0,12){\line(1,0){3}}
\put(3,12){\line(0,-1){1}}
\put(3,11){\line(1,0){2}}
\put(5,11){\line(0,-1){2}}
\put(5,9){\line(1,0){2}}
\put(7,9){\line(0,-1){1}}
\put(7,8){\line(1,0){1}}
\put(8,8){\line(0,-1){1}}
\put(8,7){\line(1,0){1}}
\put(9,7){\line(0,-1){2}}
\put(9,5){\line(1,0){2}}
\put(11,5){\line(0,-1){2}}
\put(11,3){\line(1,0){1}}
\put(12,3){\line(0,-1){3}}
\put(1,4){\circle*{.4}}
\put(.25,4.5){$\alpha$}
\put(2,7){\circle*{.4}}
\put(1.25,7.5){$\beta$}
\put(3,10){\circle*{.4}}
\put(2.25,10.5){$\gamma$}
\multiput(1,4)(0,-.5){6}{\line(0,-1){.3}}
\multiput(1,4)(.5,0){6}{\line(1,0){.3}}
\multiput(2,7)(.5,0){10}{\line(1,0){.3}}
\multiput(2,7)(0,-.5){10}{\line(0,-1){.3}}
\multiput(3,10)(.5,0){14}{\line(1,0){.3}}
\multiput(3,10)(0,-.5){14}{\line(0,-1){.3}}
\put(4,10){\circle{.4}}
\put(4,10.5){$\in\andposv$}
\put(7,10){\circle{.4}}
\put(7,10.5){$\not\in\andposv$}
\multiput(4,4)(0,.5){12}{\line(0,1){.3}}
\multiput(7,7)(0,.5){6}{\line(0,1){.3}}
\end{picture}
\end{minipage}
\\

\mysss{The definition of $\mon_C$}%
\mylabel{sss.montoC}  We will define $\mon_C$\index{SC@$\mon_C$, where $C$ is a $v$-chain} as a multiset of $\andposv$.
It is easy to see and in any case stated explicitly as part of 
Corollary~\ref{p.montoC} that it is multiplicity free and
so is actually a subset of $\andposv$.

First suppose that $C:\alpha_1=(r_1,c_1)>\cdots>\alpha_\ell=(r_\ell,c_\ell)$
is a connected $v$-chain in $\posv$.  Observe that,
if there is at all an integer $j$, $1\leq j\leq \ell$,
such that the horizontal projection $p_h(\alpha_j)$ does
not belong to $\andposv$, then $j=\ell$.   
Define
\[ \mon_C:=
\left\{\begin{array}{ll}
	\{p_v(\alpha_1),\ldots,p_v(\alpha_\ell)\}
		& \textrm{if $\ell$ is even}\\
	\{p_v(\alpha_1),\ldots,p_v(\alpha_\ell)\}\cup\{p_h(\alpha_\ell)\}
		& \textrm{if $\ell$ is odd and $p_h(\alpha_\ell)\in\andposv$}\\
	\{p_v(\alpha_1),\ldots,p_v(\alpha_{\ell-1})\}
		\cup\{\alpha_\ell,\alpha_\ell^\hash\}
	   & \textrm{if $\ell$ is odd and $p_h(\alpha_\ell)\not\in\andposv$}
				\end{array}\right.
							\]

For a $v$-chain $C$ that is not necessarily connected, 
let $C=C_1\cup C_2\cup\cdots$ be the partition of $C$ into its connected
components,  and set
\[ \mon_C:= \mon_{C_1}\cup\mon_{C_2}\cup\cdots	\]

\mysubsubsection{The type of an element $\alpha$ of a $v$-chain $C$, and the
set $\mon_{C,\alpha}$}\mylabel{sss.type} 
We introduce some terminology and notation.    Their usefulness
may not be immediately apparent.

Suppose that~$C:\alpha_1>\cdots>\alpha_\ell$ is a connected $v$-chain.
We define the {\em type in $C$\/}%
\index{type (V, H, S), of an element in a $v$-chain|(}
of an element $\alpha_j$, $1\leq j\leq \ell$, of~$C$
to be V,~H, or~S, accordingly as: 
\begin{enumerate}
\item[V:] \hspace{2pt}$j\neq\ell$, or $j=\ell$ and $\ell$ is even.
\item[H:] \hspace{2pt}$j=\ell$, $\ell$ is odd,
			and $p_h(\alpha_\ell)\in\andposv$.
\item[S:] \hspace{2pt}$j=\ell$, $\ell$ is odd,
			and $p_h(\alpha_\ell)\not\in\andposv$.
\end{enumerate}
The {\em type\/}
of an element in a $v$-chain that is not necessarily connected
is defined to be its type in its connected component. 
\index{type (V, H, S), of an element in a $v$-chain|)}

The set $\mon_{C,\alpha}$%
\index{SCalpha@$\mon_{C,\alpha}$, for $\alpha$ in a \vchain~$C$}
of elements of $\andposv$ {\em generated\/} by an element~$\alpha$
of~$C$ is defined to be:
\[
\mon_{C,\alpha}:=\left\{
\begin{array}{ll}
\{p_v(\alpha)\}  & \textrm{if $\alpha$ is of type~V in~$C$};\\
\{p_v(\alpha), p_h(\alpha)\}  & \textrm{if $\alpha$ is of type~H in~$C$};\\
\{\alpha, \alpha^\hash\}  & \textrm{if $\alpha$ is of type~S in~$C$};
\end{array}\right.			\]
Observe that, for a $v$-chain $C$, the monomial $\mon_C$ defined
in~\S\ref{sss.montoC} is the union, over all elements $\alpha$ of $C$,
of $\mon_{C,\alpha}$.

For an element $\alpha$ of a $v$-chain $C$,  we define $\myrho_{C,\alpha}$%
\index{qcalpha@$\myrho_{C,\alpha}$, for $\alpha$ in a \vchain~$C$}
to be $p_v(\alpha)$ if $\alpha$ is of type~V or~H and to be $\alpha$ if
it is of type~S.

If the horizontal projection of an element in a \vchain\ does not belong
to~$\andposv$,  then clearly the same is true for every succeeding element.
The first such element of a $v$-chain is called the {\em critical\/}\index{critical element (of a $v$-chain)} element.

\begin{proposition}\mylabel{p.type}
\begin{enumerate}
\item\mylabel{p.type.odd}
The cardinality is odd of a connected component that has an element of type~H
or~S. Conversely, if the cardinality of a component is odd,  then it
has an element of type~H or~S.
\item\mylabel{p.type.hs}
An element of type H or S can only be the last element in its connected
component. 
\item\mylabel{p.type.cr}
The critical element has type either V or S.  No element before
it can be of type S and every element after it is of type S.  In particular,
any element that succeeds an element of type~S is of type~S.
\end{enumerate}
\end{proposition}
\begin{myproof} Clear from definitions.\end{myproof}
\begin{proposition}\mylabel{p.connect}
Let $\alpha>\gamma$ be elements of a $v$-chain~$C$ (we are not assuming
that they are consecutive).
\begin{enumerate}
\item
If $\alpha>\gamma$ is connected as a $v$-chain in its own right,
then $\alpha$
is connected to its next member in~$C$;  that is,  $\alpha$ cannot
be the last element in its connected component in~$C$.
\item
If $\alpha>\gamma$ is not connected as a $v$-chain in its own right
and the legs of $\alpha$ and $\gamma$ intertwine,
then the connected component of $\gamma$ in~$C$ is the singleton~$\{\gamma\}$,
and $\gamma$ has type~S in~$C$.
\end{enumerate}
\end{proposition}
\begin{myproof}
Clear from definitions.
\end{myproof}
\newcommand\dprime{D^\prime}
\newcommand\cprime{C^\prime}
\begin{proposition}\mylabel{p.concat}
Let $E:\alpha>\ldots>\zeta$ be a \vchain, $D$ and~$\dprime$ two
$v$-chains with tail~$\alpha$, and~$C$,~$\cprime$ the concatenations
of~$D$,~$\dprime$ respectively with~$E$.    Then
\begin{enumerate}
\item\mylabel{p.concat.1}
The last element in the connected component containing~$\alpha$ is the
same in~$C$ and~$\cprime$ (and this is the same as in~$E$). 
	\newcounter{myenumi}     
	\setcounter{myenumi}{\theenumi}  
\end{enumerate}
Let $\lambda$ denote this element.
\begin{enumerate}
\setcounter{enumi}{\value{myenumi}}  
\item\mylabel{p.concat.2}
The only element among $\alpha$, \ldots, $\zeta$ that possibly has
different types in $C$ and $\cprime$ is~$\lambda$.
\end{enumerate}
\end{proposition}
\begin{myproof}
(\ref{p.concat.1}): Whether or not two successive elements in a \vchain\ are
connected is independent of other elements in the \vchain.

(\ref{p.concat.2}): The type of an element in a \vchain\ is~V unless it is the last element
in its connected component.    And the type of the last element in a
component depends on the cardinality of the component.   The components
of~$E$ not containing~$\alpha$ are still components in~$C$ and~$\cprime$.
In contrast, the component containing~$\alpha$ could possibly be larger
in $C$ (respectively~$\cprime$) and hence its cardinality could be different.
\end{myproof}

For an element $\alpha=(r,c)$
of $\andposv$,   we define $\alpha\up$\index{alphaup@$\alpha\up$, for
$\alpha\in\andposv$} to be $\alpha$ itself if
$\alpha$ is either on or above the diagonal~$\diagv$ (more precisely, if
$r\leq c^*$),   and to be its ``reflection'' in the diagonal
(more precisely, $(c^*,r^*)$) if $\alpha$ is below the
diagonal (more precisely, if $r>c^*$).    For a monomial $\mon$ of~$\andposv$,
$\mon\up$\index{Sup@$\mon\up$, for a monomial $\mon$}
is defined to be the intersection of~$\mon$ (as a multiset) with the
subset~$\posv\cup\diagv$ of~$\andposv$. 
The notations $\alpha\down$%
\index{alphadown@$\alpha\down$, for $\alpha\in\andposv$} and
$\mon\down$\index{Sdown@$\mon\down$, for a monomial $\mon$}
have similar meanings.
\begin{quote}
{\footnotesize
Caution:   It is not true that $\mon\up=\{\alpha\up|\alpha\in\mon\}$ (in the obvious sense one would make of the right hand side). In particular, for a singleton monomial $\{\alpha\}$, it is not always true that $\{\alpha\}\up=\{\alpha\up\}$.}
\end{quote}

\begin{proposition}\mylabel{p.type.>}
Let $\alpha$ and $\beta$ be elements of a $v$-chain~$C$.    Let us use
$\alpha'$ and $\beta'$ respectively to denote elements of
$\mon_{C,\alpha}\textup{(up)}$
and $\mon_{C,\beta}\textup{(up)}$.   
\begin{enumerate}
\item If $\alpha>\beta$ (these elements are not necessarily consecutive
in~$C$),   then,  given $\beta'$,  there exists
$\alpha'$ such that $\alpha'>\beta'$.   In fact, this is true for every
choice of $\alpha'$ except when
\begin{equation*}\textup{(*)\quad
$\alpha$ is of type~H,  
and $p_h(\alpha)\not>\beta'$ for some
$\beta'\in\mon_{C,\beta}$.}
\end{equation*}
In particular, $\myrho_{C,\alpha}>\beta'$
and $\myrho_{C,\alpha}>\myrho_{C,\beta}$.
\item
Conversely,  suppose that $\alpha'>\beta'$ for some choice of $\alpha'$ and
$\beta'$.  Then $\alpha\geq\beta$;  if equality occurs, then
$\alpha$ is of type~H, $\alpha'=p_v(\alpha)$ and $\beta'=p_h(\alpha)$.
In particular,  if $\alpha'>
\myrho_{C,\beta}$ (or more specially $\myrho_{C,\alpha}>\myrho_{C,\beta}$),
then $\alpha>\beta$.
\item
If (*) holds for $\alpha>\beta$ in $C$,  then
\begin{enumerate}
\item the critical element of $C$ is the one just after $\alpha$;
in particular, $\alpha$ is uniquely determined.
\item all elements of $C$ succeeding $\alpha$ are of type~S;
in particular,  $\beta$ is of type~S and $\beta'=\beta$.
\item (*) holds for $\gamma$ in place of $\beta$ for every $\gamma$ in $C$ that
	succeeds $\alpha$.
\end{enumerate}
\end{enumerate}
\end{proposition}
\begin{myproof}
(1) If $\alpha$ is of type~V or~H,  we need only take
$\alpha'=p_v(\alpha)$, for
$p_v(\alpha)>p_v(\beta)$, $p_v(\alpha)>p_h(\beta)$,  and
$p_v(\alpha)>\beta$. 
Now suppose
that $\alpha$ is of type~S.  Then
$\beta$ too is of type~S (Proposition~\ref{p.type}~(\ref{p.type.cr})),
so $\beta'$ can only be~$\beta$, and the first
part of~(1) is proved.

It follows from the above that if $\alpha'=p_v(\alpha)$ or if $\alpha$
has type~S,  then $\alpha'>\beta'$ independent of the choice of
$\alpha'$.   So if $\alpha'\not>\beta'$,  then~(*) holds and 
$\alpha'=p_h(\alpha)$.

(3)~Let $\lambda$ be the immediate successor of $\alpha$ in~$C$.
Then $\alpha$ is not connected to $\lambda$
(Proposition~\ref{p.type}~(\ref{p.type.hs})).
Since $p_h(\alpha)\not>\beta'$, it follows that 
$\alpha$ and $\beta$ have intertwining legs.  Therefore so do
$\alpha$ and $\lambda$.    By Proposition~\ref{p.connect}~(2), 
$\lambda$ has type~S in~$C$. 

Since $\alpha$ has type~H and $\lambda$ type~S,   it follows
immediately from the definition of the critical element that $\lambda$ is
the critical element.   This proves~(a).    Assertion~(b) now follows from
Proposition~\ref{p.type}~(\ref{p.type.cr}).   For~(c),  write 
$p_h(\alpha)=(a,a^*)$,  $\lambda=(R,C)$,  and $\gamma=(r,c)$.
Then $R<a^*$,  for $\alpha$ and $\lambda$ have intertwining legs but
are not connected.   So $c<r\leq R<a^*$.   This means $p_h(\alpha)\not>\gamma$.
And $\gamma$ being of type~S (by (b)),   we can take $\gamma'=\gamma$.

(2) Suppose that $\alpha\not\geq\beta$.
Then $\beta>\alpha$.
By the second part of~(1) above, $\beta$ is of type~H and
$\beta'=p_h(\beta)$; by item~(b) of~(3), $\alpha$ is of type~S, so
$\alpha'=\alpha$.   This leads to the contradiction~$\beta>\alpha>p_h(\beta)$.
\end{myproof}

\begin{corollary}\mylabel{p.montoC}
The multiset $\mon_C$ attached to a $v$-chain~$C$ is
a distinguished subset of~$\andposv$ in the sense 
of~\ref{sss.ddist}.\end{corollary}
\begin{myproof}
If $\alpha$ in $C$ is of type~V or~S,  then $\mon_{C,\alpha}$ is a singleton;
if it is of type~H,  then $\mon_{C,\alpha}=\{p_v(\alpha),p_h(\beta)\}$.  So
there can be no violation of conditions~A and~B of \S\ref{sss.ddist} by
elements of $\mon_{C,\alpha}$.

Suppose $\alpha>\beta$.   By Proposition~\ref{p.type.>}~(1),  we have
$\alpha'>\beta'$ for any choice of $\alpha'\in\mon_{C,\alpha}$ and
$\beta'\in\mon_{C,\beta}$ except when the condition~(*) holds.  
By~(3) of the same proposition,  if (*) holds,   then $\beta'=\beta$,
and writing $\beta=(r,c)$, $p_h(\alpha)=(a,a^*)$,  we have
$r<a$ (since $\alpha>\beta$) and $c<r<a^*$ (see proof of item~3(c) of the
proposition).   Thus there can be no violation of conditions~A and~B of
\S\ref{sss.ddist}.
\end{myproof}
\begin{corollary}\mylabel{c.odtod}
Let $\mon$ be a $v$-chain in $\posv$ and $w$ an element of $I(d)$.
If $w$ \ortho-dominates $\mon$,  then $w$ dominates in the sense
of~\textup{\cite{kr}} the monomial $\mon\cup\mon^\hash$ of~$\andposv$.
\end{corollary}
\begin{myproof}
By \cite[Proposition~5.15]{gr},  it is enough to show that $w$ dominates~$\mon$.
Let $C:\alpha_1>\ldots>\alpha_t$ be a $v$-chain in $\mon$.
Writing $\alpha_j=(r_j,c_j)$ and $q_{C,\alpha_j}=(R_j,C_j)$ we have
$r_j\leq R_j$ and $C_j\leq c_j$.     By Proposition~\ref{p.type.>}~(1),
we have $q_{C,\alpha_1}>\ldots>q_{C,\alpha_t}$.    Since $w$ 
\ortho-dominates $\mon$, it in particular dominates $q_{C,\alpha_1}>\ldots>
q_{C,\alpha_t}$ and so also~$C$.\end{myproof}
%

\mysection{\ortho-depth}\mylabel{s.odepth}
The concept of \odepth\ 
defined in~\S\ref{ss.odepth} below
plays a key role in this paper.  As the name suggests, 
it is the orthogonal analogue of the concept of depth of~\cite{kr}.   
In~\S\ref{ss.odandd} below, it is observed that the \odepth\ is no smaller
than depth in the sense of~\cite{kr}.
In~\S\ref{ss.od-lemmas}, some observations about the relation between
\ortho-depths and types of elements in $v$-chains are recorded. 
\mysubsection{Definition of \odepth}\mylabel{ss.odepth}
\index{odepth@$\ortho$-depth|(}
The {\em \odepth\ of an element~$\alpha$
in a $v$-chain $C$ in $\posv$\/} is the depth in~$\mon_C$
in the sense of~\cite{kr} of~$\myrho_{C,\alpha}$:  in other words, it is
 the depth in~$\mon_C$ 
of~$p_v(\alpha)$ in case~$\alpha$ is of
type V or H,  and of~$\alpha$ (equivalently of $\alpha^\hash$) in case
$\alpha$ is of type S.     It is denoted $\odepthinC{\alpha}$.
The {\em \odepth\ of an element~$\alpha$ in
a monomial~$\mon$ of~$\posv$\/} is the maximum,
over all $v$-chains~$C$ in~$\mon$ containing~$\alpha$,
of the \odepth\ of $\alpha$ in $C$.  
It is denoted $\odepthinmon{\alpha}$.  Finally,  the 
{\em \odepth\ of a monomial $\mon$ in $\posv$\/}
is the maximum of the \odepths\ in $\mon$ of all the elements of~$\mon$.
\index{odepth@$\ortho$-depth|)}

There is a conflict in the above definitions:  Is the \odepth\ of an
element of a $v$-chain $C$ the same as its depth as an element of the
monomial~$C$?   In other words,  could the \odepth\ of an element in
a $v$-chain be exceeded by its \odepth\ in a sub-chain?   The conflict
is resolved by the first item of the following proposition.

\begin{proposition}\mylabel{p.odepth}
\begin{enumerate}
\item\mylabel{p.odepth.sub}
For $v$-chains $C\subseteq D$,  the \odepth\ in~$C$ of an element of~$C$ 
is no more than its \odepth\ in $D$.
\item\mylabel{p.odepth.ini}
If a \vchain~$C$ is an initial segment of a~\vchain~$D$,  then the
\odepths\ in~$C$ and~$D$ of an element of~$C$ are the same.
\end{enumerate}
\end{proposition}
\noindent
\begin{myproof}
(\ref{p.odepth.sub}):~By an induction on the difference in the cardinalities
of~$D$ and~$C$,  we may assume that~$D$ has one more element than~$C$.
Call this extra element~$\delta$.   Suppose that~$\delta$ lies between
successive elements~$\alpha$ and~$\beta$ of~$C$ (the modifications needed
to cover the extreme cases when it 
goes at the beginning or the end are being left to the reader).

The only elements of~$C$ that could possibly undergo changes of type on addition
of~$\delta$ are $\alpha$ and the last element in the connected component
of $\beta$,  which let us call~$\beta'$.   If there are no type changes,
then $\mon_C\subseteq\mon_D$ and the assertion is immediate.
The only type change that $\alpha$ can undergo is from~H to~V.
The type changes that $\beta'$ can undergo are:  H to V; V to H; S to V; V to S.
An easy enumeration of cases shows that only one of $\alpha$ and $\beta'$
can undergo a type change.

We need not worry about changes from V to H 
for in this case $\mon_C\subseteq\mon_D$. 

First let us suppose
that $\alpha$ undergoes a change of type (from H to V).   
Then $\delta$ is connected to~$\alpha$.
It follows from Proposition~\ref{p.type}~(\ref{p.type.odd}) that 
$\delta$ has type~V in $D$:
the connected component of $\alpha$ in $C$ has odd number of elements,
so if $\delta$ happens to be the last element in its connected
component in~$D$,   the number of elements in that component will be even.
Replacing an occurrence of 
$p_h(\alpha)$ in a $v$-chain of $\mon_C$ by $p_v(\delta)$
would result in a $v$-chain in $\mon_D$ (by Proposition~\ref{p.type.>}~(1)),
and this case is settled.

Now suppose that $\beta'$ undergoes a type change. 
Then $\delta$ is connected to $\beta$ and $\delta$ is of type~V in~$D$
(Proposition~\ref{p.type}~(\ref{p.type.hs})).
Replacing by $p_v(\delta)$ any occurrence in a $v$-chain in $\mon_C$
of $p_v(\beta')$,
$p_h(\beta')$, $\beta'$ accordingly as the type of $\beta'$ in~$C$
is~V,~H, or~S,
(not necessarily in the same place but at an
appropriate place) would result in a $v$-chain in~$\mon_D$
(by Proposition~\ref{p.type.>}~(1)),  
and we see that the \odepth\ cannot decrease.

(\ref{p.odepth.ini}):
It follows from Proposition~\ref{p.type.>}~(2) that,
for an element~$\alpha$ of~$C$,
contributions to~$\mon_D$ from elements beyond $\alpha$
(in particular from those not in~$C$)
do not affect the depth in $\mon_D$ of~$\myrho_{D,\alpha}$.
Looking for the possibility of differences in types in~$C$ and~$D$
of elements of~$C$, we see that the only element
of~$C$ that has possibly a different
type in~$D$ is its last element.  And this too can change type
only from~H to~V.

The above two observations imply that the calculations of \odepths\ in~$C$
and~$D$ of an element~$\alpha$ of~$C$ are no different:
we would be considering the depth in $\mon_C$ and $\mon_D$ respectively
of the same element (either $p_v(\alpha)$ or~$\alpha$),  and the differences
in $\mon_D$ and $\mon_C$ have no effect on this consideration.
\end{myproof}
\begin{corollary}\mylabel{c.vchain1}
If $C\subseteq D$ are $v$-chains in $\posv$,  then $w_C\leq w_D$ (although
it is not always true that $\mon_C\subseteq \mon_D$). 
\end{corollary}
\begin{myproof}
By \cite[Lemma~5.5]{kr},  it is enough to show that every $v$-chain in
$\mon_C$ is dominated by $w_D$.    Let $\beta_1=(r_1,c_1)>\cdots
>\beta_t=(r_t,c_t)$ be an arbitrary $v$-chain in $\mon_C$.
To show that it is dominated by $w_D$,  it is enough,
by~\cite[Lemma~4.5]{kr},  to show the existence of 
a $v$-chain $(R_1,C_1)>\cdots>(R_t,C_t)$ in $\mon_D$ with
$r_j\leq R_j$ and $C_j\leq c_j$ for $1\leq j\leq t$.   Such a $v$-chain
exists by the proof of~(1) of Proposition~\ref{p.odepth}.
\end{myproof}
\begin{corollary}\mylabel{c.p.odepth}
\begin{enumerate}
\item
Let $\mon$ be a monomial in $\posv$ and $\alpha\in\mon$.  Then
there exists a \vchain~$C$ in~$\mon$ with tail~$\alpha$ such that
$\odepthinmon{\alpha}=\odepthinC{\alpha}$.
\item For elements $\alpha>\gamma$ in a $v$-chain~$C$ (these 
need not be consecutive), we have $\odepthinC{\alpha}
<\odepthinC{\gamma}$.
\item For elements $\alpha>\gamma$ of a monomial~$\mon$ in $\posv$,
we have $\odepthinmon{\alpha}
<\odepthinmon{\gamma}$.
\item No two elements of the same \odepth\ in a monomial in $\posv$
are comparable.
\end{enumerate}
\end{corollary}
\begin{myproof}
(1) This follows from (2) of the Proposition above and the
definition of \odepth.

(2) This follows from Proposition~\ref{p.type.>}~(1) and the definition of
\odepth.

(3) By (1),  there exists a $v$-chain $C$ with tail $\alpha$ such that
$\odepthinmon{\alpha}=\odepthinC{\alpha}$.   Concatenate $C$ with
$\alpha>\gamma$ and let~$D$ denote the resulting $v$-chain.   
By (2) of the Proposition above,  $\odepthinC{\alpha}=\odepthin{D}{\alpha}$.
By (2) above, $\odepthin{D}{\alpha}<\odepthin{D}{\gamma}$.   And finally,
$\odepthin{D}{\gamma}\leq\odepthinmon{\gamma}$ by the definition of
$\odepthinmon{\gamma}$.

(4) Immediate from~(3).
\end{myproof}
\begin{corollary}\mylabel{c.vchain2}
Let $\beta>\gamma$ be elements of a \vchain~$C$ of elements of~$\posv$.
Let $E$ be a \vchain\ in~$\mon_C$ with tail $q_{C,\gamma}$ and
length $\odepthinC{\gamma}$.   Then $q_{C,\beta}$ occurs in~$E$.
\end{corollary}
\begin{myproof}  It is enough to show that for $\alpha'\neq q_{C,\beta}$
in~$E$,  either $\alpha'>q_{C,\beta}$ or $q_{C,\beta}>\alpha'$.
Let $\alpha$ be in $C$ such that $q_{C,\beta}\neq\alpha'\in\mon_{C,\alpha}$. 
If $\beta\geq\alpha$,  then $q_{C,\beta}>\alpha'$ by 
Proposition~\ref{p.type.>}~(1).    If $\alpha>\beta$ and $\alpha'\not>
q_{C,\beta}$,  then, by~(1) and~(3) of the same proposition,
$\alpha'\not>q_{C,\gamma}$, a contradiction.
\end{myproof}
\mysubsection{\odepth\ and depth}\mylabel{ss.odandd}
\begin{lemma}\mylabel{c.od>d}\mylabel{l.od>d}
The \odepth\ of an element $\alpha$ in a monomial $\mon$ of $\posv$ is
no less than its depth (in the sense of~\textup{\cite{kr}}) in $\mon\cup\mon^\hash$.
\end{lemma}
\begin{myproof}   Let $C:\alpha_1>\ldots>\alpha_t$ be a $v$-chain in
$\mon\cup\mon^\hash$ with tail $\alpha_t=\alpha$,  where~$t$ is
the depth of $\alpha$ in $\mon\cup\mon^\hash$.
We then have $\alpha_1\textup{(up)}>\ldots>\alpha_t\textup{(up)}$,
so we may assume~$C$ to be in~$\mon$.
By Proposition~\ref{p.type.>}~(1),   $q_{C,\alpha_1}>\ldots>
q_{C,\alpha_t}$ in $\mon_C$.     So
$\textup{depth}_{\mon\cup\mon^\hash}(\alpha)=t\leq \textup{depth}_{\mon_C}(q_{C,\alpha_t})\leq
\odepthinmon{\alpha}$.\end{myproof}
\mysubsection{\odepth\ and type}\mylabel{ss.od-lemmas}
We begin by defining some useful terminology.   Let $(r,c)$ and
$(R,C)$ be two elements of $\roots$.  To say that
$(R,C)$ {\em dominates\/}\index{domination (among elements of $\roots$)} $(r,c)$ means that $r\leq R$ 
and $C\leq c$ (in terms of pictures, 
$(r,c)$ lies (not necessarily strictly) to the Northeast of $(R,C)$).
To say that they are 
{\em comparable\/}\index{comparability, of elements of $\roots$} means that 
either $(R,C)>(r,c)$ or $(r,c)>(R,C)$.    
While this is admittedly
strange,  there will arise no occasion for confusion.

For an integer~$i$, we let $i\odd$\index{i(odd)@$i\odd$, for an integer $i$} 
be the largest odd 
integer not bigger than $i$ and 
$i\even$\index{i(even)@$i\even$, for an integer $i$}
the smallest even integer
not smaller than~$i$.

\begin{lemma}\mylabel{l.odepths}
\begin{enumerate}\item \mylabel{l.i.odepths}
For consecutive elements $\alpha>\beta$ of a $v$-chain $C$, 
\[	\odepthinC{\beta}=
	\left\{\begin{array}{ll}
		\odepthinC{\alpha}+2 & \text{if and only if
		$\alpha$ is of}\\
		&\text{type~H and $p_h(\alpha)>\beta$}\\
		\odepthinC{\alpha}+1 & \text{otherwise}\\
	\end{array}\right.
			\]
\item\mylabel{l.i.parity}
For an element of a \vchain~$C$ such that either its horizontal
projection belongs to $\andposv$ or it is connected to its predecessor,  the
parity of its \odepth\ in~$C$ is the same as that of its ordinality
in its connected component in~$C$.
\item\mylabel{l.i.hone}
The \odepth\ in a $v$-chain of an element of type H is odd.
\item\mylabel{l.i.vone}
If in a \vchain\ an element of type~V is the last in its connected
component,  then its \odepth\ is even.
\item\mylabel{l.i.hits}
If in a $v$-chain~$C$ there is an element of \odepth\ $d$, then
\begin{enumerate}
\item
for every odd integer~$d'$ not exceeding~$d$,  there is in~$C$
an element of \odepth~$d'$.
\item 
if, for an even integer~$d'$ not exceeding $d$,  there is no element in~$C$
of \odepth~$d'$,  then the element~$\alpha$ in~$C$ of \ortho-depth~$d'-1$
is of type~H, and $p_h(\alpha)>\beta$,  where $\beta$ denotes the immediate
successor of~$\alpha$ in~$C$.
\end{enumerate}
\item  \mylabel{l.i.deptheven}
Let $C$ be a \vchain\ and $\alpha$ an element of type~H in $C$.
Then the depth in $\mon_C$ of $p_h(\alpha)$ equals $\odepthin{C}{\alpha}+1$.
In particular, this depth is even.
\end{enumerate}
\end{lemma}
\begin{myproof}
(\ref{l.i.odepths}):  From items~1 and~3(a) of Proposition~\ref{p.type.>},
it follows that, for $\gamma$ in $C$ with $\gamma>\alpha$,  if
$\gamma'\not>q_{C,\alpha}$ for some $\gamma'$ in~$\mon_{C,\gamma}$,  then
$\gamma'\not>q_{C,\beta}$.    Thus $\odepthinC{\beta}$ exceeds
$\odepthinC{\alpha}$ by the number of elements in $\mon_{C,\alpha}$ 
that dominate $q_{C,\beta}$.     This number is $1$ if $\alpha$ is of
type~V,  or of type~S,  or of type~H and $p_h(\alpha)\not>\beta$;  it
is~$2$ if $\alpha$ is of type~H and $p_h(\alpha)>\beta$ (note that
$p_h(\alpha)>\beta$ if and only if $p_h(\alpha)>q_{C,\beta}$).

(\ref{l.i.parity}):
Let $\lambda$ be such an element.
Everything preceding $\lambda$ in~$C$ is of type~H
or~V~(Proposition~\ref{p.type}~(\ref{p.type.cr})).  
Let $\lambda$ belong to the $k^\textup{th}$
connected component, and 
$n_1$, \ldots, $n_k$ be respectively
the cardinalities of the first, \ldots, $k^\textup{th}$ connected components.  
By~(1) above and item~3(b) of~Proposition~\ref{p.type.>},     
$\odepthinC{\lambda}$ is
$n_1\even+\cdots+n_{k-1}\even$ plus the ordinality of $\lambda$ in the
$k^\textup{th}$ connected component.

(\ref{l.i.hone}) and (\ref{l.i.vone}): These are special cases of~\ref{l.i.parity}.

(\ref{l.i.hits}): This follows easily from~(\ref{l.i.odepths})
and~(\ref{l.i.hone}).

(\ref{l.i.deptheven}):
It follows from Proposition~\ref{p.type.>}~(2) that there is no element~$\gamma$
in~$\mon_C$ that lies between $p_v(\alpha)$ and $p_h(\alpha)$ (meaning
$p_v(\alpha)>\gamma>p_h(\alpha)$),
so the assertion holds.
\end{myproof}
\begin{corollary}\mylabel{c.odbindsd}
For a \vchain~$C$ in $\posv$,
if the \ortho-depths of elements in~$C$ are bounded by~$k$,
then the depths of elements in $\mon_C$ are bounded by~$k\even$.
\end{corollary}
\begin{myproof}
The depth of $q_{C,\alpha}$ in $\mon_C$ for any $\alpha$ in $C$ is at most
$k$ by hypothesis.    An element of~$\mon_C$ that is not~$q_{C,\alpha}$ for
any~$\alpha$ in $C$ can only be of the form $p_h(\alpha)$ for some $\alpha$.
By Proposition~\ref{p.type.>},  $\depthin{\mon_C}{p_v(\alpha)}
=\depthin{\mon_C}{p_h(\alpha)}-1$,  which implies $\depthin{\mon_C}{p_h(\alpha)}
\leq k+1$.   If,  moreover,  $k$ is even,  then by~(\ref{l.i.hone}) of
Lemma~\ref{l.odepths} $\depthin{\mon_C}{p_h(\alpha)}=
\depthin{\mon_C}{p_v(\alpha)}+1\leq (k-1)+1=k$.
\end{myproof}
\begin{proposition}\mylabel{p.goodCexists}
Given a monomial~$\mon$ in~$\posv$ and an element~$\alpha$ in it,  there
exists a \vchain~$C$ in~$\mon$ with tail~$\alpha$ such that
$\odepthinC{\beta}=\odepthinmon{\beta}$ for every~$\beta$ in~$C$.
\end{proposition}
\begin{myproof}
Proceed by induction on~$d:=\odepthinmon{\alpha}$.   Choose a \vchain~$D$
in~$\mon$ with tail~$\alpha$ such that 
$\odepthin{D}{\alpha}=\odepthinmon{\alpha}$ (such a \vchain\ exists
by Corollary~\ref{c.p.odepth}~(1)).
Let~$\alpha'$ be the element in~$D$ just before~$\alpha$.   It follows from
item~(3) of Corollary~\ref{c.p.odepth} and item~(\ref{l.i.odepths}) of
Lemma~\ref{l.odepths} that $\odepthinmon{\alpha'}$ 
(as also $\odepthin{D}{\alpha'}$) is either~$d-1$ or~$d-2$.
By induction,  there exists a \vchain~$C'$ with tail~$\alpha'$
that has the desired property.    Let $C$ be the concatenation of~$C'$
with~$\alpha'>\alpha$.   

We claim that~$C$ has the desired property.
The only thing to be proved is that~$\odepthin{C}{\alpha}=d$.
By item~(\ref{l.i.odepths}) of Lemma~\ref{l.odepths}, we have
$\odepthinC{\alpha}\geq\odepthin{C'}{\alpha'}+1$.  In particular, the
claim is proved in case~$\odepthin{C'}{\alpha'}$ is~$d-1$,   
so let us assume that~$\odepthin{C'}{\alpha'}$ is~$d-2$. 
It now follows from the same item that $\alpha'$ has type~H in $D$
and $p_h(\alpha')>\alpha$; it further follows that it is enough to
show that~$\alpha'$ has type~H in $C$.   

Since~$\alpha'$ has type~H in $D$, it follows (from item~(\ref{p.type.hs}) 
of~Proposition~\ref{p.type}) that $\alpha'>\alpha$ is not connected
and (from item~(\ref{l.i.hone}) of Lemma~\ref{l.odepths}) that $d-2$ is odd.   
Now, by item~(\ref{l.i.vone}) of Lemma~\ref{l.odepths},
the type in $C'$ of $\alpha'$ cannot be~V, so it is H, and the claim
is proved.
\end{myproof}
\begin{corollary}\mylabel{c.existence}
Let $\mon$ be a monomial in $\posv$,
$\beta$ an element of $\mon$,  and $i$ an integer
such that $i<\odepthinmon{\beta}$.   Then
\mcomment{\tiny The notation of the statement is made use of 
when the corollary is invoked in Lemma~\ref{l.phipi.1} (which in
turn is invoked in the proof of Lemma~\ref{l.phipi}).}
\begin{enumerate}
\item[(a)] If $i$ is odd, there exists an element
$\alpha$ in $\mon$ of \odepth~$i$ such that $\alpha>\beta$.
\item[(b)] If $i$ is even and there is no element $\alpha$ in $\mon$
of \odepth~$i$ such that $\alpha>\beta$, 
then there is element $\alpha$ in $\mon$ of \odepth~$i-1$
such that $p_h(\alpha)>\beta$.
\end{enumerate}
\end{corollary}
\begin{myproof}  Choose a \vchain~$C$ in $\mon$ having tail $\beta$ 
and the good property of Proposition~\ref{p.goodCexists}.  Apply
Lemma~\ref{l.odepths}~(\ref{l.i.hits}). 
\end{myproof}
\begin{corollary}\mylabel{c.for.odom}
Let $C$ be a \vchain\ in $\posv$ with tail $\alpha$ such that 
$\odepthin{C}{\alpha}$ is odd.    Let $A$ be a \vchain\ in $\posv$
with head~$\alpha$,  and $D$ the concatenation of~$C$ with~$A$.
Let $C'$ denote the \vchain~$C\setminus\{\alpha\}$.
Then
\begin{enumerate}
\item
The type of an element of~$A$ is the same in both~$A$ and~$D$.
In particular, $\mon_A\subseteq\mon_D$ and $q_{A,\beta}=q_{D,\beta}$
for $\beta$ in~$A$.
\item
The type of an element of~$C'$ is the same in both~$C'$ and~$D$.
In particular, $\mon_{C'}\subseteq\mon_D$.
\item
$\mon_D=\mon_{C'}\cup\mon_A$ (disjoint union);   
letting $j_0:=\odepthin{C}{\alpha}$
we have $(\mon_D)^{j_0}=\mon_A$ and $(\mon_D)_1\cup\cdots\cup(\mon_D)_{j_0-1}
=\mon_{C'}$.   (For a monomial~$\mon$,  the subset of elements of depth
at least $i$ is denoted $\mon^i$,  and the subset of elements of depth
exactly~$i$ is denoted $\mon_i$.)
\end{enumerate}
\end{corollary}
\begin{myproof}
(1)~Generally (meaning without the assumption that $\odepthin{C}{\alpha}$ is
odd), the only element of~$A$ that could possibly have a different type in~$D$
is the last one in the first connected
component of~$A$;   whether or not
it changes type depends exactly upon whether or not the parity of the 
cardinality of its connected component in~$D$ is different from that in~$A$.
Under our hypothesis,  this parity does not change, for, by~(\ref{l.i.vone}) of
Lemma~\ref{l.odepths},  the type of~$\alpha$ in~$C$ is~H or~S, and so
the cardinality of the connected component of~$\alpha$ in~$C$ is odd.

(2)~Generally (meaning without the assumption that $\odepthin{C}{\alpha}$ is
odd), the only element of~$C'$ that could possibly have a different type in~$D$
is the last one of~$C'$;    it changes type if and only if 
it is connected to $\alpha$ and the cardinality of its connected component
in $C'$ is odd.
Under our hypothesis,  this cardinality is even, for the same reason
as in~(1).

(3)~That $\mon_D=\mon_{C'}\cup\mon_A$ (disjoint union) is an
immediate consequence of (1) and (2).
By Lemma~\ref{l.odepths}~(1),  $q_{A,\alpha}=q_{D,\alpha}$ dominates
every element of $\mon_A$,  so $\mon_A\subseteq(\mon_D)^{j_0}$
($\depthin{\mon_D}{q_{D,\alpha}}=\odepthin{D}{\alpha}
=\odepthin{C}{\alpha}=j_0$).     It is enough to prove the following claim:
every element of $\mon_{C'}$ has depth less than~$j_0$ in~$\mon_D$.
Let $\gamma'$ be an element of $\mon_{C'}$.   If $\gamma'>q_{D,\alpha}$
then the claim is clear.     If not, then,  by Proposition~\ref{p.type.>}~(1),
$\gamma'=p_h(\gamma)$.
By Lemma~\ref{l.odepths}~(\ref{l.i.hone}),
$\odepthin{D}{\gamma}$ is odd. 
Since the claim is already true for $q_{D,\gamma}=p_v(\gamma)$,    we have
$\odepthin{D}{\gamma}=\depthin{\mon_D}{p_v(\gamma)}\leq j_0-2$.
By~(\ref{l.i.deptheven}) of the same lemma, $\depthin{\mon_D}{\gamma'}=\odepthin{D}{\gamma}+1$,
so $\depthin{D}{\gamma'}\leq j_0-1$,  and the claim is proved.\nolinebreak\end{myproof}
\begin{proposition}\mylabel{p.od.shift}
Let $\mon$ be a monomial in~$\posv$ and $j$ an odd integer.
For~$\beta$ in
$\monsupjjpone(:=\{\alpha\in\mon\st\odepthinmon{\alpha}\geq j\})$%
\index{Supjjpone@$\monsupjjpone$, for monomial $\mon$ in $\posv$},
we have 
\[
\odepthin{\monsupjjpone}{\beta}
=\odepthinmon{\beta}-j+1\quad
\]
\end{proposition}
\begin{myproof}
Proceed by induction on~$j$.  For $j=1$, the assertion reduces to
a tautology.     Suppose that the assertion has been proved upto~$j$.
By the induction hypothesis,  we have $\monsupjptwojpthree
=(\monsupjjpone)^{3,4}$,  and we are reduced to proving the
assertion for $j=3$.

Let $A$ be a \vchain\ in~$\monsupthreefour$ with tail~$\beta$
and $\odepthin{A}{\beta}=\odepthin{\monsupthreefour}{\beta}$.
Let $\alpha$ be the head of~$A$.   We may assume that
$\odepthin{\mon}{\alpha}=3$ for, if $\odepthinmon{\alpha}>3$, we can find,
by Lemma~\ref{l.odepths}~(\ref{l.i.hits}),  $\alpha'$ of
\ortho-depth~$3$ in $\mon$ with $\alpha'>\alpha$,   and extending
$A$ by $\alpha'$ will not decrease the \ortho-depth in $A$ of $\beta$
(Proposition~\ref{p.odepth}~(1)).
Let $E$ be a \vchain\ in~$\mon_A$ with tail~$q_{A,\beta}$ 
and length $\odepthin{A}{\beta}$.     The head of $E$
is then $q_{A,\alpha}$ (see Proposition~\ref{p.type.>}~(1)).

Choose $C$ in $\mon$ with tail~$\alpha$ such that $\odepthinC{\alpha}=3$.
Let $D$ be the concatenation of $C$ with $A$.  By Corollary~\ref{c.for.odom},
$E$ is contained in $\mon_D$, $q_{D,\alpha}=q_{A,\alpha}$,  
and $q_{D,\beta}=q_{A,\beta}$.  
By Proposition~\ref{p.odepth}~(2),
the \ortho-depth of $\alpha$ is the same in~$D$ as in~$C$.
Choose a \vchain~$F$ in $\mon_D$ with tail $q_{D,\alpha}=q_{A,\alpha}$.
Concatenating $F$ with $E$ we get a \vchain\ in $\mon_D$ with tail
$q_{D,\beta}=q_{A,\beta}$ of length~$\odepthin{\monsupthreefour}
{\beta}+2$.   This proves that $\odepthin{\mon}{\beta}
\geq\odepthin{\monsupthreefour}{\beta}+2$.

To prove the reverse inequality,  we need only turn the
above proof on its head.
Let $D$ be a \vchain\ in~$\mon$ with tail~$\beta$ such that
$\odepthinmon{\beta}=\odepthin{D}{\beta}$.  
Let $G$ be a \vchain\ in~$\mon_D$ with tail $q_{D,\beta}$ and
length~$\odepthinmon{\beta}$.
There exists an
element $\alpha$ in $D$ of \ortho-depth~$3$ in $D$ (by
Lemma~\ref{l.odepths}~(\ref{l.i.hits})).  Let $C$ be the part
of $D$ upto and including $\alpha$,  and $A$ the part $\alpha>\ldots>\beta$.
By Proposition~\ref{p.odepth}~(2), $\odepthinC{\alpha}=3$ and,
as above, Corollary~\ref{c.for.odom} applies.

By Corollary~\ref{c.vchain2},  $q_{A,\alpha}=q_{D,\alpha}$ occurs
in~$G$.  The part~$F$ of $G$ upto and including $q_{A,\alpha}$
is of length at most~$3$,  and the part $E: q_{D,\alpha}>\ldots>q_{D,\beta}$
belongs also to $\mon_A$ (Proposition~\ref{p.type.>}~(2)).
Thus the length of $G$ is at most
$2$ more than the the length of $E$ which is
at most $\odepthin{\monsupthreefour}{\beta}$.
\end{myproof}
\begin{corollary}\mylabel{c.p.od.shift}
For odd integers $i$,~$j$, we have
$(\mon^{i,i+1})^{j,j+1}=\mon^{i+j-1,i+j}$.\hfill$\Box$
\end{corollary}
\begin{corollary}\mylabel{c1.p.concat}
Let $E:\alpha>\ldots>\zeta$ be a \vchain, $D$ and~$\dprime$ two
$v$-chains with tail~$\alpha$, and~$C$,~$\cprime$ the concatenations
of~$D$,~$\dprime$ respectively with~$E$.    Then
\begin{enumerate}
\item
$\odepthin{C}{\zeta}-\odepthin{C}{\alpha}\leq
\odepthin{\cprime}{\zeta}-\odepthin{\cprime}{\alpha}+1$;
\item
equality holds if and only if
the type of~$\lambda$ is~H in~$C$ and~V in~$\cprime$, and $p_h(\lambda)>\mu$,
where $\lambda$ is the last element in the connected component containing
$\alpha$ of~$E$ and~$\mu$ is the immediate successor in~$E$ of~$\lambda$.
\end{enumerate}
\end{corollary}
\begin{myproof}
These assertions follow from combining~(\ref{p.concat.2}) of
Proposition~\ref{p.concat}
with (\ref{l.i.odepths}) of Lemma~\ref{l.odepths}.
\end{myproof}
\begin{corollary}\mylabel{c.p.concat}
Let $\zeta$ be an element of a monomial~$\mon$ in $\posv$.  
Let $C$ be a \vchain\ in~$\mon$ with tail~$\zeta$ such that
$\odepthinC{\zeta}=\odepthinmon{\zeta}$.   Then
\begin{enumerate}
\item\mylabel{c.p.i.concat.1}
$\odepthinC{\alpha}\geq\odepthinmon{\alpha}-1$ for any $\alpha$ in $C$.
\item\mylabel{c.p.i.concat.2}
If $\odepthinC{\alpha}=\odepthinmon{\alpha}-1$ for some $\alpha$ in $C$,
then	\begin{enumerate}
	\item\mylabel{c.p.i.concat.2.1}
	letting $\lambda$ be the last element in the connected component
	containing~$\alpha$ and $\mu$ the element next to $\lambda$,
	the type of $\lambda$ in~$C$ is~H and $p_h(\lambda)>\mu$.
	\item\mylabel{c.p.i.concat.2.2}
	$\odepthinC{\gamma}=\odepthinmon{\gamma}-1$ for all $\gamma$ in~$C$
	between $\alpha$ and $\lambda$ (both inclusive).   
	\end{enumerate}
\end{enumerate}
\end{corollary}
\begin{myproofnobox}  
(\ref{c.p.i.concat.1})  Let $\alpha$ be in $C$.  Let $E$ denote the
part of~$C$ beyond (and including)~$\alpha$.    Let $\dprime$ be a \vchain\ in $\mon$
with tail $\alpha$ such that $\odepthin{\dprime}{\alpha}=\odepthinmon{\alpha}$.
Let $\cprime$ be the concatenation of $\dprime$ and $E$.   Applying
Proposition~\ref{c1.p.concat}~(1), we have 
\[ \odepthinC{\alpha}\geq\odepthinC{\zeta}-\odepthin{\cprime}{\zeta}
			+\odepthin{\cprime}{\alpha}-1.	\]
But
$\odepthinC{\zeta}-\odepthin{\cprime}{\zeta}=\odepthinmon{\zeta}
		-\odepthin{\cprime}{\zeta}\geq 0$, and,
by the choice of $\dprime$ and Proposition~\ref{p.odepth}~(\ref{p.odepth.ini}),
$\odepthin{\cprime}{\alpha}=\odepthin{\dprime}{\alpha}
		=\odepthinmon{\alpha}$. 

(\ref{c.p.i.concat.2}) Assertions~(a) and~(b) follow respectively 
from the ``only if '' and ``if'' parts of item~(2)
of~Proposition~\ref{c1.p.concat}.\hfill$\Box$
\end{myproofnobox}

\mysection{The map $\opi$}\mylabel{s.opi}
The purpose of this section is to describe the map~$\opi$\index{opi@$\opi$}.
The description is given in~\S\ref{ss.opi}.
It relies on certain claims which are proved 
in~\S\S\ref{ss.p.monjjp1},~\ref{ss.pf.p.opidef}.  Those proofs
in turn refer to results from~\S\ref{s.lemmas},  but there is
no circularity---to postpone the definition of $\opi$ until all
the results needed for it have been proved would hurt rather
than help readability.    The observations
in~\S\ref{ss.ortho.spec} are required only in~\S\ref{s.proof}.

The symbol~$j$ will be reserved
for an odd positive integer throughout this section.  
\mysubsection{Description of~$\opi$}\mylabel{ss.opi}
The map~$\opi$ takes as input a
monomial~$\mon$ in~$\posv$ and produces as output a pair~$(w,\mon')$,%
\index{wSprime@$(w,\mon')(:=\opi(\mon))$, for $\mon$ in $\posv$}
where $w$ is an element of~$\id$ such that $w\geq v$
and~$\mon'$ is 
a ``smaller'' monomial, possibly empty,
in~$\posv$.   If the input $\mon$ is empty, 
no output is produced (by definition).    So now suppose that~$\mon$
is non-empty.

We first partition $\mon$ into subsets according to the \odepths\ 
of its elements.  
Let $\monkpr$\index{Skpr@$\monkpr$, for monomial $\mon$ in $\posv$} be the sub-monomial of $\mon$ consisting of those
elements of $\mon$ that have \odepth~$k$---the superscript ``\textup{pr}'' is
short for ``preliminary''.   It follows from Corollary~\ref{c.p.odepth}~(4)
that there are no comparable elements in $\monkpr$ and so we can arrange
the elements of $\monkpr$ in ascending order of both row and column indices.
Let $\sigma_k$\index{sigmak@$\sigma_k$} be the last element of $\monkpr$ in this arrangement.

Let now $j$ be an odd integer.
We set \[ \monjpr:=\mon_j^\textup{pr}\cup\mon_{j+1}^\textup{pr}.\]%
\index{Sjjponepr@$\monjpr$, for $\mon$ in $\posv$, $j$ odd}
We say that {\em $\mon$ is truly orthogonal\index{truly orthogonal at $j$ ($j$ odd)} at $j$\/}
if $p_h(\sigma_j)$ belongs to $\andposv$
(that is, if $r>r^*$ where $\sigma_j=(r,c)$),

Let~$\monjjpone$\index{Sjjpone@$\monjjpone$, for $\mon$ in $\posv$, $j$ odd} denote the monomial in $\andposv$ defined by
$\monjjpone :=$ \[\left\{\begin{array}{lr}
		\left(\monjpr\setminus\{\sigma_j\}\right)
		\cup\left(\monjpr\setminus\{\sigma_j\}\right)^\hash
		\cup\{p_v(\sigma_j),p_h(\sigma_j)\}
		&
		\begin{array}{l}
		 \textrm{if $\mon$ is truly}\\
		 \textrm{orthogonal at $j$}\end{array}
					 \\
		\monjpr\cup\left(\monjpr\right)^\hash \quad\quad\quad
		 \textrm{otherwise}
		\end{array}\right.
		\]
Here $\monjpr\setminus\{\sigma_j\}$ and other terms on the right are to be
understood as multisets.
As proved in Corollary~\ref{c.p.monjjp1}~(\ref{i.1.c.p.monjjp1}) below,  $\monjjpone$ has
depth at most~$2$.  
Let~$\monj$%
\index{Sk@$\mon_k$, for monomial $\mon$ in $\posv$}
(respectively~$\monjpone$%
\index{Sk@$\mon_k$, for monomial $\mon$ in $\posv$}
be the subset (as a multiset) of elements of
depth~$1$ (respectively~$2$) of~$\monjjpone$.

Now, for every integer $k$,  we apply the map of~$\pi$ of~\cite[\S4]{kr}
to~$\mon_k$
to obtain a pair $(w(k),\mon_k^\prime)$, where~$w(k)$\index{wk@$w(k)$} 
is an element of~$\idd$ and~$\mon_k^\prime$
is a monomial in~$\andposv$.   Let~$\mon_{w(k)}$
be the distinguished monomial in $\andposv$ associated
to~$w(k)$---see~\S\ref{sss.montow}.
\begin{proposition}\mylabel{p.opidef}
\begin{enumerate}
\item\mylabel{i.1.p.opidef}
$\mon_{w(k)}$ and $\mon_k^\prime$ are symmetric.   And therefore so
are $\cup_k\mon_{w(k)}$ and $\cup_k\mon_k^\prime$.
\item\mylabel{i.2.p.opidef}
$\cup_k\mon_{w(k)}$ is a distinguished subset of~$\andposv$
(in particular, the $\mon_{w(k)}$ are disjoint).
\item\mylabel{i.3.p.opidef}
For $j$ an odd integer,  either
\begin{itemize}
\item both $\mon_{w(j)}$ and $\mon_{w(j+1)}$ meet the diagonal,  or
\item neither of them meets the diagonal,
\end{itemize} 
precisely as whether or not $\mon$ is truly orthogonal at~$j$.
And therefore $\cup_k\mon_{w(k)}$ has evenly many diagonal elements.
\item\mylabel{i.4.p.opidef}
No $\mon_k'$ intersects the diagonal.  And therefore neither does
$\cup_k\mon_k'$.
\end{enumerate}
\end{proposition}
\noindent
The proposition will be proved below in \S\ref{ss.pf.p.opidef}.

\begin{figure}[!t]
\begin{center}
\includegraphics[height=7.75cm,width=7.75cm]{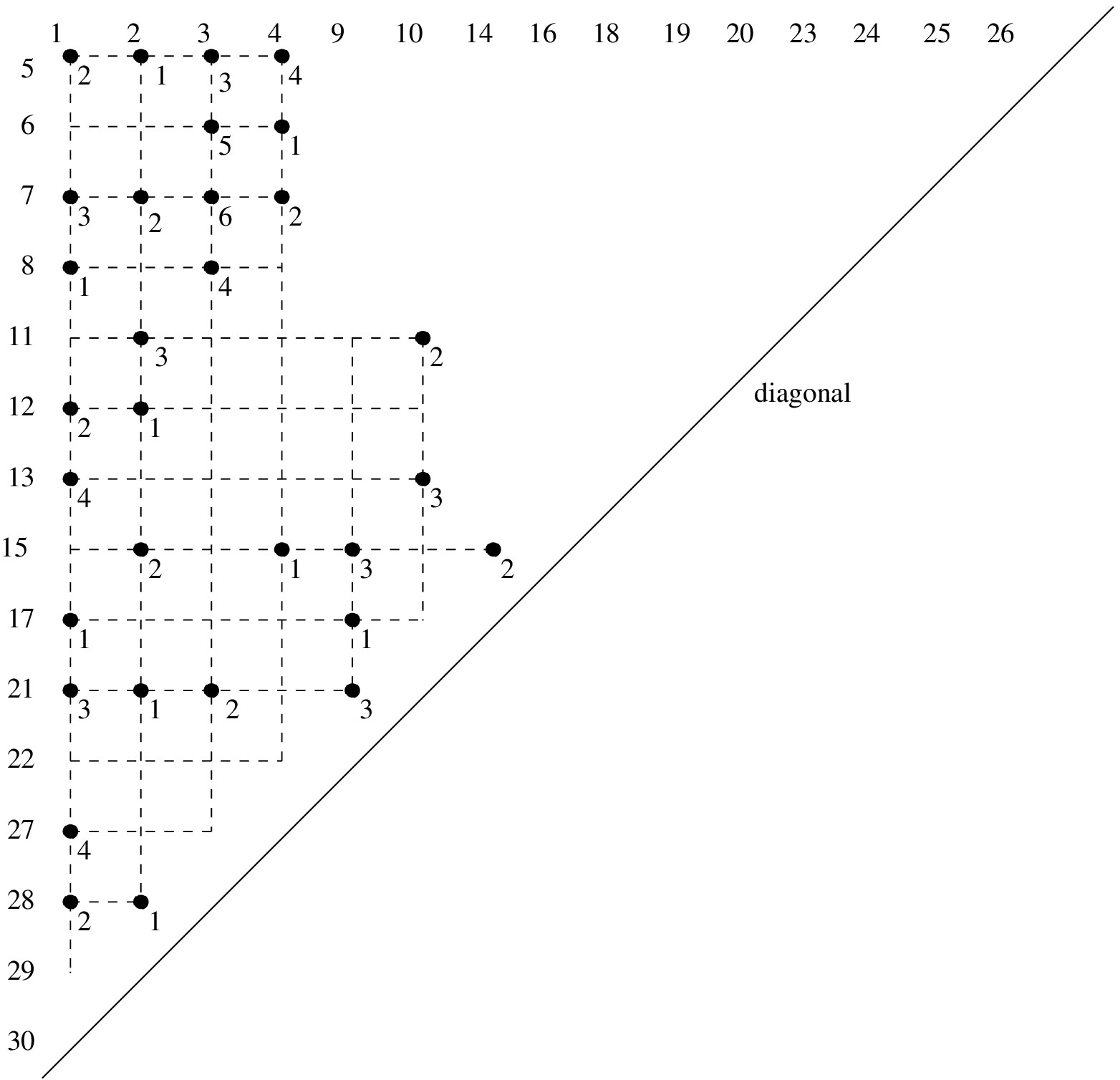}
\end{center}
\caption{\label{f.fig1}The monomial $\mon$}
\end{figure}

\index{wSprime@$(w,\mon')(:=\opi(\mon))$, for $\mon$ in $\posv$}
Finally we are ready to define the image $(w,\mon')$ of~$\mon$ under~$\opi$\index{opi@$\opi$}.
We let~$w$ be the element of~$\idd$ associated to the distinguished
subset~$\cup_k\mon_{w(k)}$ of~$\andposv$;    since $\cup_k\mon_{w(k)}$ is
symmetric and has evenly many diagonal elements,  it follows from
Proposition~\ref{p.distid} that~$w$ is in fact an element of~$\id$.
And we take $\mon':=\cup_k\mon_k^\prime\cap\posv$%
\index{Sprime@$\mon'$, for monomial $\mon$ in $\posv$}.
\index{wSprime@$(w,\mon')(:=\opi(\mon))$, for $\mon$ in $\posv$}
\bremark\mylabel{r.opidef}
Setting
\[\pi(\monjjpone):=(w_{j,j+1},\monjjpone'), \quad
       \mon':=\cup_{\textup{$j$ odd}}\monjjpone'\cap\posv,\]
and defining~$w$ to be the
element of~$\idd$ associated to $\cup_\textup{$j$ odd}\mon_{w_{j,j+1}}$
would give an equivalent definition of~$\opi$.
\eremark
\mysubsection{Illustration by an example}\mylabel{ss.example}
We illustrate the map $\opi$ by means of an example. 
Let $d=15$, and $v=(1,2,3,4,9,10,14,16,18,19,20,23,24,25,26)$.
A monomial $\mon$ in $\posv$ is shown in Figure~\ref{f.fig1}.
Solid black dots indicate the elements that occur in $\mon$ with
non-zero multiplicity. 
Integers written near the solid dots indicate multiplicities.
\begin{figure}[!b]
\begin{center}
\includegraphics[height=73.5mm,width=77.5mm]{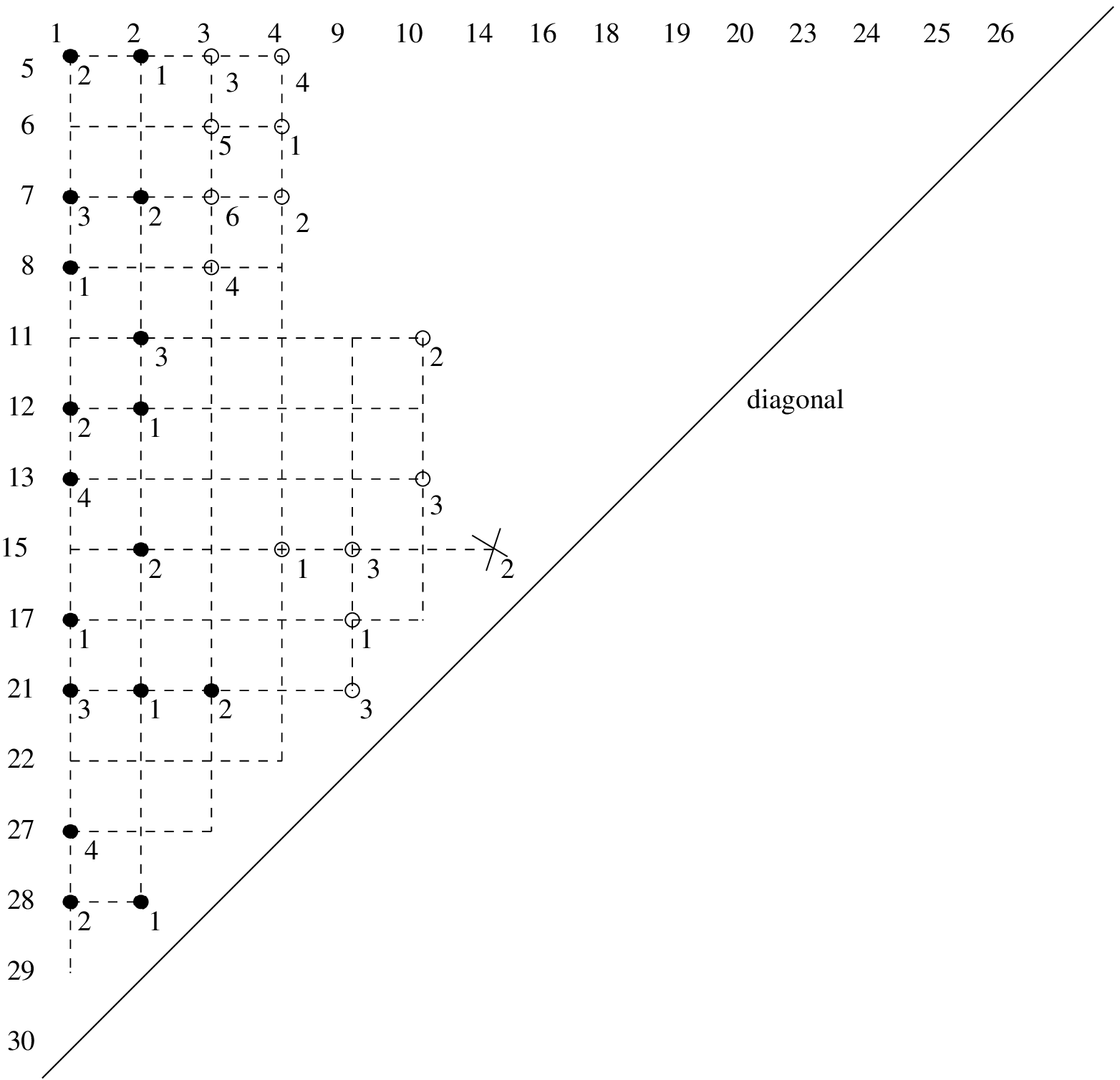}
\end{center}
\caption{\label{f.fig2} $\mon_{1,2}^\textup{pr}$, $\mon_{3,4}^\textup{pr}$, and
$\mon_{5,6}^\textup{pr}$}
\end{figure}
The \odepth\ of~$\mon$ is~$5$.    The element $(21,9)$ has \odepth~$3$
although it has depth~$2$ in~$\mon$.
Figure~\ref{f.fig2} shows the monomials $\mon_{1,2}^\textup{pr}$,
$\mon_{3,4}^\textup{pr}$, and $\mon_{5,6}^\textup{pr}$.
Solid dots, open dots, and crosses indicate elements of these monomials
respectively.
The monomial~$\mon$ is truly orthogonal at~$1$ and~$3$
but not at~$5$:
 $\sigma_1=(28,2)$, $\sigma_3=(21,9)$, and $\sigma_5=(15,14)$. 

\begin{figure}[!t]
\begin{center}
\includegraphics[height=7.5cm,width=7.75cm]{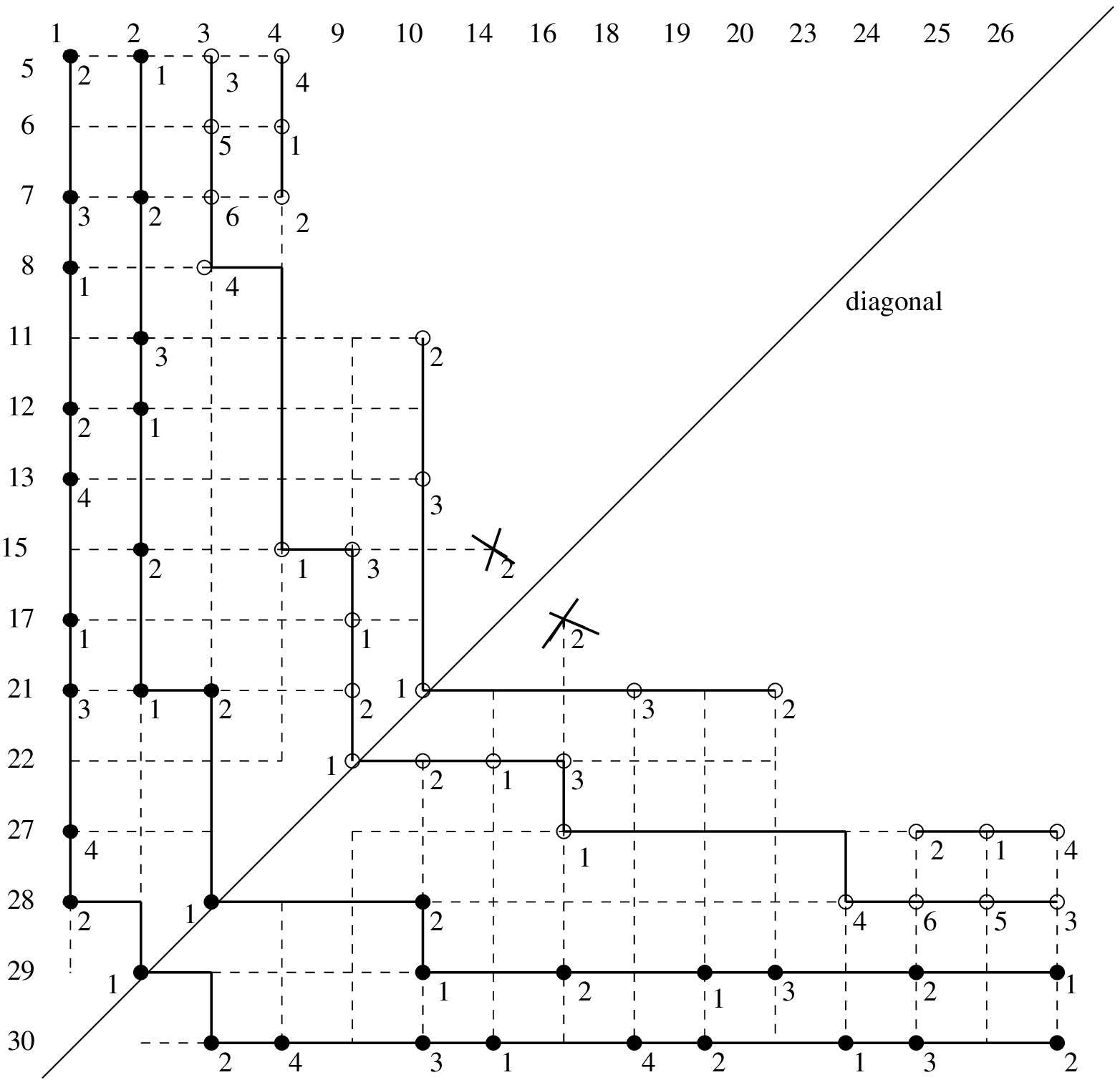}
\end{center}
\caption{\label{f.fig3} $\mon_{1,2}$, $\mon_{3,4}$, and $\mon_{5,6}$}
\end{figure}

\begin{figure}[!b]
\begin{center}
\includegraphics[height=7.5cm,width=7.75cm]{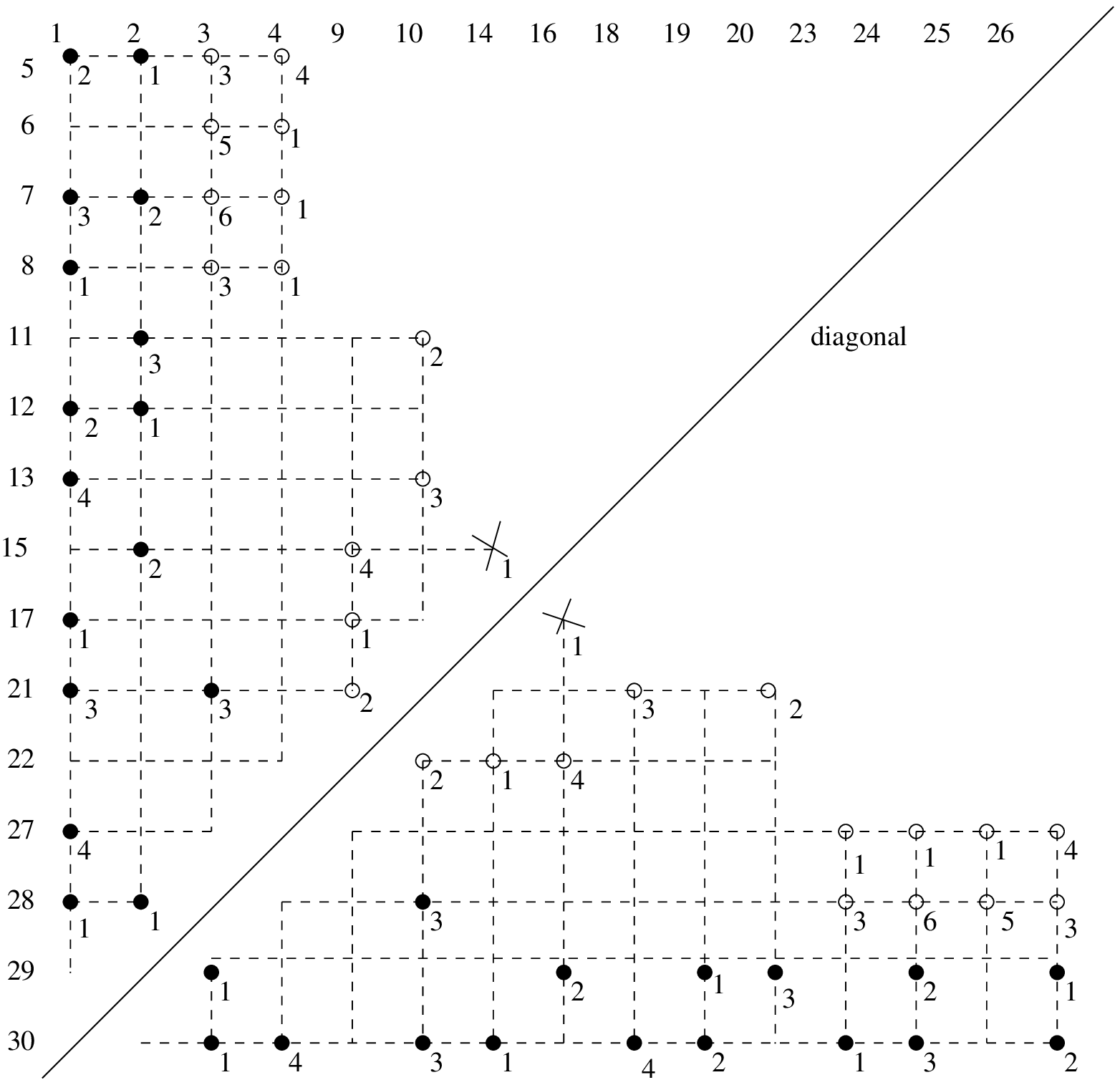}
\end{center}
\caption{\label{f.fig4} $\mon_{1,2}'$, $\mon_{3,4}'$, and $\mon_{5,6}'$}
\end{figure}

Figure~\ref{f.fig3} shows the monomials
$\mon_{1,2}$, $\mon_{3,4}$, and $\mon_{5,6}$ of~$\andposv$ and also their
decomposition into blocks,
and Figure~\ref{f.fig4} the monomials
$\mon'_{1,2}$, $\mon'_{3,4}$, and $\mon'_{5,6}$.

We have \[\monw=\{(15,14),(17,16),(21,10),(7,4),(27,24),(28,3),(30,1),(29,2)\}\]
hence $w=(7,9,15,17,18,19,20,21,23,25,26,27,28,29,30)$. 
It is easy to check that $w\in I(d)$. 
The monomial $\mon'$ is the intersection with $\posv$ of the union
of $\mon'_{1,2}$, $\mon'_{3,4}$, and $\mon'_{5,6}$---in other words
it is just the monomial lying above $\diag$ in Figure~\ref{f.fig4}.

\mysubsection{A proposition about $\monjjpone$}%
\mylabel{ss.p.monjjp1}
The aim of this subsection is to show that $\monjjpone$ has depth no
more than~$2$---see item~(\ref{i.1b.p.monjjp1})
of~Proposition~\ref{p.monjjp1}. 
This basic fact was mentioned above in the description of~$\opi$
and is necessary (psychologically although not logically) to make sense of
the definitions of~$\monj$ and $\monjpone$.       We prepare the way
for Proposition~\ref{p.monjjp1} by way of two preliminary propositions.
The first of these is about elements of \odepth~$j$ and $j+1$ in~$\mon$,
the second about the relation of these elements with~$\sigma_j$.

\begin{proposition}\mylabel{p.monkpr}
\begin{enumerate}
\item\mylabel{i.1.p.monkpr} $\monkpr$ has no comparable elements.
\item\mylabel{i.2.p.monkpr}
For $j$ an odd integer and $\beta$ an element of $\monjponepr$,  there
exists~$\alpha$ in $\monjonlypr$ such that $\alpha>\beta$.   In particular,
the row index of~$\sigma_{j+1}$ (if $\sigma_{j+1}$ exists) is less than
the row index of~$\sigma_j$.
\end{enumerate}
\end{proposition}
\begin{myproof} (\ref{i.1.p.monkpr})~follows from
Corollary~\ref{c.p.odepth}~(4); 
(\ref{i.2.p.monkpr})  follows from~Proposition~\ref{p.goodCexists}
and Lemma~\ref{l.odepths}~(\ref{l.i.hits}).
\end{myproof}
\begin{proposition}\mylabel{p.sigmaj}
Let $j$ be an odd integer and let $\mon$ be truly orthogonal at $j$.  Then
\begin{enumerate}
\item\mylabel{i.1.p.sigmaj}
$p_v(\sigma_j)>p_h(\sigma_j)$;   
if $\alpha>p_v(\sigma_j)$, then $\alpha>\sigma_j$; 
if $\alpha>\sigma_j$, then $\alpha>p_h(\sigma_j)$.
\item\mylabel{i.2.p.sigmaj}
No element of~$\monjonlypr$ is comparable to $p_v(\sigma_j)$ or $p_h(\sigma_j)$.
\item\mylabel{i.3.p.sigmaj}
No element of $\monjponepr$ is comparable to $p_h(\sigma_j)$.
\item\mylabel{i.4.p.sigmaj}
The following is not possible: $\alpha\in\monjonlypr$,
$\beta\in\monjponepr$, and $p_h(\alpha)>\beta$.
\end{enumerate}
\end{proposition}
\begin{myproof}
(\ref{i.1.p.sigmaj}) is trivial. (\ref{i.2.p.sigmaj}) follows immediately
from the definition of $\sigma_j$.   We now prove (\ref{i.3.p.sigmaj}).
First suppose $\beta>p_h(\sigma_j)$ for some $\beta$ in $\monjponepr$. 
By~(\ref{i.2.p.monkpr}) of Proposition~\ref{p.monkpr},
there exists $\alpha$ in~$\monjonlypr$ such that
$\alpha>\beta$. But then the row index of~$\alpha$
exceeds that of~$\sigma_j$,  a contradiction to the choice of $\sigma_j$.

We claim that it is not possible for $\beta\in\monjponepr$ to satisfy
$p_h(\sigma_j)>\beta$.      This being a special case of~(\ref{i.4.p.sigmaj}),  
we need only prove that statement.   So suppose that $\alpha$ belongs
to $\monjonlypr$ and that $p_h(\alpha)>\beta$.   
Let $C$ be a $v$-chain in $\mon$ with tail~$\alpha$ such that
$\odepthinC{\alpha}=j$ (see~Proposition~\ref{c.p.odepth}~(1)).
Concatenate $C$ with $\alpha>\beta$ and call
the resulting \vchain~$D$.    Then, by~Lemma~\ref{l.odepths}~(\ref{l.i.vone}),
$\alpha$ is of type~H in~$D$, so that,
by Lemma~\ref{l.odepths}~(\ref{l.i.odepths}), we have
$\odepthin{D}{\beta}=\odepthin{D}{\alpha}+2$.   But, by
Proposition~\ref{p.odepth}~(\ref{p.odepth.ini}), $\odepthin{D}{\alpha}
=\odepthinC{\alpha}=j$,  so that $\odepthinmon{\beta}\geq j+2$, 
a contradiction.
\end{myproof}

Let $\monjjponeext$\index{Sjjponeext@$\monjjponeext$, for $\mon$ in $\posv$, $j$ odd} denote the set---not multiset---defined by:
\[
\monjjponeext:=\left\{\begin{array}{ll}
		\monjjpone\cup\{\sigma_j,\sigma_j^\hash\} 
		&\textrm{ if $\mon$ is truly orthogonal at $j$}\\
		\monjjpone
		&\textrm{ otherwise}
		\end{array}\right.
					\]
Here $\monjjpone$ on the right stands for the underlying set of the
multiset~$\monjjpone$ defined above.
The set $\monjjponeext$ is the disjoint union of the sets $\monjext$
and~$\monjponeext$
defined as follows (here again the terms on the
right hand side denote the underlying sets of the corresponding multisets):
\[	\monjext:=\left\{\begin{array}{ll}
\index{Skext@$\mon_k\ext$, for monomial $\mon$ in $\posv$|(}
		\monjonlypr\cup\left(\monjonlypr\right)^\hash\cup\{p_v(\sigma)\}
		&\textrm{ if $\mon$ is truly orthogonal at $j$}\\
		\monjonlypr\cup\left(\monjonlypr\right)^\hash
		&\textrm{ otherwise}
		\end{array}\right.
					\]
\[	\monjponeext:=\left\{\begin{array}{ll}
		\monjponepr\cup\left(\monjponepr\right)^\hash\cup\{p_h(\sigma)\}
		&\textrm{ if $\mon$ is truly orthogonal at $j$}\\
		\monjponepr\cup\left(\monjponepr\right)^\hash
		&\textrm{ otherwise}
		\end{array}\right.
\index{Skext@$\mon_k\ext$, for monomial $\mon$ in $\posv$|)}
\]
\begin{proposition}\mylabel{p.monjjp1}
\begin{enumerate}
\item\mylabel{i.1.p.monjjp1}
$\monjext$ (respectively $\monjponeext$) is precisely the set of
elements of depth~$1$ (respectively~$2$) in $\monjjponeext$.  In particular,
	\begin{enumerate}
	\item\mylabel{i.1a.p.monjjp1}
	Neither $\monjext$ nor $\monjponeext$ contains comparable elements.
	\item\mylabel{i.1b.p.monjjp1}
	The length of a \vchain\ in~$\monjjponeext$ is at most~$2$.
	\item\mylabel{i.1c.p.monjjp1}
	There is a \vchain\ of length~$2$ in~$\monjjpone$
	unless~$\monjponeext$ is empty.
	\end{enumerate}
\item\mylabel{i.2.p.monjjp1}
Let $k$ be a positive integer, not necessarily odd.
If there is in $\mon$ an element of \odepth\ at least~$k$,
then $\mon_k\textup{(ext)}$ is non-empty.    The converse also holds
except possibly if $k$~is even and $\mon$ is truly orthogonal at $k-1$.
In particular,  if $\mon_k\textup{(ext)}$ is non-empty,   then there
is an element of \odepth~at least~$k-1$.
\end{enumerate}
\end{proposition}
\begin{myproof}
(\ref{i.1.p.monjjp1}):~It is enough to show that every element 
of~$\monjext\textup{(up)}$ (respectively~$\monjponeext\textup{(up)}$)
is of depth~$1$ (respectively~$2$) in $\monjjponeext\textup{(up)}$,  for 
\begin{itemize}
\item
$\alpha>\beta$ implies
$\alpha\textup{(up)}>\beta\textup{(up)}$ for elements $\alpha$,~$\beta$
of~$\andposv$. 
\item
$\monjjponeext=\monjext\cup\monjponeext$.
\item
$\monjjponeext$, $\monjext$, and $\monjponeext$ are symmetric.
\end{itemize}
In turn,  it is enough to show the following:
	\begin{enumerate}
	\item[(i)]
	Every element of $\monjext\textup{(up)}$ has depth~$1$.
	\item[(ii)]
	$\monjponeext\textup{(up)}$ has no comparable elements.
	\item[(iii)]
	Every element of $\monjponeext\textup{(up)}$ has depth at least~$2$.
	\end{enumerate}
Item~(i) follows from~Proposition~\ref{p.monkpr} and
Proposition~\ref{p.sigmaj}~(\ref{i.2.p.sigmaj});   
item~(ii) from~Pro\-po\-sition~\ref{p.monkpr}~(\ref{i.1.p.monkpr}) and
Proposition~\ref{p.sigmaj}~(\ref{i.3.p.sigmaj});   
item~(iii) from~Proposition~\ref{p.monkpr}~(\ref{i.2.p.monkpr}) and
Proposition~\ref{p.sigmaj}~(\ref{i.1.p.sigmaj}).

(\ref{i.2.p.monjjp1}):~The first assertion follows from
Lemma~\ref{l.odepths}~(\ref{l.i.hits}): if $k$ is odd there is an
element of \odepth~$k$ in~$\mon$;   if $k$ is even and there is no
element of \odepth~$k$ in~$\mon$,  then there is in~$\mon$ an element
of \odepth~$k-1$ and of type~H,  so $\mon$ is truly orthogonal at $k-1$.  
The second
assertion is clear from the definition of~$\mon_k\textup{(ext)}$.
\end{myproof}
\begin{corollary}\mylabel{c.p.monjjp1}
\begin{enumerate}
\item\mylabel{i.1.c.p.monjjp1}
No element of $\monjjpone$ has depth more than~$2$.
\item\mylabel{i.2.c.p.monjjp1}
$\monjponeext=\monjpone$ and $\monjext\cap\monjjpone=\monj$ (as sets).
In particular,  $\monjpone=\monjjpone\cap\monjponeext$ and 
$\monj=\monjjpone\cap\monjext$ as multisets defined by the intersection of
a multiset with a subset.
\end{enumerate}
\end{corollary}
\begin{myproof}
(\ref{i.1.c.p.monjjp1}):  Since $\monjjpone\subseteq\monjjponeext$ (as sets),
this follows immediately from~(\ref{i.1b.p.monjjp1}) of the proposition above.

(\ref{i.2.c.p.monjjp1}): 
Since the union of $\monjponeext$ (which always is contained in $\monjjpone$)
and $\monjext\cap\monjjpone$ is all of $\monjjpone$,  and since $\monj$,
$\monjpone$ are disjoint,  it it enough to show that $\monjponeext\subseteq
\monjpone$ and $\monjext\cap\monjjpone\subseteq\monj$.

Now, since elements of $\monjext$ have depth~$1$ even in $\monjjponeext$
(by item~(\ref{i.1.p.monjjp1}) of the proposition above),  it is immediate that
$\monjext\cap\monjjpone\subseteq\monj$.   And it follows from the proof
of item~(iii) in the proof of item~(\ref{i.1.p.monjjp1}) of the proposition
above that an element of $\monjponeext$ has depth~$2$ even in $\monjjpone$
(not just in~$\monjjponeext$),
so that $\monjponeext\subseteq\monjpone$.
\end{myproof}
\mcomment{Definitions of ``up'' and ``down'' to be put in.}
\mysubsection{Proof of Proposition~\ref{p.opidef}}\mylabel{ss.pf.p.opidef}
(\ref{i.1.p.opidef})
The monomials $\monjjpone$ are clearly symmetric.  Observe that 
$\alpha$ in $\monjjpone$ has the same depth as $\alpha^\hash$, for
$\alpha_1>\alpha_2$ implies $\alpha\textup{(up)}>\alpha_2\textup{(up)}$
and $\alpha\textup{(down)}>\alpha_2\textup{(down)}$ for 
$\alpha_1$,~$\alpha_2$ in~$\andposv$.      Thus the monomials $\mon_k$
are symmetric.    Since the map $\pi$ of~\cite{kr} respects~$\hash$---see
Proposition~5.7 of~\cite{gr}---it follows that $\mon_{w(k)}$ and $\mon_k'$
are symmetric.    Therefore so are $\cup_k\mon_{w(k)}$ and $\cup_k\mon_k'$.

(\ref{i.2.p.opidef})   This follows from Corollary~\ref{c.4.13.kr}.

(\ref{i.3.p.opidef})   If $\mon$ is truly orthogonal at $j$,  then
$p_v(\sigma_j)$ and $p_h(\sigma_j)$ are diagonal elements respectively
in $\monj$ and $\monjpone$---see
Corollary~\ref{c.p.monjjp1}~(\ref{i.2.c.p.monjjp1}).    Thus both $\monj$
and $\monjpone$ have diagonal blocks in the sense of Proposition~5.10~(A)
of~\cite{gr}.    It follows from the result just quoted that both $\mon_{w(j)}$
and $\mon_{w(j+1)}$ meet the diagonal.    It is of course clear that each
$\mon_{w(k)}$ meets the diagonal at most once since diagonal elements
are clearly comparable but elements of $\mon_{w(k)}$ are not
by Lemma~4.9 of~\cite{kr}.

Suppose that $\mon$ is not truly orthogonal at $j$.   Then $\sigma_j$
and $\sigma_j^\hash$ belong to different blocks---this is equivalent
to the definition of $\mon$ being not truly orthogonal at $j$.   By
Proposition~\ref{p.monkpr}~(\ref{i.2.p.monkpr}),  it follows that
$\sigma_{j+1}$ and $\sigma_{j+1}^\hash$ also belong to different blocks.
So neither $\monj$ nor $\monjpone$ has a diagonal block.

(\ref{i.4.p.opidef})    If $\mon$ is not truly orthogonal at $j$,
then neither $\monj$ nor $\monjpone$ has a diagonal block (as has just
been said above),   and it follows from Proposition~5.10~(A) of~\cite{gr}
that neither $\monj'$ nor $\monjpone'$ meets the diagonal.

So suppose that $\mon$ is truly orthogonal at $j$.   Then both $\monj$
and $\monjpone$ have a diagonal entry each of multiplicity~$1$, namely
$p_v(\sigma_j)$ and $p_h(\sigma_j)$ respectively.  
It is clear from the definition
of $\sigma_j$ that no element of $\monj\textup{(up)}$ shares
its row index with $p_v(\sigma_j)$.    And it follows from
Proposition~\ref{p.monkpr}~(\ref{i.2.p.monkpr}) that 
no element of $\monjpone\textup{(up)}$ shares its row
index with $p_h(\sigma_j)$.   It now follows from the proof of
Proposition~5.10~(B) of~\cite{gr}---see the last line of that proof---that
neither $\monj'$ nor $\monjpone'$ meets the diagonal.~\hfill$\Box$
\subsection{More observations}%
\mylabel{ss.ortho.spec}
\begin{proposition}\mylabel{p.d<3}
The length of any \vchain\ in~$\mon_{j,j+1}\cup\mon'_j\cup\mon'_{j+1}$
is at most~$2$.
\end{proposition}
\begin{myproof}
By Corollary~\ref{c.p.monjjp1}~(\ref{i.1.c.p.monjjp1}), 
the length of any \vchain\
in~$\monjjpone$ is at most~$2$.    Applying Lemma~\ref{l.d.monmon'}
to~$\monjjpone$,  we get the desired result.
\end{myproof}
\begin{proposition}\mylabel{p.prepare}
\begin{enumerate}
\item
For an element $\alpha'=(r,c)$ of $\monkprimeup$,  there exists
an element $\alpha=(r,C)$ of $\monkpr$ with $C\leq c$.
\item \mcomment{
It is not clear that the second item in Prop.~\ref{p.prepare}
is really needed.   It has been put in for reasons of balance.}
For an element $\alpha'=(r,c)$ of $\monjponeprimeup$,  there exists
an element $\alpha=(R,c)$ of $\mon_{j+1}\up$ with $r\leq R$.
\item
For an element $\alpha'$ of $\monjponeprimeup$,  there exists
an element $\alpha$ of $\monjonlypr$ with $\alpha>\alpha'$.
\end{enumerate}
\end{proposition}
\begin{myproof} (1) That there exists $\alpha$ in $\monk\up$ with
$C\leq c$ follows from the definition of $\monkprimeup$.  
Clearly such an $\alpha$ cannot be on the diagonal,   so 
$\alpha$ belongs to $\monkpr$.

(2) As in the proof of~(1),  it follows from the definition of
$\mon_{j+1}'$ that there exists $\alpha=(R,c)$ in $\mon_{j+1}$
with $r\leq R$.    If $\alpha$ lies strictly below the diagonal,
then $c>R^*$,  so that $\alpha^*=(c^*,R^*)>\alpha'=(r,c)$, 
a contradiction to~Lemma~\ref{l.d.monmon'} ($\alpha^*$ belongs
to $\mon_{j+1}$ by the symmetry of $\mon_{j+1})$.
Thus $\alpha$ belongs to $\mon_{j+1}\up$.

(3) Writing $\alpha'=(r,c)$, by~(1), we can find an $\beta=(r,C)$
in $\monjponepr$ with $C\leq c$.    
By Proposition~\ref{p.monkpr}~(\ref{i.2.p.monkpr}),  there exists
$\alpha$ in~$\monjpr$ such that $\alpha>\beta$.
\end{myproof}
\begin{corollary}\mylabel{c.p.prepare}
If in $\mon_j'\up\cup\mon_{j+1}'\up$ there exists an element with
horizontal projection in $\andposv$,  then $\mon$ is truly orthogonal
at~$j$.
\end{corollary}
\begin{myproof}
 Follows directly from Proposition~\ref{p.prepare} (1) and (3).
\end{myproof}

\begin{proposition}\mylabel{p.od.2}
The \odepth\ of an element in $\monjonlypr\cup\monjponepr$ is at most~$2$.
More strongly,  
the \odepth\ of an element in $\monjpr\cup\monjprimeup\cup\monjponeprimeup$
is at most~$2$.
\end{proposition}
\begin{myproof}
It is enough to show that no element in 
$\monjprimeup\cup\monjponeprimeup$ has \odepth\ more than~$2$, for
we may assume by increasing multiplicities that
$\monjonlypr\subseteq\monjprimeup$
and $\monjponepr\subseteq\monjponeprimeup$ (as sets).
It follows from Proposition~\ref{p.d<3} that a $v$-chain in $\monjprimeup
\cup\monjponeprimeup$ has length at most~$2$.  
Let $\alpha'_1=(r_1,c_1)>
\alpha'_2=(r_2,c_2)$ be such a $v$-chain.  
It follows from the proof of Corollary~4.14~(2)
of \cite{kr} that $\alpha'_1\in\monjprimeup$ and $\alpha'_2\in\monjponeprimeup$.
By item~(\ref{l.i.odepths}) of Lemma~\ref{l.odepths},
it is enough to rule out the following possibility:
$\alpha'_1$ is of type~H in $\alpha'_1>\alpha'_2$
and $p_h(\alpha'_1)>\alpha'_2$. 

Suppose that this is the case.
By Proposition~\ref{p.prepare}~(1) and~(2), it follows that
there exist elements $\alpha_1=(r_1,C_1)\in\monjpr$ and
$\alpha_2=(R_2,c_2)\in\monjpone\textup{(up)}$
with $C_1\leq c_1$ and $r_2\leq R_2$.  Since $p_h(\alpha'_1)>\alpha'_2$,
it follows that $\alpha_1>\alpha_2$.
Now, if $\alpha_2=p_h(\sigma_j)$,   then 
Proposition~\ref{p.sigmaj}~(\ref{i.2.p.sigmaj}) is contradicted;
if $\alpha_2$ belongs to $\monjponepr$,    
Proposition~\ref{p.sigmaj}~(\ref{i.4.p.sigmaj}) is contradicted
(because $p_h(\alpha_1)>\alpha_2$).    
\end{myproof}

\mysection{The map $\ophi$}\mylabel{s.ophi}
The purpose of this section is to describe the map~$\ophi$\index{ophi@$\ophi$} and prove
some basic facts about it.  Certain proofs here
refer to results from~\S\ref{s.lemmas},  but there is
no circularity---to postpone the definition of $\ophi$ until all
the results needed for it have been proved would hurt rather
than help readability. As in~\S\ref{s.opi}, the symbol~$j$ will
be reserved for an odd integer throughout this section.
\mysubsection{Description of $\ophi$}\mylabel{ss.ophi}%
\checked{09 November;  nov-09 version should therefore be fine}
The map $\ophi$ takes as input a
pair~$(w,\mont)$, where~$\mont$%
\index{tmon@$\mont$, fixed monomial in~$\posv$ in~\S\ref{s.ophi}|ignore}
is a monomial, possibly empty,
in~$\posv$ and~$w\geq v$
an element of $\id$ that \ortho-dominates $\mont$,  and produces as output
a monomial~$\mont^*$ of~$\posv$.      To describe~$\ophi$,  we first
partition $\mont$ into subsets $\montwjjpone$\index{twjjpone@$\montwjjpone$}.   As the subscript~$w$ 
in~$\montwjjpone$ suggests, this partition depends on~$w$.


For an odd integer $j$,  let $\mon_w^j$%
\index{Swupj@$\mon_w^j$, $w$ in~$\id$, $j$ odd} 
(respectively $\monwsubjjpone$%
\index{Swsubjjpone@$\monwsubjjpone$, $w$ in~$\id$, $j$ odd})
denote the subset of $\mon_w$
consisting of those elements that are $j$-deep
(respectively that are $j$~deep but not $j+2$~deep, or equivalently
of depth $j$ or $j+1$) in $\mon_w$ in the sense of~\cite[\S4]{kr}.
Since $\mon_w$ is distinguished, symmetric, and has evenly many elements on the
diagonal~$\diagv$,    it follows that~$\monwsupj$ and~$\monwsubjjpone$
too have these properties, and that, in fact, the number of diagonal
elements of~$\monwsubjjpone$ is either~$0$ or $2$ (in the latter case,
the elements have to be distinct since $\mon_w$ is distinguished and
so is multiplicity free).   Let us denote by $\wsupj$\index{wsupj@$\wsupj$, $j$ odd} and $\wsubjjpone$\index{wsubjjpone@$\wsubjjpone$, $j$ odd} the 
elements of $I(d)$ corresponding to $\mon_w^j$ and $\mon_{w,j,j+1}$ 
by Proposition~\ref{p.distid}.

Let $\montwjjpone$\index{twjjpone@$\montwjjpone$}  denote the subset of $\mont$
consisting of those elements $\alpha$ such that
\begin{itemize}
\item
every $v$-chain in
$\mont$ with head $\alpha$ is \ortho-dominated by $\wsupj$,  and
\item there
exists a $v$-chain in $\mont$ with head $\alpha$ that is not \ortho-dominated
by~$\wsupjptwo$.   
\end{itemize}
It is evident that the subsets $\montwjjpone$ are disjoint (as $j$ varies
over the odd integers) and that their union is all of $\mont$ (for $w=w^1$
\ortho-dominates all \vchains\ in~$\mont$ by hypothesis and  
$\mon_w^j$ is empty for large~$j$ and so $w^j=v$).   In other words,
the $\montwjjpone$ form a partition of~$\mont$.
\begin{lemma}\mylabel{l.forphi}
\begin{enumerate}
\item The length of a $v$-chain in $\montwjjpone\cup\montwjjpone^\hash$ is at most~$2$. 
In fact, the \odepth\ of any element in $\montwjjpone$ is at most~$2$.
\item
$\wsubjjpone$ \ortho-dominates $\montwjjpone$.
\end{enumerate}
\end{lemma}
\begin{myproof}
The lemma follows rather easily from Corollary~\ref{c.odomination} as we now
show.   Let $C$ be a \vchain\ in~$\montwjjpone$.    Let $\tau$ be the tail
of $C$.    Choose a \vchain~$D$ in~$\mont$ with head $\tau$ that is not
\ortho-dominated by~$\wsupjptwo$.    Let $E$ be the concatenation of 
$C$ with $D$.     Since the head of $E$ belongs to~$\montwjjpone$,  
it follows that $E$ is \ortho-dominated by $\wsupj$.   It
follows from (the only if part of) Corollary~\ref{c.odomination} (applied
with $\mon=E$ and $x=\wsupj$) 
that $\wsubjjpone$ \ortho-dominates $E_1^\textup{pr}\cup E_2^\textup{pr}$
and $w^{j+2}$ \ortho-dominates $E^{3,\textup{pr}}$.    This means
$\tau\not\in E^{3,\textup{pr}}$,  so $\tau\in E_1^\textup{pr}\cup E_2^\textup{pr}$,
and so $C\subseteq E_1^\textup{pr}\cup E_2^\textup{pr}$.    
This proves~(2).
By Proposition~\ref{p.odepth}~(\ref{p.odepth.ini}),  the \odepths\ of
elements of $C$ are the same in $C$ and $E$,  so $C\subseteq C_1^\textup{pr}
\cup C_2^\textup{pr}$, which proves the second assertion of~(1).
The first assertion of~(1) follows from
the second (see Lemma~\ref{c.od>d}).\end{myproof}
\begin{corollary}\mylabel{c.l.forphi}
$\wsubjjpone$ dominates $\montwjjpone\cup\montwjjpone^\hash$ in the sense of~\cite{kr}.
\end{corollary}
\begin{myproof}
  This follows from~(2) of Lemma~\ref{l.forphi} and Corollary~\ref{c.odtod}
(the latter applied with $\mon=\montwjjpone$ and $w=\wsubjjpone$).
\end{myproof}

\mcomment{From this point on,  the section is almost entirely rewritten.}
We may therefore apply
the map~$\phi$ of~\cite[\S4]{kr} to the pair $(\wsubjjpone,\montwjjponehash)$ to
obtain a monomial~$\montwjjponehashstar$\index{twjjponehashstar@$\montwjjponehashstar$} in~$\andposv$.
In applying~$\phi$, 
there is the partitioning
of $\montwjjponehash$ into ``pieces'',  these being indexed by elements of
$\mon_{\wsubjjpone}=\monwsubjjpone$---observe that the elements of depth~$1$
(respectively~$2$) of $\monwsubjjpone$ are 
precisely those of $\monw$ of depth~$j$
(respectively~$j+1$).   We denote by $\pbeta$\index{pbeta@$\pbeta$}
the piece of $\montwjjponehash$ corresponding to $\beta$ in $\mon_{\wsubjjpone}$.
We also use the notation $\pbstar$\index{pbetastar@$\pbstar$}
as in \cite{kr}.   Moreover,
we will use the phrase {\em piece of\/ $\mont$\/}\index{pieceoft@piece of $\mont$ (see also caution)} (with respect to $w$
being implicitly understood) to refer to a piece of $\montwjjponehash$
for some odd integer~$j$.    
\begin{quote}\begin{footnotesize}
Caution:   Thinking of $\mont$ as a monomial
in~$\andposv$ and $w$ as an element of $I(d,2d)$ that dominates it,
there is, as in~\cite{kr}, the notion
of ``piece of $\mont$'' (with respect to~$w$).     
The two notions of ``piece'' are different.\end{footnotesize}
\end{quote}
\begin{lemma}\mylabel{l.forphi.2}
\begin{enumerate}
\item
The monomial $\montwjjponehashstar$ is symmetric and has either none
or two distinct diagonal elements depending exactly on whether
$\mon_{\wsubjjpone}=\monwsubjjpone$
has~$0$ or~$2$ elements on the diagonal.
\item  The depth of $\montwjjponehashstar$ is~$2$; and
$\cup_{\beta\in(\mon_w)_j}\pbstar$, 
$\cup_{\beta\in(\mon_w)_{j+1}}\pbstar$ are respectively the elements of depth~$1$
and~$2$ in $\montwjjponehashstar$.
\end{enumerate}
\end{lemma}
\begin{myproof}
(1)~The symmetry follows by combining Proposition~5.6
of~\cite{gr}, which says that the map $\pi$ respects the involution~$\hash$, 
with Proposition~4.2 of~\cite{kr},  which says that $\pi$ and $\phi$ are
are inverses of each other.   

The assertion about diagonal elements follows
by combining item~(B) of~\cite[Proposition~5.10]{gr},  which is an assertion
about the existence and relative multiplicities of diagonal elements 
in $\block$ and $\block'$ where $\block$ is a diagonal block of a monomial
in $\andposv$,  and Proposition~4.2 of \cite{kr}.

(2)~It follows from Propositions~4.2 of~\cite{kr} that the map $\pi$
(described in \S4 of that paper) applied to $\montwjjponehashstar$ 
results in the pair $(\wsubjjpone,\montwjjpone\cup\montwjjpone)$.
It now follows from Lemma~4.16 of~\cite{kr} that the depth of 
$\montwjjponehashstar$ is exactly~$2$.    The latter assertions again
follow from the results of~\cite{kr}---in fact, the proof that
$\pi\circ\phi$ is identity on pages~47--49 of~\cite{kr} shows that
the $\pbstar$ are
the blocks in the sense of~\cite{kr} of the monomial $\montwjjponehashstar$.
\end{myproof}

Suppose that $\montwjjponehashstar$ contains the pair
$(a,a^*)$, $(b,b^*)$ of diagonal elements with $a>b$.
We call the pair $(b,a^*)$, $(a,b^*)$
the ``twists,'' and set $\delta_j:=(b,a^*)$.
In other words, $\delta_j$\index{deltaj@$\delta_j$, for $j$ odd} is the element of the twisted pair
that lies above the diagonal---observe that the twisted elements
are reflections of each other.   We allow ourselves the
following ways of expressing the condition
that $\montwjjponehashstar$
has diagonal elements: {\em $\deltaj$ exists;   $w$ is diagonal
at $j$} (the latter expression is justified by the lemma above).

With notation as above,   consider the new monomial defined as 
\[\left\{
        \begin{array}{cl}
          \montwjjponehashstar & \textup{if $w$ is not diagonal at $j$}\\
 \left(\montwjjponehashstar\setminus\diagv\right)\cup\{\delta_j,\delta_j^\hash\}
       &   \textup{if $w$ is diagonal at $j$}\\
                  \end{array}
                                       \right.\]
This new monomial is 
symmetric and contains no diagonal elements.   Its intersection 
with~$\posv$ is denoted $\montwjjponestar$\index{twjjponestar@$\montwjjponestar$}.   In other words,~$\montwjjponestar$
is the intersection of the new monomial with the subset of~$\andposv$
of those elements that lie strictly above the diagonal.

The union of~$\montwjjponestar$ over all
odd integers $j$ is defined to be $\montwstar$\index{twstar@$\montwstar(:=\ophi(w,\mont))$}, the result of~$\ophi$ applied to $(w,\mont)$.    This finishes the description of the map~$\ophi$.

For $\beta$ in $\monwsubjjpone\up$,
we define the ``orthogonal piece-star'' $\opbstar$
\index{opbetastar@$\opbstar$}
corresponding to $\beta$ as
\begin{equation}\label{n.opstar}
\opbstar:=\left\{
  \begin{array}{ll}
    \pbstar=\pbstar\up & \textup{if $\beta$ is not on the diagonal}\\
    \pbstar\cap\posv & \textup{if $\beta\in(\mon_w)_{j+1}$ is on 
      the diagonal}\\
    \{\pbstar\cap\posv\}\cup\{\deltaj\} &
      \textup{if $\beta\in(\mon_w)_{j}$ is on 
      the diagonal}\\
  \end{array}\right.
\end{equation}
With this, we can say that $\montwstar$ is the union of
$\opbstar$ as $\beta$ varies over $\monw\up$.
\begin{lemma}\mylabel{l.forphi.3}
Suppose that $\montwjjponehashstar$ contains the pair $(a,a^*)$, $(b,b^*)$
of diagonal elements with $a>b$.     Let 
\[ \ldots,~(r_1,c_1),~(a,a^*),~(c_1^*,r_1^*),~\ldots; \quad\quad
\ldots,~(r_2,c_2),~(b,b^*),~(c_2^*,r_2^*),~\ldots\]
be respectively the elements of depth~$1$ and~$2$ of $\montwjjponehashstar$
arranged in increasing order of row and column indices.
Then 
\begin{enumerate}
\item $c_1\leq a^*$ and $r_1\leq b$ (assuming $(r_1,c_1)$ exists); and
\item $r_2<b$ and $c_2\leq b^*$ (assuming $(r_2,c_2)$ exists).
\end{enumerate}
\end{lemma}
\begin{myproof}
  (1)~ Suppose that $(r_1,c_1)$ exists.
It is clear that $c_1\leq a^*$.   
From way the map $\phi$ of~\cite{kr} is defined, it follows
that $(r_1,a^*)$ is an element of $\montwjjpone$.  Suppose that $r_1>b$.
Then $p_h(r_1,a^*)=(r_1,r_1^*)$ belongs to~$\andposv$.
We consider two cases.   

 If $(r_2,c_2)$ exists,  then,
again from the definition of the map~$\phi$, it follows that $(r_2,b^*)$
is an element of $\montwjjpone$.   But then $p_h(r_1,a^*)=(r_1,r_1^*)>(b,b^*)$
and $(b,b^*)$ dominates $(r_2,b^*)$,  which means that the \vchain\
$(r_1,a^*)>(r_2,b^*)$ (note that $a^*<b^*$ because $a>b$ by hypothesis)
in $\montwjjpone$ has \odepth\ more than~$2$,  a contradiction to
Lemma~\ref{l.forphi}~(1).

Now suppose that $(r_2,c_2)$ does not exist.
(Then~$(b,b^*)$ is the diagonal element in~$(\monw)_{j+1}$.)
Consider the singleton
\vchain~$C:=\{(r_1,a^*)\}$ in $\montwjjpone$.   
Then $\mon_C=\{(a,a^*),(r_1,r_1^*)\}$ which is not dominated 
by $\wsubjjpone$, a contradiction to Lemma~\ref{l.forphi}~(2).

(2)~ Suppose that $(r_2,c_2)$ exists.  Then there exists, by the
definition of the map~$\phi$, an element $(r_2,b^*)$ in $\montwjjpone$.
Since $(r_2,b^*)$ lies above the diagonal, it follows that $r_2<b$.
That $c_2\leq b^*$ is clear.
\end{myproof}
\mysubsection{Basic facts about $\montwjjpone$ and $\montwjjponestar$}%
\mylabel{ss.phibasic}
\begin{lemma}\mylabel{l.jjp}
  \begin{enumerate}
  \item 
Let $\alpha'>\alpha$ be elements of $\mont$.   Let $j$ and $j'$ be
the odd integers such that $\alpha'\in\mont_{w,j',j'+1}$ and
$\alpha\in\montwjjpone$.   Then $j'\leq j$.
  \item
If, further, either 
\begin{enumerate}
\item there exists $\mu$ in $\mont$ such that $\alpha'>\mu>\alpha$, or
\item $\alpha'\in\pbetap$ for $\beta'$ in $(\mon_w)_{j'+1}$,
\end{enumerate}
then $j'<j$.
  \end{enumerate}
\end{lemma}
\begin{myproof}
(1)  By hypothesis, every \vchain\ with head $\alpha'$ is \ortho-dominated
by~$w^{j'}$.    This implies, by Corollary~\ref{c.vchain1}, that every
\vchain\ with head $\alpha$ is \ortho-dominated by~$w^{j'}$.    This
shows $j'\leq j$.

(2a) Suppose that $j'=j$.   It follows from~(1) that $\alpha'$, $\mu$,
and $\alpha$ all belong to $\montwjjpone$.   But then $\alpha'>\mu>\alpha$
is a \vchain\ of length~$3$ in $\montwjjpone$,  a contradiction to
Lemma~\ref{l.forphi}~(1).

(2b) Suppose that $j'=j$.   Then $\alpha'>\alpha$ is a \vchain\ in
$\montwjjpone$.   Being of length~$2$,  it cannot be dominated by
$(\mon_w)_{j+1}$,  which means, by the definition of $\pbetap$, that 
$\alpha'$ cannot belong to $\pbetap$,  a contradiction.
\end{myproof}
\begin{proposition}\mylabel{p.t.od<2}
  \begin{enumerate}
   \item The length of a \vchain\ in~$\montwjjponestar$ is at most~$2$.
  \item The \odepth~of $\montwjjponestar$ is at most~$2$.
\item
$\cup_{\beta\in(\monw)_j\up}\opbstar$ is precisely the set of depth~$1$
elements of $\montwjjponestar$ (in particular,  no two elements there
are comparable); if $\delta_j$ exists, then it is the last element of 
$\cup_{\beta\in(\monw)_j\up}\opbstar$ when
the elements are arranged in increasing order of row and column indices.
\item
$\cup_{\beta\in(\monw)_{j+1}\up}\opbstar$ is precisely the set of
depth~$2$ elements of $\montwjjponestar$ (in particular,  
no two elements there are comparable); if $\delta_j$ exists, 
then its row index exceeds the row index of any element 
in~$\cup_{\beta\in(\monw)_{j+1}\up}\opbstar$.
\end{enumerate}
\end{proposition}
\begin{myproof}  
For~(1), it is enough,  given Lemma~\ref{l.forphi.2}~(2),  to show that
$\deltaj$ is not comparable to any element of depth~$1$ 
of~$\montwjjponehashstar$,  and this follows from Lemma~\ref{l.forphi.3}~(1).
In fact, the above argument proves also~(3).  

For~(4), it is enough,  given Lemma~\ref{l.forphi.2}~(2),
the symmetry of the monomials involved in that lemma, and
the observation that $\alpha>\beta$ implies~$\alpha\up>\beta\up$ for 
elements $\alpha$, $\beta$ of $\andposv$,  to show the following:
if $(a,a^*)>\gamma=(e,f)$ for~$\gamma$ an element of~$\montwjjponehashstar$
lying (strictly) above the diagonal,
then $\deltaj>\gamma$.    But this follows from Lemma~\ref{l.forphi.3}~(2):
$\gamma$ is a depth~$2$ element in $\montwjjponehashstar$, and we have
$e\leq r_2<b$ (and $a^*<f$ since $(a,a^*)>\gamma$).  In fact, the above argument proves also~(2): observe that $f\leq b^*$ (Lemma~\ref{l.forphi.3}~(2)).~\end{myproof}

\mysection{Some Lemmas}\mylabel{s.lemmas}
The main combinatorial results of this paper are~Propositions~\ref{p.4.1.kr}
and~\ref{p.4.2.kr}.   They are analogues respectively of 
Propositions~4.1 and~4.2 of~\cite{kr}.    We have tried 
to preserve the structure of the proofs in~\cite{kr}
of those propositions.   The proofs in~\cite{kr} rely on certain
lemmas and it is natural therefore to first establish
the orthogonal analogues of those.    The purpose of this
section is precisely that.
Needless to say that the lemmas (especially those in \S\ref{ss.pf.lemmas}) 
may be unintelligible until one tries to read~\S\ref{s.proof}.

The division of this section into four subsections is also suggested
by the structure of the proofs in~\cite{kr}.   
Each subsection 
has at its beginning a brief description of its contents.

\mysubsection{Lemmas from the Grassmannian case}\label{ss.l.grass}
In this subsection,  the terminology and notation of~\cite[\S4]{kr}
are in force.   The statements here could have been made in~\cite[\S4]{kr}
and would perhaps have improved the efficiency of the proofs there,
but do not appear there explicitly.

Let $\mon$ 
be a monomial in $\andposv$.
Recall from~\cite{kr} the notion of {\em depth\/}%
\index{depthelement@depth (of an element $\alpha$ in a monomial $\mon$ in~$\andposv$) = $\depthin{\mon}{\alpha}$} 
of an element~$\alpha$ in~$\mon$:  
it is the largest possible length of a
\vchain\ in $\mon$ with tail~$\alpha$ and denoted $\depthin{\mon}{\alpha}$%
\index{depthelement@depth (of an element $\alpha$ in a monomial $\mon$ in~$\andposv$) = $\depthin{\mon}{\alpha}$}.
The {\em depth\/}%
\index{depthmonomial@depth (of a monomial $\mon$ in $\andposv$)}
of $\mon$ is the maximum of
the depths in it of all its elements.   We denote by $\mon_k$%
\index{Sk@$\mon_k$, for monomial $\mon$ in $\andposv$} the set
of elements of depth~$k$ of~$\mon$ (as in ~\cite{kr}) and by $\mon^k$%
\index{Supk@$\mon^k$, for monomial $\mon$ in $\andposv$} the set of elements
of depth at least~$k$ of $\mon$.
\begin{quote}
{\footnotesize
Caution:    For a monomial $\mon$ of $\posv$, we have introduced 
in~\S\ref{ss.opi} the notation $\mon_k$. 
That is different from the~$\mon_k$ we have just defined.
}
\end{quote}

\begin{lemma}\mylabel{l.d.monmon'}
Let $\mon$ be a monomial in $\andposv$, and let $\pi(\mon)=(w,\mon')$,
where~$\pi$ is the map defined in~\cite[\S4]{kr}.   Then the maximum 
length of a \vchain\ in $\mon\cup\mon'$ is the same as the maximum length
of a \vchain\ in $\mon$.
\end{lemma}
\begin{myproof}
We use the notation of~\cite[\S4]{kr} freely.    Let~$d$ be the maximum
length of a \vchain\ in $\mon$.    Suppose~$\alpha_1>\ldots>\alpha_\ell$
is a \vchain\ in $\mon\cup\mon'$.    Let~$i_1,\ldots,i_\ell$ be such 
that~$\alpha_j$ belongs to~$\mon_{i_j}\cup\mon_{i_j}'$ (the integers $i_j$ are
uniquely determined---see~Corollary~5.4 of~\cite{gr}).  We
claim that $i_1<\ldots<i_\ell$.   This suffices to prove the lemma, for
$\mon_k\cup\mon_k'$ is empty for $k>d$.

To prove the claim, it is enough to show $i_1<i_2$.   
It follows from Lemma~4.10 of~\cite{kr} that $i_1\neq i_2$.  
We now assume that~$i_1>i_2$ and arrive at a contradiction. 
First suppose that $\alpha_1\in\mon_{i_1}$.    
Then,  by the definition of $\mon_{i_1}$,
there exists $\beta$ in~$\mon_{i_2}$ with $\beta>\alpha_1$.   Now
$\beta>\alpha_2$ and both~$\beta$, $\alpha_2$ belong 
to~$\mon_{i_2}\cup\mon'_{i_2}$, a contradiction to~\cite[Lemma~4.10]{kr}.
If $\alpha_1=(r,c)$ belongs to $\mon_{i_1}'$,  then,
by the definition of $\mon_{i_1}'$,  there exists $(r,a)$ in $\mon_{i_1}$
with $a\leq c$,  and there exists $\beta$ in $\mon_{i_2}$ with $\beta>(r,a)$.
This leads to the same contradiction as before.
\end{myproof}
\newcommand\montemp{\blockb}
\begin{lemma}\mylabel{l.kr.!}
Let $\montemp$ and $\monu$ be monomials in~$\andposv$.
Assume that
\begin{itemize}
\item
the elements of $\montemp$ form
a single block\index{block!in the sense of ~\cite{kr}} (in the sense of \cite[Page~38]{kr}).
\item 
$\monu$ has depth~$1$ (equivalently, 
there are no comparable elements in~$\monu$).
\item for every $\beta=(r,c)$ in $\blockb$, there exist
$\gamma^1(\beta)=(R^1,C^1)$, and $\gamma^2(\beta)=(R^2,C^2)$
in $\monu$ such that
\[
C^1<c, \quad C^2<R^1, \quad r<R^2	\]
(this holds, for example, when there exists $\gamma(\beta)$ in
$\monu$ such that $\gamma(\beta)>\beta$: take $\gamma^1(\beta)=\gamma^2(\beta)
=\gamma(\beta)$).
\end{itemize}
Then there exists a unique block\index{block!in the sense of ~\cite{kr}}~$\blockc$ of $\monu$ such that
$w(\blockc)>w(\montemp)$.
\end{lemma}
\begin{myproof}
It is useful to isolate the following observation:
\newcommand\rone{r_1}
\newcommand\rtwo{r_2}
\newcommand\cone{c_1}
\newcommand\ctwo{c_2}
\newcommand\goneone{\gamma_1^1}
\newcommand\gonetwo{\gamma_1^2}
\newcommand\gtwoone{\gamma_2^1}
\newcommand\gtwotwo{\gamma_2^2}
\newcommand\coneone{C_1^1}
\newcommand\conetwo{C_1^2}
\newcommand\ctwoone{C_2^1}
\newcommand\ctwotwo{C_2^2}
\newcommand\roneone{R_1^1}
\newcommand\ronetwo{R_1^2}
\newcommand\rtwoone{R_2^1}
\newcommand\rtwotwo{R_2^2}
\begin{lemma}\mylabel{l.kr.!.1}
Let $(\rone,\cone)$ and $(\rtwo,\ctwo)$ be elements of~$\andposv$ with
$\ctwo<\rone\leq \rtwo$.   Let $\goneone=(\roneone,\coneone)$,~$\gonetwo=(\ronetwo,\conetwo)$ and
$\gtwoone=(\rtwoone,\ctwoone)$,~$\gtwotwo=(\rtwotwo,\ctwotwo)$ be elements of~$\andposv$ such that
\begin{enumerate}
\item\mylabel{i.1.l.kr.!.1}
$\coneone\leq\cone, \quad \conetwo<\roneone, \quad \rone\leq\ronetwo$.
\item\mylabel{i.2.l.kr.!.1}
$\ctwoone\leq\ctwo, \quad \ctwotwo<\rtwoone, \quad \rtwo\leq\rtwotwo$.
\item\mylabel{i.3.l.kr.!.1}
No two of $\goneone$, $\gonetwo$, $\gtwoone$, $\gtwotwo$ are comparable 
(they could well be equal and this is important for us---see our
definition of comparability).
\end{enumerate}
Then the monomial $\{\goneone,\gonetwo,\gtwoone,\gtwotwo\}$
consists of a single block.
\end{lemma}
\begin{myproofnobox}
It follows from assumption~(\ref{i.1.l.kr.!.1}) that $\goneone$ and~$\gonetwo$
belong to a single block:
\begin{itemize}
\item if $\roneone<\ronetwo$,  then $\conetwo<\roneone$ becomes relevant;
\item if $\ronetwo<\roneone$,  then the other two inequalities
in~(\ref{i.1.l.kr.!.1}) become relevant:
\newline
$\coneone\leq\cone<\rone\leq\ronetwo$.
\end{itemize}
Similarly
it follows from assumption~(\ref{i.2.l.kr.!.1}) that $\gtwoone$ and~$\gtwotwo$
belong to a single block.

We therefore need only consider the cases when, in the arrangement of the
elements $\{\goneone,\gonetwo,\gtwoone,\gtwotwo\}$ in increasing order of
row indices,   both $\goneone$,~$\gonetwo$ come before or after
$\gtwoone$,~$\gtwotwo$.      In the former case,  the first sequence of
inequalities below shows that $\gonetwo$ and $\gtwoone$ belong to the
same block, and we are done;   in the latter case,  the second sequence
of inequalities below shows that $\gtwotwo$ and~$\goneone$ belong to
the same block, and we are done:
\begin{itemize}
\item
$\ctwoone\leq\ctwo<\rone\leq\ronetwo$.
\item
$\coneone\leq c_1<r_1\leq r_2\leq\rtwotwo$.\hfill$\Box$
\end{itemize}
\end{myproofnobox}

Continuing with the proof of Lemma~\ref{l.kr.!},
we first prove the existence part.
Arrange the elements of $\montemp$ in non-decreasing order of row
numbers as well as column numbers (this is possible since there are
no comparable elements in $\montemp$).    If $\beta_1=(r_1,c_1)$ and
$\beta_2=(r_2,c_2)$ are successive elements,   then
$c_2<r_1\leq r_2$ (since $\montemp$ is a single
block).   
Apply Lemma~\ref{l.kr.!.1} with
$\goneone=\gamma^1(\beta_1)$,
$\gonetwo=\gamma^2(\beta_1)$,  and
$\gtwoone=\gamma^1(\beta_2)$,
$\gtwotwo=\gamma^1(\beta_2)$.
We conclude that $\{\goneone,\gonetwo,\gtwoone,\gtwotwo\}$ 
belongs to a single block, say~$\blockc$, of $\monu$.     Continuing
thus, we conclude that all $\gamma^1(\beta)$ and $\gamma^2(\beta)$,
as $\beta$ varies over
$\montemp$, belong to~$\blockc$.      Since the row (respectively column)
index of $w(\blockc)$ is the maximum (respectively minimum) of all
row (respectively column) indices of elements of $\blockc$ (and similarly
for $\montemp$),   it follows that $w(\blockc)>w(\montemp)$.  

To prove uniqueness,  let $\blockc_1$ and $\blockc_2$ be two blocks of
$\monu$ with $w(\blockc_1)>w(\montemp)$ and $w(\blockc_2)>w(\montemp)$.
Apply the lemma with $(r_1,c_1)=(r_2,c_2)=w(\montemp)$ and 
$\goneone=\gonetwo=w(\blockc_1)$ and $\gtwoone=\gtwotwo=w(\blockc_2)$;  
it follows from~\cite[Lemma~4.9]{kr} that $w(\blockc_1)$ and~$w(\blockc_2)$
are not comparable.     
But, unless $\blockc_1=\blockc_2$, neither is the monomial
$\{w(\blockc_1),w(\blockc_2)\}$ a single block,  again by~\cite[Lemma~4.9]{kr}.
\end{myproof}
\begin{lemma}\mylabel{l.domination}
Let $\mon$ be a monomial in $\andposv$ and $x$ an element of $I(d,n)$.
For $x$ to dominate $\mon$ it is necessary and sufficient that for every
$\alpha=(r,c)$ in $\mon$ there exist $\beta=(R,C)$ in $\mon_x$ with
$C\leq c$, $r\leq R$, and $\depthin{\mon_x}{\beta}\geq\depthin{\mon}{\alpha}$.
(Here $\mon_x$
\index{Sw@$\monw$, $w$ in $\idd$ or $I(d,n)$}
denotes the distinguished  monomial in $\andposv$ associated
to $x$ as in \cite[Proposition~4.3]{kr}.)
\end{lemma}
\begin{myproof} 
The lemma is a corollary of \cite[Lemma~4.5]{kr} as we now show.

First suppose that $x$ dominates $\mon$.  Let $\alpha=(r,c)$
be an element of $\mon$, and
$C$ a \vchain\ in $\mon$ with tail $\alpha$ and length $\depthin{\mon}{\alpha}$.
Since $x$ dominates $C$,   there exists, by \cite[Lemma~4.5]{kr},
a chain in $D$ in $\mon_x$ of length $\depthin{\mon}{\alpha}$ and tail
$\beta=(R,C)$ with $C\leq c$ and $r\leq R$, and we are done with the
proof of the necessity.

To prove the sufficiency,  let $C:\alpha_1=(r_1,c_1)>\ldots>\alpha_k=(r_k,c_k)$
be a \vchain\ in $\mon$.    By hypothesis,  there exist $\beta_1=(R_1,C_1)$,
\ldots, $\beta_k=(R_k,C_k)$ in $\mon_x$ with $C_i\leq c_i$, $r_i\leq R_i$,
and $\depthin{\mon_x}{\beta_i}=i$ for $1\leq i\leq k$  (observe that
replacing the $\geq$ in the latter condition of the statement by an
equality yields an equivalent statement).    We claim that $\beta_1>\ldots
>\beta_k$.    By~\cite[Lemma~4.5]{kr}, it suffices to prove the claim.

Since $\beta_k$ has depth $k$ in $\mon_x$,  there exists a
$\beta'_{k-1}=(R'_{k-1},C'_{k-1})$
of depth $k-1$ in $\mon_x$ such that $\beta'_{k-1}>\beta_k$.  It follows
from the distinguishedness of $\mon_x$ that 
that $\beta'_{k-1}=\beta_{k-1}$:   if not,   then we have two
distinct elements of the same depth (namely $k-1$) in $\mon_x$ both
dominating $\alpha_k$,  a contradiction.   So $\beta_{k-1}>\beta_k$,
and the claim is proved by continuing in a similar fashion.
\end{myproof}

Let $x$ be an element of $I(d,n)$. Let $\mon_x$ denote the distinguished  monomial in $\andposv$ associated
to $x$ as in \cite[Proposition~4.3]{kr}.
For $k$ a positive integer,
let $x_k$%
\index{xkkpone@$x_k$, $x^k$, $x_{k,k+1}$, for $x\in I(d,n)$|(}
denote the element of $I(d,n)$ corresponding
to the distinguished subset $(\mon_x)_k$. For a monomial $\mon$ of $\andposv$, let $\mon_{k,k+1}:=
\mon_k\cup\mon_{k+1}$%
\index{Skkpone@$\mon_{k,k+1}$, for monomial $\mon$ in $\andposv$}. 
Let $x_{k,k+1}$ denote the element of $I(d,n)$ corresponding
to the distinguished monomial $(\mon_x)_{k,k+1}$;
let $x^k$ denote the element of $I(d,n)$ corresponding
to the distinguished subset $(\mon_x)^k$.
\index{xkkpone@$x_k$, $x^k$, $x_{k,k+1}$, for $x\in I(d,n)$|)}
\begin{quote}
{\footnotesize
Caution:    For a monomial $\mon$ of $\posv$ and an odd integer $j$, we have
introduced in~\S\ref{ss.opi} the notation $\mon_{j,j+1}$. 
That is different from the~$\mon_{k,k+1}$ just defined.
}
\end{quote}

\begin{corollary}\mylabel{c.l.domination}
$x$ dominates $\mon$  $\Leftrightarrow$ 
		$x_k$ dominates $\mon_k$ $\forall$ $k$
	 $\Leftrightarrow$  $x_{1,2}$ dominates $\mon_{1,2}$ and 
		$x^{3}$ dominates $\mon^3$.
\end{corollary}
\begin{myproof}  The first equivalence is a restatement of
the lemma:   in the statement of the lemma we could equally well
have written $\depthin{\mon_x}{\beta}=\depthin{\mon}{\alpha}$.
The second follows from the first and the following observations:
$(\mon_{1,2})_1=\mon_1$, 
$(\mon_{1,2})_2=\mon_2$, 
$(\mon^{3})_k=\mon_{k+2}$;   and
$(x_{1,2})_1=x_1$, 
$(x_{1,2})_2=x_2$,  
$(x^{3})_k=x_{k+2}$.
\end{myproof}
\mysubsection{Orthogonal analogues of Lemmas of~\ref{ss.l.grass}}%
\mylabel{ss.oa.l.grass}  
Lemma~\ref{l.odomination} below is the orthogonal analogue of
Lemma~\ref{l.domination} (more precisely, that of the first assertion
of Corollary~\ref{c.l.domination}).  The following proposition will be
used in its proof.  
\begin{proposition}\mylabel{p.l.odomination}
Let $x$ be an element of $I(d)$ and $\mon$ a monomial in $\posv$.
Then $x$ \ortho-dominates $\mon_1^\pr\cup\mon_2^\pr$ if and only
if it \ortho-dominates every \vchain\ in~$\mon$ of \ortho-depth at most~$2$.
\end{proposition}
\begin{myproof}
The ``if'' part is immediate from definitions (in any case, see
also Proposition~\ref{p.od.2}).    For the ``only if'' part,  let 
$C$ be a \vchain\ in~$\mon$ of \ortho-depth at most~$2$.  
Our goal is to show that $x$ dominates $\mon_C$.
For this, it is enough, by Corollary~\ref{c.l.domination}, to show
that $x_1$ dominates $(\mon_C)_1$ and $x_2$ dominates $(\mon_C)_2$
(by choice of~$C$, $(\mon_C)_k$ is empty for $k\geq 3$).  

Let $\alpha'\in(\mon_C)_1$.   Choose $\alpha$ in $C$ such that
$\alpha'\in\mon_{C,\alpha}$.    Choose $\alpha_0$ in $\mon_1^\pr$
such that $\alpha_0$ dominates $\alpha$.    
Since $x$ \ortho-dominates the singleton \vchain~$\{\alpha_0\}$,
it follows that $x_1$ dominates $q_{\{\alpha_0\},\alpha_0}$.
We claim that $q_{\{\alpha_0\},\alpha_0}$ dominates $\alpha'$.
To prove the claim, we need only rule out the possibility that
$\alpha_0$ is of type~S in $\{\alpha_0\}$ and
$\alpha$ of type~V in $C$.  
Since $\alpha'\in(\mon_C)_1$,  it follows from Proposition~\ref{p.type.>}~(1)
that $\alpha$ is the first element of~$C$.   In particular, if $\alpha$
is of type~V in $C$,  then $p_h(\alpha)\in\andposv$,  so $p_h(\alpha_0)\in
\andposv$,  and $\alpha_0$ is of type~H in $\{\alpha_0\}$.  The claim
is thus proved.

Now consider an element of $(\mon_C)_2$.  Observe that the length
of~$C$ is at most~$2$ (Lemma~\ref{c.od>d}).    So our element is
either the horizontal projection $p_h(\alpha)$ of the head $\alpha$ of~$C$, 
or it is 
$q_{C,\beta}$ where $\beta$ is the tail of~$C$.     In the first case,
let $\alpha_0$ be as  in the previous paragraph, and proceed similarly.
It is clear that $p_h(\alpha_0)\in\andposv$ (because $p_h(\alpha)\in
\andposv$);  $x_2$ dominates $p_h(\alpha_0)$ and so also $p_h(\alpha)$.

Now we handle the second case.   If $\beta\in\mon_2^\pr$,  then $C$ is
contained in $\mon_1^\pr\cup\mon_2^\pr$ and there is nothing to prove.
So assume that $\odepthin{\mon}{\beta}\geq 3$.  Choose a \vchain~$D$
in $\mon$ with tail $\beta$, $\odepthin{D}{\beta}\geq 3$,  and with
the good property as in Proposition~\ref{p.goodCexists}.   There
occurs in $D$ an element of \ortho-depth~$3$,  say $\delta$.
(Lemma~\ref{l.odepths}~(\ref{l.i.hits})).    Let $A$ denote the part
$\delta>\ldots$ of~$D$ and $C'$ the part up to but not including $\delta$.
There clearly is an element---call it~$\mu$---of depth~$2$ in~$\mon_D$
that dominates $q_{D,\beta}$.  This element $\mu$ belongs to $\mon_{C'}$
(Corollary~\ref{c.for.odom}~(3)).     Since $D$ has the good property
of Proposition~\ref{p.goodCexists},   $C'\subseteq\mon_1^\pr\cup
\mon_2^\pr$,   so $\mu$ is dominated by an element in $(\mon_x)_2$.
In particular, $q_{D,\beta}$ is dominated by the same element of $(\mon_x)_2$.

We are still not done, for it is possible that $q_{D,\beta}$ be $\beta$
and $q_{C,\beta}$ be $p_v(\beta)$.   Suppose that this is the case.
Then $\alpha>\beta$ is connected.   So
$p_h(\alpha)\in\andposv$ and the legs of $\alpha$ and $\beta$ intertwine.
As seen above in the third paragraph of the present proof,  there is
an element of $(\mon_x)_2$ that dominates $p_h(\alpha)$.    By the
distinguishedness of $\mon_x$,  it follows that the element in $(\mon_x)_2$
dominating $\beta$ is the same as the one dominating~$p_h(\alpha)$.  By the
symmetry of $\mon_x$,   this element lies on the diagonal and so dominates
$p_v(\beta)$,  and, finally, we are done with the proof in the second case.
\end{myproof}
\begin{lemma}\mylabel{l.odomination}
Let $\mon$ be a monomial in $\posv$ and $x$ an element of $I(d)$.
For $x$ to \ortho-dominate $\mon$ it is necessary and sufficient that, for
every odd integer~$j$, every \vchain\ in $\monjonlypr\cup\monjponepr$
is \ortho-dominated by $x_{j,j+1}$.
\end{lemma}
\begin{myproof}
First suppose that $x$ dominates $\mon$. 
Let $j$ be an odd integer and let $A$ a \vchain\ in 
$\monjonlypr\cup\monjponepr$.    We need to show that $\xjjpone$ dominates
$\mon_A$.  
For this, we may assume that $A$ is maximal (by
Corollary~\ref{c.vchain1}).    By Corollary~\ref{c.p.odepth}~(3),
the length of $A$ is at most~$2$.  By
Lemma~\ref{l.odepths}~(\ref{l.i.hits})~(b),  for every $\beta$ in
$\monjponepr$ there exists $\alpha$ in $\monjonlypr$ with $\alpha>\beta$.
Thus we may assume that the head~$\alpha$ of $A$ belongs to~$\monjonlypr$.

It is enough to show (see \cite[Lemma~4.5]{kr})
that for any \vchain~$E$ in $\mon_A$
\begin{itemize}
\item
the length of~$E$ is at most~$2$;
\item 
there exists an $x$-dominated monomial in~$\andposv$ containing~$E$
and the head of~$E$ is an element of depth at least~$j$ in that monomial.
\end{itemize}
The first of these conditions holds by Proposition~\ref{p.od.2}.
We now show that the second holds.

We may assume that $E$ is maximal in $\mon_A$.
By Proposition~\ref{p.type.>}~(1),  the head of $E$ is $q_{A,\alpha}$.
Let $C$ a \vchain\ in $\mon$ with tail $\alpha$ such that
$\odepthinC{\alpha}=j$.   Let $D$ be the concatenation of $C$ with $A$.
We claim that the monomial $\mon_D$ has the desired properties.
That $\mon_D$ is $x$-dominated is clear (since
$x$ \ortho-dominates $\mon$).
By Corollary~\ref{c.for.odom},  it follows that $q_{D,\alpha}=q_{A,\alpha}$
and $\mon_A\subseteq \mon_D$ (in particular that $E\subseteq\mon_D$). 
By Proposition~\ref{p.odepth}~(2),  $\odepthin{D}{\alpha}=\odepthin{C}{\alpha}=j$, that is, $\depthin{\mon_D}{q_{D,\alpha}}=j$.  
The proof of the necessity is thus complete.


To prove the sufficiency, proceed by induction on the largest odd integer~$J$
such that $\mon_J^\pr\cup\mon_{J+1}^\pr$ is non-empty.
When $J=1$,   there is nothing to prove, for $\mon_1^\pr\cup
\mon_2^\pr=\mon$ and $x_{1,2}$ \ortho-dominates $\mon_1^\pr\cup
\mon_2^\pr$.    So suppose that $J\geq 3$.   We implicitly use
Corollary~\ref{c.p.od.shift} in what follows.
By induction, $x^{3}$ \ortho-dominates~$\monsupthreefour$.

Let $D$ be a \vchain\ in $\mon$.   Our goal is to show that 
$x$ dominates $\mon_D$.   Let $\alpha$ be the element of~$D$ with
$\odepthin{D}{\alpha}=3$---such an element exists, by
Lemma~\ref{l.odepths}~(\ref{l.i.hits}) (if there exists in $D$
an element of \ortho-depth in~$D$ exceeding~$2$);   the following
proof works also in the case when $\alpha$ does not exist.
Let $A$ be the part $\alpha>\ldots$ of $D$,  and~$C'$ the part
up to but not including $\alpha$.  
By Proposition~\ref{p.odepth}~(2),
the \ortho-depth (in $C'$) of elements of $C'$ is at most~$2$.   
By Proposition~\ref{p.l.odomination},  $x_{1,2}$ dominates $\mon_{C'}$.
By Corollary~\ref{c.for.odom}~(3),  $(\mon_D)_{1,2}=\mon_{C'}$
and $(\mon_D)^3=\mon_A$.     Since $A\subseteq\monsupthreefour$,
it follows that $x^{3}$ dominates $\mon_A$ (induction hypothesis).   
Finally, by an application of Corollary~\ref{c.l.domination},
we conclude that $x$ dominates $\mon_D$.
\end{myproof}
\begin{corollary}\mylabel{c.odomination}
Let $\mon$ be a a monomial in $\posv$ and $x$ an element of $I(d)$.
For $x$ to \ortho-dominate $\mon$ it is necessary and sufficient
that  $x_{1,2}$ \ortho-dominate $\mon_1^\pr\cup\mon_2^\pr$
and $x^3$ \ortho-dominate $\mon^{3,4}$.
\end{corollary}
\begin{myproof}
It is easy to see that $(x^3)_{j,j+1}=x_{j+2,j+3}$;  it follows from
Proposition~\ref{p.od.shift} that $(\mon^{3,4})^\pr_j\cup(\mon^{3,4})^\pr_{j+1}=
\mon^\pr_{j+2}\cup\mon^\pr_{j+3}$.  The assertion follows from the lemma.  
\end{myproof}
\mysubsection{Orthogonal analogues of some lemmas in~\cite{kr}}%
\mylabel{ss.l.ortho}
The proofs of Propositions~4.1 and~4.2
of~\cite{kr} are based on assertion~4.9--4.16 (of that paper).
Assertion~4.9 being a statement about a single $\mon_k$,
it is applicable in the present situation.   Since references to it
are frequent,  we recall it below as Lemma~\ref{l.4.9.kr}.
As to assertions~4.10--4.16
of~\cite{kr}, assertions~\ref{l.4.10.kr},~\ref{l.4.11.kr}--\ref{c.4.16.kr} 
below are their respective analogues.


A {\em block\/}\index{block!of a monomial $\mon$ in $\posv$} of a monomial $\mon$ in~$\posv$
means a block of $\monjjpone$ in the sense of~\cite{kr} for some
odd integer~$j$.
\begin{quote}{\footnotesize
Caution: Considering $\mon$ as a monomial in $\andposv$,  there is the
notion of a ``block'' of~$\mon$ as in~\cite{kr}, which has in fact
been used in \S\ref{ss.l.grass}, and which is different from the
notion just defined.   Both notions are used and it will be
clear from the context which is meant.}\end{quote}
Throughout this section $\mon$%
\index{S@$\mon$, fixed monomial in $\posv$ in~\S\ref{s.opi},~\S\ref{ss.l.ortho}|ignore}
denotes a monomial in~$\posv$ and~$j$ an integer (not necessarily odd).

\begin{lemma}\mylabel{l.4.9.kr}
If $\block_1,\ldots, \block_l$ are the blocks in order
 from left to right of  some $\mon_k$,   and
$w(\block_1)=(R_1,C_1)$,
$w(\block_2)=(R_2,C_2)$,
$\ldots$,
$w(\block_l)=(R_l,C_l)$,
then \[ C_1<R_1<C_2<R_2<\ldots<R_{l-1}<C_l<R_l	\]
\end{lemma}
\begin{myproof}   This is merely a recall Lemma~4.9 of~\cite{kr}.
In any case it follows easily from the definitions.
\end{myproof}
\begin{lemma}\mylabel{l.4.10.kr}
No two elements of $\mon_k\ext\cup\mon_k'$ are comparable.    More precisely,
it is not possible to have elements $\alpha>\beta$ both belonging to
$\mon_k\ext\cup\mon_k'$.
\end{lemma}
\begin{myproof} It follows from Lemma~\ref{l.4.9.kr} that 
$\monk\cup\monkprime$ contains no comparable elements.
If $k$ is even,  then $\monkext=\monk$ (Corollary~\ref{c.p.monjjp1}~(2));
if $k$ is odd, we may assume $\monkext=\monk$ (as sets) by increasing
the multiplicity of $\sigma_k$ in $\monkpr$.\end{myproof}
\begin{lemma}\mylabel{l.phipi.1}
For integers $i\leq k$,   there cannot exist $\gamma\in\mon'_i\up$ and
$\beta\in\mon_k^\pr$ such that $\beta>\gamma$.
For integers $i<k$,   there cannot exist $\gamma\in\mon'_i\up$ and
$\beta\in\mon_k^\pr$ such that $\beta$ dominates $\gamma$.
\end{lemma}
\begin{myproof}
Let $\gamma\in\mon'_i\up$ and $\beta\in\mon_k^\pr$.
If $i=k$ and $\beta>\gamma$,  then
we get a contradiction immediately to Lemma~\ref{l.4.10.kr}.
Now suppose that $i<k$ and that $\beta$ dominates $\gamma$.
Apply Corollary~\ref{c.existence} (the notation of the corollary
being suggestive of how exactly to apply it).    Let $\alpha$ be
as in its conclusion.     The chain $\alpha>\gamma$ contradicts
Lemma~\ref{l.4.10.kr} in case $i$ is odd
and either Lemma~\ref{l.4.10.kr} or 
Proposition~\ref{p.od.2} in case $i$ is even.\end{myproof}
\begin{lemma}\mylabel{l.4.11.kr}
For $(r,c)$ in $\mon'$,  there exists a unique block~$\blockb$ of~$\mon$
with $(r,c)$ in~$\blockb'$. 
\end{lemma}
\begin{myproof}  The existence is clear from the definition of $\mon'$.
For the uniqueness,  suppose that~$\blockb$ and~$\blockc$ are two distinct
blocks of~$\mon$ with $(r,c)$ in both $\blockb'$ and~$\blockc'$. 
We will show that this leads to a contradiction.

Let $i$ and $k$ be such that $\blockb\subseteq\mon_i$ and
$\blockc\subseteq\mon_k$.  
From Lemma~4.11 of~\cite{kr} (of which the present lemma
is the orthogonal analogue) it follows that $i\neq k$,  so we can assume without
loss of generality that~$i<k$.
By applying the involution~$\hash$ if necessary,
we may assume 
that $(r,c)\in\mon'_i\up$.    
Now there exists an element~$(r,a)$
in~$\blockc$ with $a\leq c$ (this follows from the definition of $\blockc'$).
Clearly $(r,a)\in\mon_k^\pr$.
Taking $\beta=(r,a)$ and $\gamma=(r,c)$,  we get a contradiction to
Lemma~\ref{l.phipi.1}.
%
\end{myproof}
\begin{lemma}\mylabel{l.4.12.kr}
Let $i<j$ be positive integers.
\begin{enumerate}
\item 
Given a block~$\blockb$ of~$\monj$,  there exists
a unique block~$\blockc$ of~$\mon_i$ such that $w(\blockc)>w(\blockb)$.
\item Given an element $\beta$ in
$\mon_j\ext\cup\mon'_j$,  there exists $\alpha$ in $\mon_i$
such that $\alpha>\beta$.
\end{enumerate}
\end{lemma}
\begin{myproof}
(1):
The assertion follows by applying Lemma~\ref{l.kr.!} with $\blockb=\blockb$
and $\monu=\mon_i$.     We need to make sure however that the lemma
can be applied.  More precisely,  we need to check that for every $\beta=(r,c)$
in $\blockb$ there exist $\gamma^1(\beta)=(R^1,C^1)$
and~$\gamma^2(\beta)=(R^2,C^2)$
in $\mon_i$ such that $C^1<c$, $C^2<R^1$, and $r<R^2$.
We may assume
$\beta=\beta\textup{(up)}$, for,  if $\beta=\beta\textup{(down)}$,
then $\beta\textup{(up)}$ also belongs to $\mon_j$ because $\mon_j$
is symmetric,  and we can set
$\gamma^1({\beta})=\gamma^2(\beta\textup{(up)})\textup{(down)}$, and
$\gamma^2({\beta})=\gamma^1(\beta\textup{(up)})\textup{(down)}$---note
that these two belong to~$\mon_i$ since $\mon_i$ is symmetric.

We consider three cases: 
\begin{enumerate}
\item\mylabel{i.1.l.4.12.kr}
$\beta$ belongs to $\mon$.
\item\mylabel{i.2.l.4.12.kr}
$\beta=p_h(\sigma_{j-1})$ (in particular, $j$ is even and $\mon$ is truly
				orthogonal at $j-1$).
\item\mylabel{i.3.l.4.12.kr}
$\beta=p_v(\sigma_j)$ (in particular, $j$ is odd and
				$\mon$ is truly orthogonal at $j$).
\end{enumerate}
Define $\beta'$ to be $\beta$ in case~\ref{i.1.l.4.12.kr},
$\sigma_{j-1}$ in case~\ref{i.2.l.4.12.kr}, and $\sigma_{j}$ in
case~\ref{i.3.l.4.12.kr}.    
Let $C$ be a \vchain\ in $\mon$ with tail $\beta'$ 
and having the good property as in Proposition~\ref{p.goodCexists}.  

First suppose that there exists in $C$ an element of \odepth~$i$
and denote it by~$\gamma$.    
If $p_h(\gamma)\not\in\andposv$ (this can happen only in 
case~\ref{i.1.l.4.12.kr}),  then set 
$\gamma^1(\beta)=\gamma^2(\beta)=\gamma$.
Now suppose $p_h(\gamma)\in\andposv$.  Then
$\gamma\in\mon_i$ except when $\gamma=\sigma_i$ with $i$ odd
and $\sigma_i$ has multiplicity~$1$ in $\mon$.  If $\gamma\in\mon_i$,
take
$\gamma^1(\beta)=\gamma$ and
$\gamma^2(\beta)=\gamma^\hash=\gamma\textup{(down})$;
if $\gamma\not\in\mon_i$,  then take
$\gamma^1(\beta)=\gamma^2(\beta)=p_v(\gamma)$.

Now suppose that $C$ has no element of \odepth~$i$.   Then,
by Lemma~\ref{l.odepths}~(\ref{l.i.hits}), $i$ is
even and there exists in~$C$ an element of \odepth~$i-1$.
This element of~$C$ is of type~H by Lemma~\ref{l.odepths}~(\ref{l.i.odepths}),  
so $\mon$ is truly orthogonal at $i-1$.   Set
$\gamma^1(\beta)=\gamma^2(\beta)=p_h(\sigma_{i-1})$.

(2):  This proof parallels the proof of~(1) above.    As in the above
proof,  we may assume that $\beta=\beta\up$.  Suppose $\beta=(r,c)$ belongs
to $\mon_j'$.  Then there exists $(r,a)\in\mon_j$ with $a\leq c$.
Since $\mon_j'$ does not meet the diagonal,  it is clear that
$(r,a)\in\posv$, and thus it is enough to prove the assertion for
$\beta\in\mon_j\ext$.

So now take $\beta\in\mon_j\ext$.
Let $\beta'$ and $C$ be in the proof of~(1).  First suppose that
there exists in~$C$ an element of \odepth~$i$.  Denote it by $\gamma$.    
If $\gamma\in\mon_i$, then take $\alpha=\gamma$.   If $\gamma\not\in\mon_i$,
then $p_v(\gamma)\in\mon_i$, and we take $\alpha=p_v(\gamma)$.   In case
there is no element in~$C$ of \odepth~$i$,  we take $\alpha=p_h(\sigma_{i-1})$
(see the above proof). 
\end{myproof}
\begin{corollary}\mylabel{c.4.13.kr}
If $\block$ and $\block_1$ are blocks of $\mon$ with $w(\block)=(r,c)$
and $w(\block_1)=(r_1,c_1)$,  then exactly one of the following holds:
\begin{eqnarray*}
        c<r<c_1<r_1,&\quad\quad\quad\quad& c_1<r_1<c<r,\\
  c<c_1<r_1<r,&\quad\quad\mbox{\begin{rm}or\end{rm}}
  \quad\quad& c_1<c<r<r_1.
\end{eqnarray*}
\end{corollary}
\begin{myproof}
This is a formal consequence of Lemmas~\ref{l.4.9.kr} and~\ref{l.4.12.kr},
just as Corollary~4.13 of~\cite{kr} is of Lemmas~4.9 and~4.12 of that paper.
\end{myproof}
\newcommand\bcor{\begin{corollary}}
\newcommand\ecor{\end{corollary}}
\bcor\mylabel{c.4.16pr.kr}
If $w(\blockb)>w(\blockc)$ for blocks $\blockb\subseteq\mon_i$ and
$\blockc\subseteq\mon_j$ of $\mon$,  then $i<j$.
\ecor
\begin{myproof}
This is a formal consequence of Lemmas~\ref{l.4.9.kr} and~\ref{l.4.12.kr}.   
It follows from the first lemma that $i\neq j$.    Suppose $i>j$.   Then
there exists by the second lemma a block~$\blockc'\subseteq\mon_j$
such that $w(\blockc')>w(\blockb)$.   But then $w(\blockc')>w(\blockc)$,
a contradiction of the first lemma.
\end{myproof}
\bcor\mylabel{c.4.14.kr}
Let $(s,t)>(s_1,t_1)$ be elements of $\mon'$, and $\block$,
$\block_1$ be blocks of $\mon$ such that $(s,t)\in\block'$,
and $(s_1,t_1)\in\block_1'$.  Then $w(\block)>w(\block_1)$.
\ecor
\begin{myproof}
Let $w(\block)=(r,c)$ and $w(\block_1)=(r_1,c_1)$.    
By Corollary~\ref{c.4.13.kr},   we have four possibilities.
Since $(r,c)$ dominates $(s,t)$ and $(r_1,c_1)$ dominates $(s_1,t_1)$,
the possibilities $c<r<c_1<r_1$ and $c_1<r_1<c<r$ are eliminated.
It is thus enough to eliminate the possibility $c_1<c<r<r_1$.
Suppose that this is the case.   Then, by Corollary~\ref{c.4.16pr.kr},
$j_1<j$, where $j_1$ and $j$
are such that $\blockb\subseteq\mon_j$ and $\blockb_1\subseteq\mon_{j_1}$.
Now, by Lemma~\ref{l.4.12.kr}~(2), there exists $\alpha$ in $\mon_{j_1}$
such that $\alpha>(s,t)>(s_1,t_1)$.     But then this contradicts
Lemma~\ref{l.4.10.kr}.
\end{myproof}
\bcor\mylabel{c.4.16.kr}
For a $\blockb\subseteq\mon_i$ of $\mon$,  
the depth of $w(\blockb)$ in $\mon_w$ is exactly $i$.
\ecor
\begin{myproof}
That the depth is at least $i$ follows from Lemma~\ref{l.4.12.kr}.   
That the depth cannot exceed~$i$ follows from Corollary~\ref{c.4.16pr.kr}.
\end{myproof}
\begin{corollary}\mylabel{c.monprime}
Let $\alpha\in\mon_k'\up$, $\beta\in\mon_m'\up$, and $\alpha>\beta$.
Then $k<m$.
\end{corollary}
\begin{myproof} Corollary~\ref{c.4.14.kr} and Corollary~\ref{c.4.16.kr}.
\end{myproof}
\mysubsection{More lemmas}\mylabel{ss.pf.lemmas}
This subsection is a collection of lemmas to be invoked
in the later subsections.    More specifically, Lemma~\ref{l.sprime}
and Corollary~\ref{c.l.sprime} are invoked in the proof of
Proposition~\ref{p.4.1.kr} in \S\ref{ss.pf.p.4.1.kr},
Lemma~\ref{l.phipi} in the proof of the first half of 
Proposition~\ref{p.4.2.kr} in \S\ref{ss.pf1.p.4.2.kr}, and
Lemma~\ref{l.piphi} in the proof of the second half of 
Proposition~\ref{p.4.2.kr} in \S\ref{ss.pf2.p.4.2.kr}. Throughout this subsection, $\mon$ denotes a monomial in $\posv$.
\begin{lemma}\mylabel{l.sprime}
Let $C$ be a $v$-chain in $\mon'$,   $\alpha$ an element of~$C$,
and $\alpha'\in\mon_{C,\alpha}$.    Then
$\depthin{\mon_C}{\alpha'}\leq k\even$,  for
$k$ the integer such that $\alpha\in\monkprimeup$.
\end{lemma}
\begin{myproof}
Proceed by induction on~$k$.   If $k=1$,  the assertion follows
from Corollary~\ref{c.monprime}, so assume $k>1$.
%
Choose a $v$-chain $C'$ in $\mon_C$ with tail $\alpha'$ and
$\depthin{C'}{\alpha'}=\depthinmonC{\alpha'}$.    The length of
a \vchain\ in $\mon_{C,\alpha}$ is clearly at most~$2$.
So, if~$\gamma'$ is the element two steps  
before $\alpha'$ 
in $C'$ (if $\gamma'$ does not exist then there is clearly
nothing to prove),   then
$\gamma'\in\mon_{C,\gamma}$ with $\gamma>\alpha$
(see Proposition~\ref{p.type.>}~(2)).
We claim that
$\depthinmonC{\gamma'}\leq k\odd-1$.   It is enough to prove the claim,
for then $\depthinmonC{\alpha'}=\depthin{C'}{\alpha'}=\depthin{C'}{\gamma'}+2
	\leq k\odd-1+2=k\even$.

The claim follows by induction from Corollary~\ref{c.monprime}
if $k$ is odd 
or more generally if $\gamma\in\mon_l'\up$ with
$l\leq k\odd-1$.   So assume that $k$ is even and $\gamma\in\mon_{k-1}'\up$.
By~\ref{p.od.2},  it is not possible that $\gamma$ is of type~H and 
$p_h(\gamma)>\alpha$.   So the only possibility is that $\alpha'=p_h(\alpha)$
and $\gamma>\alpha$ is connected.  In particular, $\gamma$ is of type~V 
and $\alpha$ of type~H
in~$C$ and $\gamma'=p_v(\gamma)$.   

Now let $\mu$ be the first element
in the connected component of~$\alpha$ in~$C$.    The cardinality of 
the part $\mu>\ldots>\gamma$ of~$C$ is even (by Proposition~\ref{p.type}~(1),
it follows that the cardinality of $\mu>\ldots>\alpha$ is odd),  say~$e$.
Letting $m$ be such that $\mu\in\mon_m'\up$, we have,  
by Proposition~\ref{c.monprime}, $m\leq k-1-(e-1)=k-e$.
If $m\even<k-e$, then,
since $\depthin{C'}{\gamma'}=\depthin{C'}{p_v(\mu)}+e-1$
(by Proposition~\ref{p.type.>}~(1), since, 
by Proposition~\ref{p.type}~(\ref{p.type.hs}),
$\mu$,~\ldots,~$\gamma$ all have type~V in~$C$)
and $\depthin{C'}{p_v(\mu)}\leq m\even$ by induction, it follows that
$\depthin{C'}{\gamma'}<k-e+e-1=k-1$, and we are done.

So suppose that $m\even=k-e$.   
Let $\nu$ be the element just before $\mu$ in~$C$ (if such an element
does not exist,  then $\depthin{C'}{\gamma'}=e\leq k-2$---observe
that $m\even\geq2$---and we are done).
Then $\nu>\mu$ is not connected (by choice of $\mu$).  
So $p_h(\nu)>\mu$.
By Proposition~\ref{p.od.2},  this means that $j\leq m\even-2$
where $j$ is the odd integer defined by
$\nu\in\monjprime\up\cup\monjponeprime\up$.
So, again by induction, $\depthin{C'}{\gamma'}
= \depthin{C'}{p_h(\nu)}+e\leq m\even-2+e=k-2$, and the claim is proved.
\end{myproof}
\begin{corollary}\mylabel{c.l.sprime}
The \odepth\ of an element $\alpha$ in $\mon'$ is at most $k$
where $k$ is such that $\alpha\in\mon_k'\up$.
\end{corollary}
\begin{myproof}
Let $C'$ be a \vchain\ in $\mon_C$ with tail $q_{C,\alpha}$.
If $k$ is even,  then, by the lemma, $\depthin{C'}{q_{C,\alpha}}\leq k$.
So suppose that $k$ is odd.    Let $\gamma'$ be the immediate predecessor
of $q_{C,\alpha}$ in $C'$.   By Proposition~\ref{p.type.>}~(2),  $\gamma>\alpha$,
and so $\gamma\in\mon_l'\up$ with $l\leq k-1$ (see the observation in the
first paragraph of the proof of the lemma).   So $\depthin{C'}{\gamma'}\leq
k-1$ (by the lemma) and $\depthin{C'}{\alpha'}=\depthin{C'}{\gamma'}+1\leq k$.
\end{myproof}
\begin{lemma}\mylabel{l.phipi}
Let $\mon$ be a monomial in $\posv$ and $\opi(\mon)=(w,\mon')$.
Let $i<k$ be integers, $\alpha$ an element of $\mon_i'\up$, and
$\delta$ an element of $(\monw)_k\up$ that dominates~$\alpha$.
\begin{enumerate}
\item  If $k$ is even,  then there exists $\beta\in\monkprime\up$
with $\alpha>\beta$.
\item If $k$ is odd and  $w_{k,k+1}$ \ortho-dominates the singleton
\vchain~$\alpha$, 
then either there exists $\beta\in\monkprime\up$ with $\alpha>\beta$
or there exists $\gamma\in\mon_{k+1}'\up$ with $p_h(\alpha)>\gamma$.
\end{enumerate}
\end{lemma}
\begin{myproof} 
Write $\alpha=(r,c)$ and $\delta=(A,B)$.
By Corollary~\ref{c.4.16.kr},  there exists a block~$\blockb$ of $\mon_k$
such that $\delta=w(\blockb)$.   Let $(D,B)$ be the first element of~$\blockb$
(arranged in increasing order of row and column indices).
We have the following possibilities:
\begin{enumerate}
\item[(i)]
$D\leq A$ and $(D,B)\in\mon_k^\pr$.
\item[(ii)]
$k$ is odd, $\mon$ is truly orthogonal at $k$,
$(D,B)=(A,B)=p_v(\sigma_k)$, and $\blockb$ consists
of the single diagonal element $(D,B)=(B^*,B)$.
\item[(iii)]
$k$ is even, $\mon$ is truly orthogonal at $k-1$,
$(D,B)=(A,B)=p_h(\sigma_{k-1})$, and $\blockb$ consists
of the single diagonal element $(D,B)=(B^*,B)$.
\end{enumerate}
We claim the following: in case~(i), $D<r$ (in particular, $D<A$); 
in case~(ii), the row index of $\sigma_k$ is less than $r$; and
case~(iii) is not possible.
The first two assertions and also the third in the case $i<k-1$ 
follow readily from Lemma~\ref{l.phipi.1};  in case~(iii) holds
and $i=k-1$,  then $\sigma_{k-1}>\alpha$, a contradiction to
Lemma~\ref{l.4.10.kr}.

First suppose that possibility~(ii) holds.
Write $\sigma_k=(s,B)$.  Since $s<r$ and $p_h(\sigma_k)\in\andposv$,
it is clear that $p_h(\alpha)=(r,r^*)$ also belongs to $\andposv$.
From the hypothesis that $w_{k,k+1}$ \ortho-dominates~$\{\alpha\}$,
it follows that there is an element of 
$(\mon_w)_{k+1}$ that dominates $p_h(\alpha)=(r,r^*)$.  Such an
element must be diagonal (because of the distinguishedness of $\mon_w$),
and so must be the $w(\blockc)$ for the unique diagonal block 
$\blockc$ of $\mon_{k+1}$.    In particular,  this means that there
are elements other than $(s,s^*)$ in $\mon_{k+1}$,  and so $\mon_{k+1}'$
is non-empty.    
In the arrangement of elements of $\mon_{k+1}'\up$ in increasing
order of row and column numbers,  let $\gamma=(e,s^*)$ be the last element.
Then $e<s<r$ and $r^*<s^*$, so $p_h(\alpha)>\gamma$, and we are done.

Now suppose that possibility~(i) holds.  
Let $(p,q)$ be the element of $\mon_k$ such that~$p$ is the
largest row index that is less than $r$,  and, among those elements with
row index $p$, the maximum possible column index is~$q$.
The arrangement of elements of $\mon_k$ (in increasing order or row
and column indices) looks like this:
\begin{center}
\begin{tabular}{c}
\ldots, $(p,q)$, $(s,t)$, \ldots
				\end{tabular}
				\end{center}
Since $p<r\leq A$ and $w(\block)=(A,B)$,  we can be sure that $(p,q)$
is not the last element of $\blockb$.    

We first consider the case $c<t$.   
Then $\alpha=(r,c)>\beta:=(p,t)\in\monkprime$.    If $\beta\in\monkprime\up$,
then we are done.  It is possible that
$(p,q)$ lies on or below the diagonal so that $\beta$ lies below the diagonal,
in which case, $\alpha>\beta\up$ and $\beta\up\in\monkprime\up$, and
again we are done.

Now suppose that $t\leq c$.    We claim that:
\begin{itemize}
\item $(s,t)$ belongs to the diagonal;
\item $k$ is odd and $\mon$ is truly orthogonal at $k$; and
\item $\sigma_k=(u,t)$  with $u<r$.
\end{itemize}
Suppose that $(s,t)$ does not belong to the diagonal.
Since $r\leq s$ (by choice of $(p,q)$),  it follows that
 $(s,t)$ dominates $(r,c)$.
This leads to a contradiction to Lemma~\ref{l.phipi.1}, for either
$(s,t)$ or its reflection $(t^*,s^*)$ (whichever is above the diagonal)
belongs to~$\monkpr$ and dominates $\alpha=(r,c)$ in~$\mon_i^\prime\up$.
This shows that $(s,t)$ belongs to the diagonal.
If $k$ is even,  then $(s,t)=p_h(\sigma_{k-1})$,  which means
$\sigma_{k-1}>\alpha$,  again contradicting Lemma~\ref{l.phipi.1}, so
$k$ must be odd.  It also follows that $\mon$ is truly orthogonal
at $k$ and that $(s,t)=p_v(\sigma_k)$.   Writing $\sigma_k=(u,t)$,
if $r\leq u$,  then~$\sigma_k$ would dominate $\alpha$,  again
contradicting Lemma~\ref{l.phipi.1}.    So $u<r$, and the claim
is proved.    

To finish the proof of the lemma, now proceed as in the proof
when possibility~(ii) holds.\end{myproof}
\begin{lemma}\mylabel{l.piphi}\mylabel{l.4.21.kr}
Let $\mont$ be a monomial in $\posv$ and $w$ an element of $I(d)$ that
\ortho-dominates $\mont$.     
Let $\beta'>\beta$ be elements $\mon_w\up$.
Let $d-1$ and $d$ be their respective depths in $\mon_w$.
Let $\alpha$ be an element of $\opbstar$ or more generally an element 
of $\posv$ such that
\begin{enumerate}
\item[(a)]
it is dominated by $\beta$,
\item[(b)] it is not comparable to any element of $\pbeta$, and
\item[(c)] in case $d$ is odd, then $\{\alpha\}\cup\mont_{w,d,d+1}$ has 
\odepth\ at most~$2$.
\end{enumerate}
Then
\begin{enumerate}
\item 
there exists $\alpha'\in\pbpstar\up$
with $\alpha'>\alpha$;
\item
for $\alpha'$ as in~(1), if $\alpha'$ is diagonal,  
then $p_h(\delta_{d-2})>\alpha$ if $d$ is odd and
$\delta_{d-1}>\alpha$ if $d$ is even.
\end{enumerate}
\end{lemma}
\begin{myproof}
Assertion~(2) is rather easy to prove.   If $d$ is odd,  then, in fact,
$p_h(\delta_{d-2})=\alpha'$;   if $d$ is even,  then $\delta_{d-1}$ has the
same column index as $\alpha'$ and, by Proposition~\ref{p.t.od<2}~(4),
has row index more than that of $\alpha$,  so $\delta_{d-1}>\alpha$.

Let us prove~(1).   
Write $\alpha=(r,c)$, $\beta=(R,C)$, and $\beta'=(R',C')$.
There exists, by the definition of $\pbpstar$,  an element in $\pbpstar$
with column index $C'$.     We have $C'<c$ (for $C'<C\leq c$).  Let $(r',c')$
be the element of $\pbpstar$ such that $c'$ is maximum possible
subject to $c'<c$ and among those elements with column index $c'$
the maximum possible row index is $r'$.    If $r<r'$,  then we
are done (if $(r',c')$ is below the diagonal,  its mirror image
would have the desired properties).     It suffices therefore
to suppose that $r'\leq r$ and arrive at a contradiction.

In the arrangement of elements of $\pbpstar$
in non-decreasing order of row and column indices,  there is a portion
that looks like this:
\[   \ldots,~(r',c'),~(a,b),~\ldots  \]
Since there is in $\pbstar$ an element with row index $R'$ (and clearly
$r'\leq r<R<R'$),   it follows that $(a,b)$ exists (that is, $(r',c')$
is not the last element in the above arrangement).   It follows from
the construction of $\pbpstar$ from $\pbetap$ that $(r',b)$ is an element
in $\pbetap$.    By the choice of $(r',c')$, we have $c\leq b$.
Thus $(r,c)$ dominates $(r',b)$.

The proof now splits into two cases accordingly as $d$ is even or odd.
First suppose that $d$ is even.    Then, since $\beta$ dominates
$(r',b)$ and yet $(r',b)$ does not belong to~$\pbeta$,   there exists
a \vchain\ in $\mont_{w,d-1,d}$ of length~$2$ and head $(r',b)$.
The tail of this \vchain\ then belongs to $\pbeta$ and is dominated
by $(r,c)$,  a contradiction to our assumption that $\alpha$ is 
not comparable to any element of $\pbeta$. 

Now suppose that $d$ is odd.
Choose a \vchain~$C$ in~$\mont$ with head $(r',b)$ that is 
not \ortho-dominated by~$w^d$.   Let $D$ be the part of $C$ consisting
of elements of \ortho-depth (in~$C$) at most~$2$.  
We claim that $D$ is \ortho-dominated by $w_{d,d+1}$.   
In fact, we claim the following:
Any \vchain~$F$ with head $(r',b)$ and \odepth\ at
most~$2$ is \ortho-dominated by $w_{d,d+1}$.   

To prove the claim, we first prove the following subclaim:
\begin{quote}
(\dag)  If the horizontal projection of $(r',b)$ belongs to $\andposv$,
then $\beta$ is on the diagonal and dominates the vertical projection
of $(r',b)$,   and the diagonal element $\beta_1$ of $(\mon_w)_{d+1}$
dominates the horizontal projection of $(r',b)$.
\end{quote}
Let $p_h(r',b)\in\andposv$.   Then $\beta$
 belongs to the diagonal because $\mon_w$ is distinguished and symmetric.    
Once $\beta$
is on the diagonal,  it is clear that it dominates $p_v(r',b)$ 
(from our assumptions, $\beta$ dominates $(r,c)$ and $(r,c)$
 dominates $(r',b)$).     It follows from Proposition~\ref{p.t.od<2}~(3) that
the row index of $\beta_1$ exceeds the row index~$r$ of $(r,c)$, so $\beta_1$
dominates $p_h(r',b)$.   This finishes the proof of the subclaim (\dag).

To begin the proof of the claim,  observe that $F$ has length at most~$2$.
Suppose first that $F$ consists only of the single element~$(r',b)$.
The type of $(r',b)$ in $F$ is either~H or~S. 
If it is~S,  then since $\beta$ dominates $(r',b)$,
the claim follows immediately.
If it is~H,  then the claim follows immediately from the subclaim~(\dag).

Continuing with the proof of the claim,  let now $F$ consist of 
two elements: $(r',b)>\mu$.   
Let $\gamma$ be the element of $\mon_w$ such that $\mu\in\piece_\gamma$,
and let $e$ be the depth of $\gamma$ in~$\mon_w$.  
From Lemma~\ref{l.jjp}~(2b) it follows that $e\geq d$.
If $e=d$,  then $\gamma=\beta$ (by the distinguishedness of $\mon_w$),
and the comparability of $(r,c)$ and $\mu$ contradicts 
our hypothesis~(b).
So $e\geq d+1$,  and there exists
$\delta$ of depth $d+1$ in $\mon_w$ that dominates $\mu$.   
We have $\beta>\delta$ (again by the distinguishedness of $\mon_w$).

The possibilities for the types of $(r',b)$ and $\mu$ in $F$ are:
S and~S, V and~V, H and~S (in the last case $p_h(r',b)\not>\mu$ by
Lemma~\ref{l.odepths}~(\ref{l.i.odepths})).  Noting the existence in
$(\mon_w)_{d,d+1}$ of the \vchain\ $\beta>\delta$ in the first case
and also of $\beta>\beta_1$ (where $\beta_1$ is as in the subclaim)
in the last case,  the proof of the claim
in these cases is over.   So suppose that the second possibility
holds.
The distinguishedness of $\mon_w$ implies that $\delta=\beta_1$.
Since $\delta$ is diagonal,
it dominates the vertical projection of $\mu$.
Noting the existence of the \vchain\ in $\beta>\delta$ in 
$(\mon_w)_{d,d+1}$,  the proof of the claim in this case too is over.

We continue with the proof of the lemma.   It follows from the claim that 
$D$ is \ortho-dominated by $w_{d,d+1}$.
From Corollary~\ref{c.odomination} it follows that 
the complement~$E$ of~$D$ in~$C$ is not \ortho-dominated by
$w_{d+2,d+3}$ (in particular, that $E$ is non-empty) 
and that every \vchain\ in~$\mont$ with head~$\epsilon$
(where $\epsilon$ denotes the head of~$E$) is \ortho-dominated 
by~$w^d$ (given such a \vchain,  the concatenation of~$D$ with
it is \ortho-dominated by $w^{d-2}$, and $\epsilon$ continues to have
\odepth~$3$ in the concatenated \vchain).   Thus $\epsilon$ belongs to $\mont_{w,d,d+1}$.   From (1) and (2b) of
Lemma~\ref{l.jjp} it follows that the element $\mu$ of~$C$
in between $(r',b)$ and $\epsilon$ (if it exists at all) also
belongs to $\mont_{w,d,d+1}$.      Now consider the \vchain\
obtained as follows:   take the part of $C$ up to (and including)~$\epsilon$
and replace its head~$(r',b)$ by $(r,c)$.     This chain has \odepth~$3$
and lives in $\{\alpha\}\cup\mont_{w,d,d+1}$,  a contradiction
to hypothesis~(c).
\end{myproof}
\newcommand\opbppstar{\ortho\piece_{\beta''}^*}
\begin{corollary}\mylabel{c1.l.piphi}
Let $\mont$ be a monomial in $\posv$ and $w$ an element of $\id$ that
\ortho-dominates~$\mont$.  Let $\beta'>\beta$ be elements of
$\monw\up$, $\alpha$ an element of $\opbstar$, 
and $d':=\depthin{\monw}{\beta'}$.
\begin{enumerate}
\item If $d'$ is odd,  there exists $\alpha'\in\opbpstar$ such that
  $\alpha'>\alpha$.
\item If there does not exist $\alpha'\in\opbpstar$ such that 
$\alpha'>\alpha$ then ($d'$ is even by (1) above and) there exists
$\alpha''\in\opbppstar$ such that $p_h(\alpha'')>\alpha$, where
$\beta''$ is the unique element of $(\monw)_{d'-1}$ such that $\beta''>\beta'$.
\end{enumerate}
\end{corollary}
\begin{myproof}
  Immediate from the lemma.
\end{myproof}
\begin{corollary}\mylabel{c2.l.piphi}
  Let $\mont$ be a monomial in $\posv$ and $w$ an element of $\id$ that
\ortho-dominates~$\mont$.   Let $\beta$, $\beta'$ be elements of
$\monw\up$, and $\alpha$, $\alpha'$ elements of $\opbstar$ and $\opbpstar$
respectively.    
\begin{enumerate}
\item If $\alpha'>\alpha$ then $\beta'>\beta$ (in particular, 
$\depthin{\monw}{\beta'}<\depthin{\monw}{\beta}$).
\item If $p_h(\alpha')>\alpha$ and $\depthin{\monw}{\beta}$ is even,
$\depthin{\monw}{\beta'}\leq\depthin{\monw}{\beta}-2$.
\end{enumerate}
    \end{corollary}
    \begin{myproof}
(1)~Writing $\beta=(r,c)$ and $\beta'=(r',c')$,   there are,
since both $\beta$ and $\beta'$ dominate~$\alpha$ and~$\monw$ is 
distinguished, the following four
possibilities:
\[ c<r<c_1<r_1,\quad c_1<r_1<c<r, \quad c<c_1<r_1<r, \quad
c_1<c<r<r_1
\]
Since $\alpha'>\alpha$,  and $\alpha$, $\alpha'$
are dominated respectively by $\beta$, $\beta'$ (this is because 
$\alpha$, $\alpha'$ belong to $\opbstar$, $\opbpstar$ respectively),
the possibilities $c<r<c_1<r_1$ and $c_1<r_1<c<r$ are eliminated (by the
distinguishedness of $\monw$).   It is thus
enough to eliminate the possibility $\beta>\beta'$.     Suppose, by
way of contradiction, that $\beta>\beta'$.   By Corollary~\ref{c1.l.piphi},
either there exists $\gamma\in\opbstar$ such that $\gamma>\alpha'$,
in which case the \vchain~$\gamma>\alpha$ in~$\opbstar$
contradicts Proposition~\ref{p.t.od<2}~(3) or~(4),
or $d:=\depthin{\monw}{\beta}$ is even and there exists (with $\beta''$
being the unique element in $\monw$ such that $\beta''>\beta$ and
$\depthin{\monw}{\beta''}=d-1$) an element $\alpha''\in\opbppstar$
with $p_h(\alpha'')>\alpha'$,  in which case the \vchain~$\alpha''>\alpha$
in $\mont_{w,d-1,d}^\star$ has \odepth~$3$ and so contradicts 
Proposition~\ref{p.t.od<2}~(2).

(2)~Set $d:=\depthin{\monw}{\beta}$.  If $\depthin{\monw}{\beta'}$
were $d-1$,  then the \vchain~$\alpha'>\alpha$ in $\mont_{w,d-1,d}^\star$
would be of \odepth~$3$ and so would contradict Proposition~\ref{p.t.od<2}~(2).
    \end{myproof}
\mysection{The Proof}\mylabel{s.proof}
The aim of this section is to prove Propositions~\ref{p.4.1.kr}
and~\ref{p.4.2.kr}.    The proof of first proposition appears in
\S\ref{ss.pf.p.4.1.kr} and that of the second in
\S\S\ref{ss.pf1.p.4.2.kr},~\ref{ss.pf2.p.4.2.kr}.
In \S\ref{ss.pf.lemmas} some lemmas are established that are used
in the proofs.     Needless to say that the lemmas maybe unintelligible
until one tries to read the proofs in the later subsections.
\mysubsection{Proof of Proposition~\ref{p.4.1.kr}}%
\mylabel{ss.pf.p.4.1.kr}
(\ref{i.1.p.4.1.kr}) By definition,  $w$ is the element of $\id$
associated to the distinguished monomial $\cup_k\mon_{w(k)}$.   By
the very definition of this association,  we have $w\geq v$.

(\ref{i.2.p.4.1.kr}) This follows from the corresponding property of
the map~$\pi$ of~\cite{kr}.   More precisely,  that property justifies
the third equality below.   The other equalities are clear from the
definitions.
\begin{eqnarray*} \vdeg{w}+\degree{\mon'} &  = &
	\frac{1}{2} \degree{\mon_w} + \frac{1}{2}\sum_k\degree{\mon_k'}\\
& = & \frac{1}{2}\sum_k\left(\degree{\mon_{w(k)}}+\degree{\mon_k'}\right) \\
& = & \frac{1}{2}\sum_k\degree{\mon_k} \\
& = & \frac{1}{2}\sum_{\textup{$j$ odd}}\degree{\monjjpone} \\
& = & \sum_{\textup{$j$ odd}}\left(\degree{\monjonlypr}+\degree{\monjponepr}
	\right)\\
& = & \degree{\mon} \end{eqnarray*}

(\ref{i.3.p.4.1.kr}) We have:\\[2mm]
\begin{tabular}{rcl}
$w$ \ortho-dominates $\mon'$
& $\Leftrightarrow$ &
	$w\geq w_C$ $\forall$ \vchain\ $C$ in~$\mon'$\\
& $\Leftrightarrow$ &
	$w$ dominates $\mon_C$ $\forall$ \vchain\ $C$ in~$\mon'$\\
& $\Leftrightarrow$ &
	$\forall$ \vchain\ $C$ in~$\mon'$, $\forall$
	$\alpha'=(r,c)\in\mon_C$,\\
	&&
	$\exists$ $\beta=(R,C)\in\mon_w$ with
	$C\leq c$, $r\leq R$,\\
	&&
	and $\depthin{\mon_w}{\beta}\geq\depthin{\mon_C}{\alpha'}$.\\
\end{tabular}\\[2mm]
The first equivalence above follows from the definition of
\ortho-domination,  the second from~\cite[Lemma~4.5]{kr},  the third
from Lemma~\ref{l.domination}.   

Now let $C$ be a \vchain\ in~$\mon'$ and $\alpha'=(r,c)$ in $\mon_C$.
We will show that there exists $\beta$ in $\monw$ that dominates $\alpha$
and satisfies $\depthin{\monw}\beta\geq\depthin{\mon_C}\alpha'$.
%
Let $\alpha$ be the element in $C$ such that $\alpha'\in\mon_{C,\alpha}$,
let $k$ be such that $\alpha\in\mon_k'\up$, and
let $\blockb$ be the block of $\mon_k$ such that $\alpha\in\blockb'$.
Writing $\alpha=(r_1,c_1)$ and $w(\blockb)=(R_1,C_1)$,  we have
$C_1\leq c_1$ and $r_1\leq R_1$ straight from 
the definition of $w(\blockb)$.
By Corollary~\ref{c.4.16.kr}, $\depthin{\mon_w}{w(\blockb)}=k$.

First suppose that $w(\blockb)$ dominates $\alpha'$ (meaning
$C_1\leq c$ and $r\leq R_1$).  
If $k\geq \depthin{\mon_C}{\alpha'}$,   we are clearly done;
by Corollary~\ref{c.l.sprime},  this is the case when $\alpha'=q_{C,\alpha}$.
So suppose that $\alpha$ is of type~H, $\alpha'=p_h(\alpha)$, and that 
$k<\depthin{\mon_C}{\alpha'}$.    By Lemma~\ref{l.sprime},  
$\depthin{\mon_C}{\alpha'}\leq k\even$.  It follows that $k$ is odd
and $\depthin{\mon_C}{\alpha'}=k+1$.      By Corollary~\ref{c.p.prepare},
$\mon$ is truly orthogonal at $k$, which means that 
$\mon_{k+1}$ has a diagonal block, say $\blockc$.  Note that
$w(\blockc)$ dominates $p_h(\sigma_k)$ which in turn dominates $p_h(\alpha)$.
Since $\depthin{\mon_w}w({\blockc})=k+1$ by Corollary~\ref{c.4.16.kr},
we are done.

Now suppose that $w(\block)$ does not dominate $\alpha'$.
Then $\blockb$ is non-diagonal and $\alpha'=p_v(\alpha)$.
Since $\blockb$ is non-diagonal, $p_h(\alpha)\not\in\andposv$,
and $\alpha$ cannot be of type~H.   So $\alpha$ is of type~V in $C$.
It follows easily (see Proposition~\ref{p.type}~(\ref{p.type.cr})) 
that $\alpha$ is the critical element in~$C$, and 
and that last element in its connected component in $C$;   by
Lemma~\ref{l.odepths}~(\ref{l.i.vone}),  
$\odepthin{C}{\alpha}=\depthin{\mon_C}{q_{C,\alpha}}=:d$ is even.
By Proposition~\ref{p.type}~(\ref{p.type.odd}),~(\ref{p.type.hs}),
the cardinality of the connected component of $\alpha$ in $C$ is even.

The immediate predecessor $\gamma$ of $\alpha$ in $C$ is connected to $\alpha$
(this follows from what has been said above).
It is of type~V in $C$, $p_h(\gamma)$ belongs to $\andposv$,
and $\depthin{\mon_C}{p_v(\gamma)}=d-1$
(see Lemma~\ref{l.odepths}~(\ref{l.i.odepths})).    Let $\ell$
be such that $\gamma\in\mon'_\ell\up$. 
Let $\blockc$ be the block of $\mon_\ell$ such that $\gamma\in\blockc'$.
Since $p_h(\gamma)\in\andposv$,   $\blockc$ is diagonal.
Note that $w(\blockc)$ dominates $p_v(\gamma)$ and that
$p_v(\gamma)>p_v(\alpha)$.
By Corollary~\ref{c.4.16.kr},
$\depthin{\mon_w}{w(\blockc)}=\ell$.    Thus if $d\leq \ell$ we are done.
On the other hand, $d-1\leq \ell$ by Corollary~\ref{c.l.sprime}.

\newcommand\blockd{\mathfrak D}
So we may assume that $\ell=d-1$.
By Corollary~\ref{c.p.prepare},
$\mon$ is truly orthogonal at $d-1$.    This implies that
$\mon_{d}$ has a diagonal block, say $\blockd$.
Note that $w(\blockd)$ dominates $p_h(\sigma_{d-1})$ which
in turn dominates $p_h(\gamma)$.
Writing $\gamma=(r_2,c_2)$,   since $\gamma>\alpha$ is connected,
it follows that $(r_1,r_2^*)$ belongs to $\posv$.   
Now both $w(\blockb)$ and $w(\blockd)$ dominate $(r_1,r_2^*)$.
Since $\mon_w$ is distinguished and symmetric
and $w(\blockb)$ is not on the diagonal,  it follows that
$w(\blockd)>w(\blockb)$.    This implies, since $w(\blockd)$
is on the diagonal, $w(\blockd)>p_v(\alpha)$.   
Since $\depthin{\mon_w}{w(\blockd)}=d$ by Corollary~\ref{c.4.16.kr},
we are done.

(4)~Let $x$ be an element of $I(d)$ that \ortho-dominates $\mon$.
We will show that $x\geq w$.   By~\cite[Lemma~5.5]{kr},  it is enough
to show that $x$ dominates $\mon_w$.  By Lemma~\ref{l.domination},
it is enough to show the following:  for every block $\blockb$ of
$\mon$,  there exists $\beta$ in $\mon_x$ such that $\beta$ dominates
$w(\blockb)$ and $\depthin{\mon_x}{\beta}\geq\depthin{\mon_w}{w(\blockb)}$.

Let $\blockb$ be a block of $\mon$.   By Corollary~\ref{c.4.16.kr},
$\depthin{\mon_w}{w(\blockb)}=k$ where $\blockb\subseteq\mon_k$.
Let $\mon_x^k$ denote the set of elements of $\mon_x$ of depth at least~$k$.
Our goal is to show that there exists $\beta$ in $\mon^k_x$ that
dominates $w(\blockb)$.   
It follows easily from the distinguishedness of $\mon_x$ and
the fact that $\blockb$ is a block,  that it suffices to show the 
following: given $\alpha\in\blockb$, there exists $\beta$ in $\mon_x^k$
(depending upon $\alpha$) that dominates~$\alpha$.
Moreover,  since $\blockb$ and $\mon_x^k$ are symmetric,   we may assume that
$\alpha=\alpha\up$.

So now let $\alpha=\alpha\up$ belong to $\blockb$.   Then either
\begin{enumerate}
\item
$\alpha$ belongs to $\monkpr$, or
\item
$k$ is odd, $\mon$ is truly orthogonal at $k$, and
$\alpha=p_v(\sigma_k)$, or
\item
$k$ is even, $\mon$ is truly orthogonal at $k-1$, and
$\alpha=p_h(\sigma_{k-1})$.
\end{enumerate}
The proofs in the three cases are similar.   In the first case,
choose a \vchain~$C$ in~$\mon$ with tail $\alpha$ such that
$\odepthinC{\alpha}=k$ (see Corollary~\ref{c.p.odepth}~(1)).
Then $\depthin{\mon_C}{q_{C,\alpha}}=k$ and, clearly, $q_{C,\alpha}$ dominates
$\alpha$.   Since $x$ dominates $\mon_C$,   there exists, by~Lemma~4.5 
of~\cite{kr},  $\beta$ in $\mon^k_x$ that dominates $q_{C,\alpha}$
(and so also $\alpha$).

In the second case,  
choose a \vchain~$C$ in $\mon$ with tail $\sigma_k$
with the property that $\odepthinC{\sigma_k}=k$.
Then $\depthin{\mon_C}{q_{C,\sigma_k}}=k$.  Since $p_h(\sigma_k)$
belongs to $\andposv$,  $\sigma_k$ is of type~V or~H in $C$,
so $q_{C,\sigma_k}=\alpha$.
Since $x$ dominates $\mon_C$,   there exists, by
\cite[Lemma~4.5]{kr},  $\beta$ in $\mon^k_x$ that dominates
$q_{C,\sigma_k}=\alpha$.

In the third case,  
choose a \vchain~$C$ in $\mon$ with tail $\sigma_{k-1}$
such that the \odepth\ in~$C$ of $\sigma_{k-1}$ is $k-1$. 
Then $\depthin{\mon_C}{q_{C,\sigma_{k-1}}}={k-1}$. 
Since $p_h(\sigma_{k-1})$
belongs to $\andposv$,  $\sigma_{k-1}$ is of type~V or~H in $C$,
so $q_{C,\sigma_{k-1}}=p_v(\sigma_{k-1})$.
From Lemma~\ref{l.odepths}~(\ref{l.i.vone}), it follows, since $k-1$ is odd,
that $\sigma_{k-1}$ is of type~H. 
Since $p_v(\sigma_{k-1})>p_h(\sigma_{k-1})=\alpha$,
it follows that $\depthin{\mon_C}{p_h(\sigma_{k-1})}\geq k$ (in fact equality
holds as is easily seen).
Since $x$ dominates $\mon_C$,   there exists, by
\cite[Lemma~4.5]{kr},  $\beta$ in $\mon^k_x$ that dominates 
$p_h(\sigma_{k-1})=\alpha$.\hfill$\Box$


\mysubsection{Proof that $\ophi\opi=identity$}\mylabel{ss.pf1.p.4.2.kr}%
\mcomment{Pre-requisites: \ref{p.4.1.kr}, \ref{l.phipi}}
Let $\mon$ be a monomial in $\posv$ and let $\opi=(w,\mon')$.
We need to show that $\ophi$ applied to the pair $(w,\mon')$ gets us back to~$\mon$.
We know from~(3) of Proposition~\ref{p.4.1.kr} that $w$ \ortho-dominates $\mon'$,
so $\ophi$ can indeed be applied to the pair~$(w,\mon')$.   

The main ingredients of the proof are
the corresponding assertion in the case of Grassmannian~\cite[Proposition~4.2]{kr}
and the following claim which we will presently prove:
\[	(\mon')_{w,j,j+1}= \monj'\textup{(up)}\cup\mon_{j+1}'\textup{(up)}
\quad \textrm{ for every odd integer $j$}
						\]
Let us first see how the assertion follows assuming the truth of the claim,
by tracing the steps involved in applying $\ophi$ to $(w,\mon')$.
From the claim it follows that when we partition $\mon'$ into pieces
(see \S\ref{s.ophi}), we
get $\mon_j'\up\cup\mon_{j+1}'\up$ (for odd integers $j$).   Adding the
mirror images will get us to $\monjprime\cup\monjponeprime$.  
From Corollary~\ref{c.4.16.kr} it follows that $w_{j,j+1}$ is exactly
the element of $I(d,2d)$ obtained by acting $\pi$ on $\mon_j\cup\mon_{j+1}$. Now, since
$\phi\circ\pi=\textrm{identity}$,  it follows that on application of $\phi$
to $(w_{j,j+1},\monjprime\cup\monjponeprime)$ we obtain $\mon_j\cup\mon_{j+1}$.
By twisting the two diagonal elements in $\mon_j\cup\mon_{j+1}$ (if they exist
at all) and removing the elements below the diagonal~$\diagv$,  we get back
$\monjpr$.     Taking the union of $\monjpr$ (over odd integers $j$),  we get back $\mon$.

Thus we need only prove the claim.
Since $\mon'$ is the union over all odd integers
of the right hand sides (this follows from the definition of~$\mon'$),
and the left hand sides as $j$ varies are mutually disjoint,   it is enough
to show that the right hand side is contained in the left hand side.
Thus we need only prove:  for $j$ an odd integer and~$\alpha$ an element in
$\monjprime\up\cup\monjponeprime\up$,  
\begin{itemize}
\item every \vchain\ in $\mon'$ with head $\alpha$ is \ortho-dominated
by $w^{j}$.
\item there exists a \vchain\ in $\mon'$ with head $\alpha$ that is not
\ortho-dominated by $w^{j+2}$.
\end{itemize}

To prove the first item,  write
$\mont=\monsupjjpone:=\{\alpha\in\mon|\odepthinmon{\alpha}\geq j\}$
and set $\opi(\mont)=(x,\mont')$.    
By Proposition~\ref{p.od.shift}, we have
$\mont_i^\pr\cup\mont_{i+1}^\pr=\mon_{i+j-1}^\pr\cup\mon_{i+j}^\pr$
for any odd integer $i$.
Thus, by the description of $\opi$,  we have
$\mont'=\cup_{k\geq j}\monkprimeup$.
By Corollary~\ref{c.4.16.kr} and 
the description of $\opi$, 
we have $x=w^{j}$.
By Corollary~\ref{c.monprime}, any $v$-chain in $\mon'$ with head
belonging to $\monjprime\up\cup\monjponeprime\up$  is contained
entirely in~$\cup_{k\geq j}\monkprime\up$.
Finally, by Proposition~\ref{p.4.1.kr}~(3) applied to $\mont$,
the desired conclusion follows.

To prove the second item we use Lemma~\ref{l.phipi}. 
Proceed by decreasing induction
on $j$.   For $j$ sufficiently large the assertion is vacuous, for $\monjprimeup
\cup\monjponeprimeup$ is empty.   To prove the induction step,  
assume that the assertion holds for
$j+2$.     If the \vchain\
consisting of the single element~$\alpha$ is not \ortho-dominated by $w^{j+2}$,
then we are done.    So let us assume the contrary.
Since the \odepth\ of the singleton \vchain~$\alpha$ is at most $2$,
it follows from Lemma~\ref{l.odomination} that $w_{j+2,j+3}$ \ortho-dominates
the \vchain~$\alpha$.    
Apply Lemma~\ref{l.phipi} with $k=j+2$.
By its conclusion,  either
there exists $\beta\in\mon_{j+2}'\up$ such that $\alpha>\beta$
or there exists $\gamma\in\mon_{j+3}'\up$ such that $p_h(\alpha)>\gamma$.

First suppose that a~$\gamma$ as above exists.
By induction,  there exists a \vchain\ in~$\mon'$---call it $D$---with 
head $\gamma$
that is not
\ortho-dominated by $w^{j+4}$.   Let $C$ be the concatenation of
$\alpha>\gamma$ and~$D$.      Since elements of $D$ have \odepth\
at least~$3$ in~$C$ 
(Lemma~\ref{l.odepths}~(\ref{l.i.odepths})),   it follows from 
Corollary~\ref{c.odomination} that
$C$ is not \ortho-dominated by~$w^{j+2}$, and we are done.

Now suppose that such a $\gamma$ does not exist.   Then a $\beta$ as
above exists.
If $\alpha>\beta$ is not \ortho-dominated by $w^{j+2}$ we are again done.
So assume the contrary.    Since the \odepth\ of $\beta$ in $\alpha>\beta$
is at least $2$,  it follows that there exists an element of $(\mon_w)_{j+3}$
that dominates $\beta$.    Applying Lemma~\ref{l.phipi} again,  this time
with $k=j+3$,    we find $\gamma'\in\mon'_{j+3}\up$ 
such that $\beta>\gamma'$.   Arguing as in the previous paragraph
with $\gamma'$ in place of $\gamma$,  we are done.\hfill$\Box$
\mysubsection{Proof that $\opi\ophi=identity$}\mylabel{ss.pf2.p.4.2.kr}
Let $\mont$ be a monomial in $\posv$ and $w$ an element of $I(d)$ that
\ortho-dominates $\mont$.     We can apply $\ophi$ to the pair $(w,\mont)$
to obtain a monomial $\montwstar$ in $\posv$.      We need to 
show that $\opi$ applied to $\montwstar$ results in
$(w,\mont)$.   The main step of the proof is to establish the following:
\begin{equation}
  \label{eq:piphi}
 \montwjjponestar =(\montwstar)_{j,j+1}^\pr
\end{equation}
(for the meaning of the left and right sides of the above equation,
see \S\ref{s.ophi} and \S\ref{s.opi} respectively).
Assuming this for the moment let us show that $\opi\circ\ophi=
\textup{identity}$.     

We trace the steps involved in applying
$\opi$ to $\montwstar$.     From Eq.~(\ref{eq:piphi}) it follows that
when we break up $\montwstar$ according to the \odepths\ of its
elements as in \S\ref{s.opi},
we get $\montwjjponestar$ (as $j$ varies over odd integers).
The next step in the application of $\opi$ is the passage from
$(\montwstar)_{j,j+1}^\pr$ to $(\montwstar)_{j,j+1}$.     This involves
replacing $\sigma_j$ by its projections and adding the mirror image
of the remaining elements of $(\montwstar)_{j,j+1}^\pr$.   It follows from
Proposition~\ref{p.t.od<2}~(3) that $\sigma_j=\delta_j$ and so
$(\montwstar)_{j,j+1}=\montwjjponehashstar$.    The next step is to
apply $\pi$ to $\montwjjponehashstar$.   Since $\pi$ is the inverse
of $\phi$ (as proved in \cite{kr}),    we have 
$\pi(\montwjjponehashstar)=(\wsubjjpone, \montwjjpone)$.   Since
$\mon_w$ and $\mont$ are respectively the unions, as $j$ varies over
odd integers, of $(\mon_w)_{j,j+1}$ and $\montwjjpone$,   we see that
$\ophi$ applied to $\montwstar$ results in $(w,\mont)$.

Thus it remains only to establish Eq.~(\ref{eq:piphi}). 
It is enough to show that the left hand side is contained in the right 
hand side, for the union over all odd $j$ of either side is $\montwstar$
and the right hand side is moreover a disjoint union. 
In other words,  we need only show that the \odepth\ 
in~$\montwstar$ of an element of $\montwjjponestar$ is either $j$ or $j+1$.
We will show,  more precisely,  that, for any element $\beta$ of $\mon_w$,
the \odepth\ in~$\montwstar$ of any
element of $\opbstar$ equals the depth in~$\monw$ of~$\beta$.
Lemma~\ref{l.piphi} will be used for this purpose.

Let $\alpha$ be an element of $\opbstar$ and 
set~$e:=\odepthin{\montwstar}{\alpha}$.
We first show, by induction on~$d:=\depthin{\monw}{\beta}$, that $e\geq d$.
There is nothing to prove in case $d=1$,  so we proceed to the induction step.
Let $\beta'$ be the element of $\mon_w$ of depth $d-1$ such that $\beta'>\beta$.
If there exists $\alpha'$ in $\ortho\piece_{\beta^\prime}^\star$ with $\alpha'>\alpha$,  the desired
conclusion follows from Corollary~\ref{c.p.odepth}~(3) and induction.
Lemma~\ref{l.piphi} says that such an $\alpha'$ exists in case $d$ is even.
So suppose that $d$ is odd and such an $\alpha'$ does not exist.
The same lemma now says that $p_h(\delta_{d-2})>\alpha$,  so the
desired conclusion follows from
Lemma~\ref{l.odepths}~(\ref{l.i.odepths}).

We now show, by induction on~$e$, that
$d\geq e$.   There is nothing to prove in case $e=1$,  so we proceed to
the induction step.   
Let $C$ be a \vchain\ in~$\montwstar$ with tail $\alpha$ and having the
good property of Proposition~\ref{p.goodCexists}.    Let $\alpha'$ be
the immediate predecessor in~$C$ of~$\alpha$.  
Let $\beta'$ in $\monw$ be such that $\alpha'\in\opbpstar$ 
(we are not claiming at the moment that $\beta'$ is unique although that is true
and follows from the assertion that we are proving, the distinguishedness
of $\monw$, and the fact that $\beta'$ dominates~$\alpha'$).     It follows from Corollary~\ref{c2.l.piphi} that $\beta'>\beta$.

Let $d':=\depthin{\monw}{\beta'}$.
It follows from Corollary~\ref{c.p.odepth}~(3) 
that $e'<e$ where $e':=\odepthin{\montwstar}{\alpha'}$.
We have, $d\geq d'+1\geq e'+1\geq (e-2)+1=e-1$,  the
first equality being justified because $\beta'>\beta$,  
the second by the induction
hypothesis, and the last 
by Lemma~\ref{l.odepths}~(\ref{l.i.odepths}).    It suffices
to rule out the possibility that $d=e-1$.    So assume $d=e-1$.
Then $d=d'+1$ and $d'=e'=e-2$.   It follows from~(\ref{l.i.odepths})
of~Lemma~\ref{l.odepths} that the \vchain\ $\alpha'>\alpha$ has \odepth~$3$
and from~(\ref{l.i.hone}) of the same lemma that $e'$ is odd.
But then we get a contradiction to Proposition~\ref{p.t.od<2}~(2)
($\alpha'$ and $\alpha$ belong to $\mont_{w,d',d'+1}^\star$).
The proof of Eq.~(\ref{eq:piphi}) is thus over.\hfill$\Box$
\mysubsection{Proof of Proposition~\ref{p.new}}\mylabel{ss.pf.p.new}
Observe that the condition (\ddag) makes sense also for a monomial 
of~$\andposv$.   By virtue of belonging to~$\id$, $v$ has
$f^*$ as an entry.  It follows from the description of the bijection
$w\leftrightarrow\mon_w$ of~\S\ref{sss.montow} that
for an element~$w$ of~$\id$ to satisfy~(\ddag)
it is necessary and sufficient that $\mon_w$ (equivalently all its
parts $\mon_{w,j,j+1}$) satisfy~(\ddag).

(1)~
%
Since $\mont$ satisfies (\ddag), so do its parts $\montwjjpone$ and
 $\montwjjponehash$ (adding the mirror image preserves~(\ddag)).
Since $\mon_{w,j,j+1}$ also satisfies (\ddag),  it follows 
from the description of the map $\phi$ of~\cite{kr}
(observe the passage from a piece $\piece$ to its ``star'' $\piece^*$)
that the $\montwjjponehashstar$ satisfy~(\ddag). Since
the ``twisting'' involved in the 
passage from $\montwjjponehashstar$ to $\montwjjponestar$ involves only
a rearrangement of row and column indices,  it follows that the
$\montwjjponestar$ satisfy~(\ddag).  Finally so also does their 
union~$\montwstar$.

(2)~The parts $\monjjponepr$ of $\mon$ clearly satisfy~(\ddag).   Therefore
so do the $\monjjpone$, for, first of all, adding the mirror image
preserves~(\ddag), and then the removal of $\sigma_j$ and addition of
its projections involves only a rearrangement of row and column indices.
It follows from description of the map $\pi$ of~\cite{kr} (observe the passage
from a block $\blockb$ to the pair $(w(\blockb),\blockb')$) that  both
$\mon_{w,j,j+1}$ and $\monprimejjpone$ satisfy~(\ddag).    Finally,
$\mon_w$ and $\mon'$ being the union (respectively) of $\mon_{w,j,j+1}$
and $\monprimejjpone\up$,   they satisfy~(\ddag).~\hfill$\Box$

\mypart{An Application}\mylabel{p.app}
As an application of the main theorem (Theorem~\ref{t.main}),
an interpretation of the multiplicity is presented.
\mysection{Multiplicity counts certain paths}%
\mylabel{s.multiplicity}\mylabel{s.application}
\checked{\S\ref{s.multiplicity}; 22 nov; 1445hrs}
\spellchecked{\S\ref{s.multiplicity}; 22 nov; 1445hrs}
Fix elements $v,w$ in $I(d)$ with $v\leq w$.   It follows from 
Corollary~\ref{c.main} that the multiplicity\index{multiplicity@multiplicity, of $X(w)$ at $\pointe^v$} of the Schubert variety
$X(w)$ in $\miso$ at the point $\pointe^v$ 
can be interpreted as the cardinality 
of a certain set of non-intersecting lattice paths.
We first illustrate this by means of two examples and then justify the
interpretation.   \mcomment{Put in correct references.} 
\ignore{
The interpretation leads immediately to
a formula for the multiplicity involving binomial determinants.
Specializing to the case $v=\epsilon = (1,\ldots,d)$,  we recover
Conca's result~\cite[Theorem~3.6]{conca}.
}     

\mysubsection{Description and illustration}\mylabel{ss.dandi}
The points of $\andposv$ can be represented, in a natural way, as the
lattice points of a grid.     The column indices of the points of the
grid are the entries of $v$ and the row indices are the entries of
$\{1,\ldots,2d\}\setminus v$.      In Figure~\ref{f.x.big} the points
of~$\posv$ and those of the diagonal in $\andposv$ are shown
(for the specific choice of~$v$ in Example~\ref{x.big}).       
The open circles represent
the points of $\mon_w\up$,   where $\mon_w$ is the distinguished monomial
in~$\andposv$ that is associated to $w$ as in \S\ref{sss.montow}.   
From each point $\beta$ of $\mon_w\up$ we
draw a vertical line upwards from $\beta$ 
and let $\bstart$\index{betastart@$\bstart$, for $\beta\in\mon_w\up$} denote the top most point
of $\posv$ on this line.       In case $\beta$ is not on the
diagonal,  draw also a horizontal line rightwards from $\beta$
and let $\bfinish$%
\index{betafinish@$\bfinish$, for $\beta\in\mon_w\up$} 
denote the right most point of $\posv$ on this line.
In case $\beta$ is on the diagonal,  then $\bfinish$%
\index{betafinish@$\bfinish$, for $\beta\in\mon_w\up$} 
is not a fixed point
but varies subject to the following constraints:
\begin{itemize}
\item $\bfinish$ is one step away from the diagonal (that is,  it is of
the form $(r,c)$, for some entry $c$ of $v$, 
where $r$ is the largest integer less than~$c^*$ that is not an entry
of~$v$);
\item the column index of $\bfinish$ is not less than that of $\beta$;
\item if $\depthin{\mon_w}{\beta}$ is odd,  
then the horizontal projection
of $\bfinish$ is the same as the vertical projection of $\gamma\finish$
where $\gamma$ is the diagonal element of $\mon_w$ of depth~$1$ more
than that of $\beta$. 
\end{itemize}

With $v$ and $w$ as in Example~\ref{x.big},  we have $\bstart=(6,3)$
and $\bfinish=(9,5)$ for $\beta=(9,3)$;  $\bstart=\bfinish=(21,20)$
for $\beta=(21,20)$;  $\bstart=(15,11)$ for the diagonal element
$\beta=(36,11)$; $\bstart=(6,1)$ for the diagonal element 
$\beta=(46,1)$.    In the particular case (of non-intersecting lattice
paths) drawn in Figure~\ref{f.big},   $\bfinish=(27,19)$ for 
$\beta=(36,11)$ and $\bfinish=(28,14)$ for $\beta=(46,1)$.

A {\em lattice path}%
\index{lattice path, from $\bstart$ to $\bfinish$, denoted $\path_\beta$} 
between a pair of such points $\bstart$
and $\bfinish$ is a sequence $\alpha_1,\ldots,\alpha_q$
of elements of $\posv$ with $\alpha_1=\bstart$ and
$\alpha_q=\bfinish$  
such that, for $1\leq j\leq q-1$, if we let $\alpha_j=(r,c)$, then 
$\alpha_{j+1}$ is either 
$(R,c)$ or $(r,C)$ where $R$ is the least element of $\{1,\ldots,2d\}\setminus v$
that is bigger than $r$ and $C$ the least element of~$v$ that is bigger
than~$c$.
Note that if  $\bstart=(r,c)$ and
$\bfinish=(R,C)$,    
then~$q$ equals
\[ |(\{1,\ldots,2d\}\setminus v)\cap\{r,r+1,\ldots,R\}| +
 | v\cap\{c,c+1,\ldots,C\}| - 1, \]
where $|\cdot|$ is used to denote cardinality.

Consider the set~$\pathsvw$\index{pathsvw@$\pathsvw$} of all
tuples~$\left(\path_\beta\right)_{\beta\in\mon_w\up}$ 
of paths where
\begin{itemize}
\item 
$\path_\beta$%
\index{lambdabeta@$\path_\beta$, for $\beta\in\mon_w\up$} 
is a lattice path between $\bstart$
and $\bfinish$ (if $\beta$ is on the diagonal,  then $\bfinish$
is allowed to vary in the manner described above);
\item
$\path_\beta$ and $\path_\gamma$ do not intersect for $\beta\neq\gamma$.
\end{itemize}
The number of such $p$-tuples, where $p:=|\monw\up|$, is the multiplicity of 
$X(w)$ at the point $\pointe^v$.
\bexample\mylabel{x.big}
Let $d=23$, 
\begin{align*}
v&=(1,2,3,4,5,11,12,13,14,19,20,22,
23,26,29,30,31,32,37,38,39,40,41),\\
w&=(4,5,9,10,14,17,18,21,23,25,27,28,
31,32,34,35,36,39,40,41,44,45,46),
\end{align*}
so that 
\begin{equation*}
\begin{split}
\mon_w = & \{(9,3), (10,2), (17,13), (18,12), (21,20),  (25,22), (27,26),\\
& \quad (28,19), (34,30), (35,29), (36,11), (44,38), (45,37), (46,1)\}
\end{split}\end{equation*}
and $\mon_w\up=$
\[
\{(9,3), (10,2), (17,13), (18,12), (21,20),
 (25,22),  (28,19), (36,11), (46,1)\}.\]
A particular element of $\pathsvw$ is depicted in Figure~\ref{f.big}.
\eexample
\setlength{\unitlength}{.46cm}
\begin{figure}[t]
\begin{picture}(13,23)(-1,-0.5)


\put(-1,-0.25){46}
\put(-1,0.75){45}
\put(-1,1.75){44}
\put(-1,2.75){43}
\put(-1,3.75){42}
\put(-1,4.75){36}
\put(-1,5.75){35}
\put(-1,6.75){34}
\put(-1,7.75){33}
\put(-1,8.75){28}
\put(-1,9.75){27}
\put(-1,10.75){25}
\put(-1,11.75){24}
\put(-1,12.75){21}
\put(-1,13.75){18}
\put(-1,14.75){17}
\put(-1,15.75){16}
\put(-1,16.75){15}
\put(-1,17.75){10}
\put(-1,20.75){~7}
\put(-1,19.75){~8}
\put(-1,18.75){~9}
\put(-1,21.75){~6}

\linethickness{0.1mm}
\multiput(0,22)(0,-1){5}{\line(1,0){4}}
\multiput(0,17)(0,-1){4}{\line(1,0){8}}
\put(0,13){\line(1,0){10}}
\put(0,12){\line(1,0){11}}
\put(0,11){\line(1,0){10}}
\put(0,10){\line(1,0){9}}
\put(0,9){\line(1,0){8}}
\put(0,8){\line(1,0){7}}
\put(0,7){\line(1,0){6}}
\put(0,6){\line(1,0){5}}
\put(0,5){\line(1,0){4}}
\put(0,4){\line(1,0){3}}
\put(0,3){\line(1,0){2}}
\put(0,2){\line(1,0){1}}

\put(-0.15,22.5){1}
\put(0.85,22.5){2}
\put(1.85,22.5){3}
\put(2.85,22.5){4}
\put(3.85,22.5){5}

\put(4.65,17.5){11}
\put(5.65,17.5){12}
\put(6.65,17.5){13}
\put(7.65,17.5){14}

\put(8.65,13.5){19}
\put(9.65,13.5){20}

\put(10.65,12.5){22}
\put(11.65,12.5){23}

\linethickness{0.1mm}
\put(0,1){\line(0,1){21}}
\put(1,2){\line(0,1){20}}
\put(2,3){\line(0,1){19}}
\put(3,4){\line(0,1){18}}
\put(4,5){\line(0,1){17}}

\put(5,6){\line(0,1){11}}
\put(6,7){\line(0,1){10}}
\put(7,8){\line(0,1){9}}
\put(8,9){\line(0,1){8}}

\put(9,10){\line(0,1){3}}
\put(10,11){\line(0,1){2}}

\thicklines
\color{blue}\put(0,0){\circle{0.4}}
\put(0,22){\circle*{0.3}}
\put(5,5){\circle{0.4}}
\put(5,17){\circle*{0.3}}

\put(0,22){\line(0,-1){4}}
\put(0,18){\line(1,0){2}}
\put(2,18){\line(0,-1){4}}
\put(2,14){\line(1,0){1}}
\put(3,14){\line(0,-1){1}}
\put(3,13){\line(1,0){1}}
\put(4,13){\line(0,-1){1}}
\put(4,12){\line(1,0){1}}
\put(5,12){\line(0,-1){1}}
\put(5,11){\line(1,0){1}}
\put(6,11){\line(0,-1){2}}
\put(6,9){\vector(1,0){2}}

\put(5,17){\line(0,-1){2}}
\put(5,15){\line(1,0){1}}
\put(6,15){\line(0,-1){3}}
\put(6,12){\line(1,0){1}}
\put(7,12){\line(0,-1){1}}
\put(7,11){\line(1,0){1}}
\put(8,11){\line(0,-1){1}}
\put(8,10){\vector(1,0){1}}
\color{red}
\put(9,9){\circle{0.4}}
\put(9,13){\circle*{0.3}}
\put(11,11){\circle{0.4}}
\put(11,12){\circle*{0.3}}

\put(9,13){\line(0,-1){1}}
\put(9,12){\line(1,0){1}}
\put(10,12){\vector(0,-1){1}}

\color{yellow}
\put(1,18){\circle{0.4}}
\put(1,22){\circle*{0.3}}
\put(1,22){\line(0,-1){2}}
\put(1,20){\line(1,0){1}}
\put(2,20){\line(0,-1){1}}
\put(2,19){\line(1,0){1}}
\put(3,19){\line(0,-1){1}}
\put(3,18){\vector(1,0){1}}

\color{cyan}
\put(2,19){\circle{0.4}}
\put(2,22){\circle*{0.3}}
\put(2,22){\line(1,0){1}}
\put(3,22){\line(0,-1){2}}
\put(3,20){\line(1,0){1}}
\put(4,20){\vector(0,-1){1}}

\color{green}
\put(6,14){\circle{0.4}}
\put(6,17){\circle*{0.3}}
\put(6,17){\line(0,-1){1}}
\put(6,16){\line(1,0){1}}
\put(7,16){\line(0,-1){2}}
\put(7,14){\vector(1,0){1}}

\color{magenta}
\put(7,15){\circle{0.4}}
\put(7,17){\circle*{0.3}}
\put(7,17){\line(1,0){1}}
\put(8,17){\vector(0,-1){2}}

\color{black}
\put(10,13){\circle{0.4}}
\put(10,13){\circle*{0.3}}
\linethickness{0.05mm}
\multiput(0,0)(1,1){12}{\line(0,1){1}}
\multiput(0,1)(1,1){12}{\line(1,0){1}}
\end{picture}
\caption{An element of $\pathsvw$ with $v$ and $w$ as in Example~\ref{x.big}}
\label{f.big}\label{f.x.big}
\end{figure}

\bexample\mylabel{x.two}\mylabel{x.small} Figure~\ref{f.small}
shows all the elements of $\pathsvw$ in the following simple case:
\[
d=7, \quad v=(1,2,3,4,7,9,10), \quad \mbox{ and } \quad 
w=(4,6,7,10,12,13,14). \]
We have $\mon_w=\{(6,3),(12,9),(13,2),(14,1)\}$, 
$\mon_w\up=\{(6,3),(13,2)(14,1)\}$.   
There are $15$ elements in~$\pathsvw$
and thus the multiplicity in this case is $15$.
\eexample
\setlength{\unitlength}{.5cm}
\begin{figure}
\definecolor{one}{rgb}{0,0,1}
\definecolor{two}{rgb}{0,1,0}
\definecolor{three}{rgb}{1,0,0}
\begin{picture}(24,21)

\newsavebox{\grid}  
\savebox{\grid}(0,6){
\linethickness{0.1mm}
\multiput(0,1)(0,0){1}{\line(0,1){5}}
\multiput(1,2)(0,0){1}{\line(0,1){4}}
\multiput(2,3)(0,0){1}{\line(0,1){3}}
\multiput(3,4)(0,0){1}{\line(0,1){2}}
\multiput(0,4)(0,1){3}{\line(1,0){3}}
\multiput(0,3)(0,0){1}{\line(1,0){2}}
\multiput(0,2)(0,0){1}{\line(1,0){1}}
\linethickness{0.05mm}
\multiput(0,0)(1,1){4}{\line(0,1){1}}
\multiput(0,1)(1,1){4}{\line(1,0){1}}
%
\color{one}
\put(0,0){\circle{0.30}}
\put(0,6){\circle*{0.20}}
%
\color{two}
\put(1,1){\circle{0.30}}
\put(1,6){\circle*{0.20}}
%
\color{three}
\put(2,5){\circle{0.30}}
\put(2,6){\circle*{0.20}}
}      
\multiput(0,0)(5,0){5}{\usebox{\grid}}
\multiput(0,7.5)(5,0){5}{\usebox{\grid}}
\multiput(0,15)(5,0){5}{\usebox{\grid}}

\thicklines
\color{three}
\multiput(2,21)(5,0){5}{\line(0,-1){1}}
\multiput(2,20)(5,0){5}{\vector(1,0){1}}

\multiput(2,13.5)(5,0){5}{\line(1,0){1}}
\multiput(3,13.5)(5,0){5}{\vector(0,-1){1}}
\multiput(2,6)(5,0){5}{\line(1,0){1}}
\multiput(3,6)(5,0){5}{\vector(0,-1){1}}
\color{two}
\multiput(1,21)(0,-7.5){2}{\line(0,-1){3}}
\multiput(1,18)(0,-7.5){2}{\vector(0,-1){1}}
\multiput(6,21)(0,-7.5){2}{\line(0,-1){3}}
\multiput(6,18)(0,-7.5){2}{\vector(1,0){1}}

\multiput(11,21)(0,-7.5){2}{\line(0,-1){2}}
\multiput(11,19)(0,-7.5){2}{\line(1,0){1}}
\multiput(12,19)(0,-7.5){2}{\vector(0,-1){1}}

\multiput(16,21)(0,-7.5){2}{\line(0,-1){2}}
\multiput(16,19)(0,-7.5){2}{\line(1,0){1}}
\multiput(17,19)(0,-7.5){2}{\vector(0,-1){1}}

\multiput(21,21)(0,-7.5){2}{\line(0,-1){2}}
\multiput(21,19)(0,-7.5){2}{\line(1,0){1}}
\multiput(22,19)(0,-7.5){2}{\vector(1,0){1}}

\multiput(1,6)(5,0){5}{\line(0,-1){1}}
\multiput(1,5)(5,0){5}{\line(1,0){1}}
\multiput(2,5)(5,0){5}{\line(0,-1){1}}

\multiput(2,4)(5,0){3}{\vector(0,-1){1}}
\multiput(17,4)(5,0){2}{\vector(1,0){1}}
\color{one}
\multiput(0,21)(0,-7.5){2}{\line(0,-1){4}}
\multiput(0,17)(0,-7.5){2}{\vector(0,-1){1}}

\put(0,6){\line(0,-1){2}}
\put(0,4){\line(1,0){1}}
\put(1,4){\vector(0,-1){2}}

\multiput(5,21)(0,-7.5){3}{\line(0,-1){4}}

\multiput(5,17)(0,-7.5){3}{\vector(1,0){1}}

\multiput(10,21)(0,-7.5){2}{\line(0,-1){4}}
\multiput(10,17)(0,-7.5){2}{\vector(1,0){1}}

\multiput(15,21)(0,-7.5){3}{\line(0,-1){3}}
\put(10,6){\line(0,-1){3}}

\multiput(20,21)(0,-7.5){2}{\line(0,-1){3}}

\multiput(15,18)(0,-7.5){3}{\line(1,0){1}}
\put(10,3){\line(1,0){1}}

\multiput(20,18)(0,-7.5){2}{\line(1,0){1}}

\multiput(16,18)(0,-7.5){2}{\vector(0,-1){1}}
\put(11,3){\vector(0,-1){1}}

\multiput(21,18)(0,-7.5){3}{\vector(1,0){1}}
\put(16,3){\vector(1,0){1}}

\put(20,6){\line(0,-1){2}}
\put(20,4){\line(1,0){1}}
\put(21,4){\line(0,-1){1}}
\end{picture}
\caption{All the $15$ non-intersecting lattice paths of Example~\ref{x.small}}
\label{f.x.small}\label{f.small}
\end{figure}
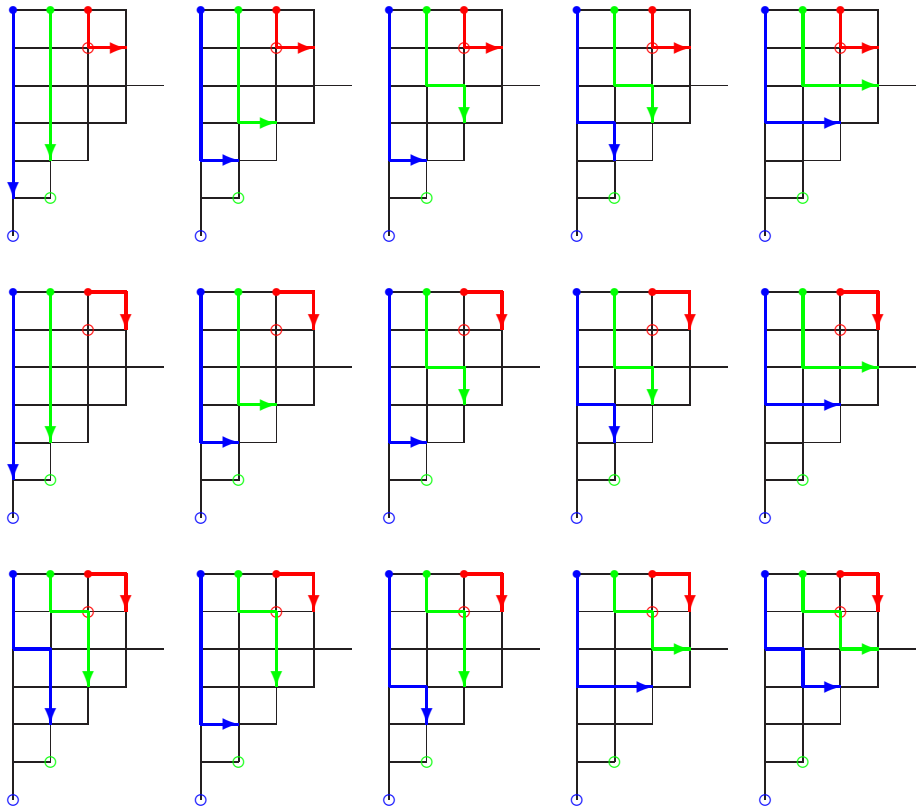
\bexample \mylabel{x.disallowedpath}Let $d=10$, 
\[v=(1,2,3,4,6,8,11,12,14,16), \quad\textup{and}\quad w=(8,9,11,14,15,16,17,18,19,20).\]
so that $\monw=\{(20,1)(19,2)(18,3),(17,4),(9,6)(15,12)\}$. 
Figure~\ref{f.disallowedpath} shows a tuple of paths that is disallowed
(meaning one that is not in $\pathsvw$).    The elements of~$\posv$
are represented as usual by a grid.
The slanted line represents the diagonal~$\diagv$. 
The solid dot represents the point of $\monw\up$ that is not on~$\diag$, 
and the crosses on~$\diag$ represent the points 
of $\monw\up$ that lie on $\diag$.    
The tuple is disallowed because the horizontal projection of the last point
of the path $\path_{\beta_1}$ is not the vertical projection of the last point
of the  path~$\path_{\beta_2}$, where $\beta_1=(20,1)$ and $\beta_2=(19,2)$
are the diagonal elements of $\monw$ of depths~$1$ and~$2$ respectively.~\eexample

\begin{figure}[!t]
\begin{center}
\includegraphics[height=86mm,width=120mm]{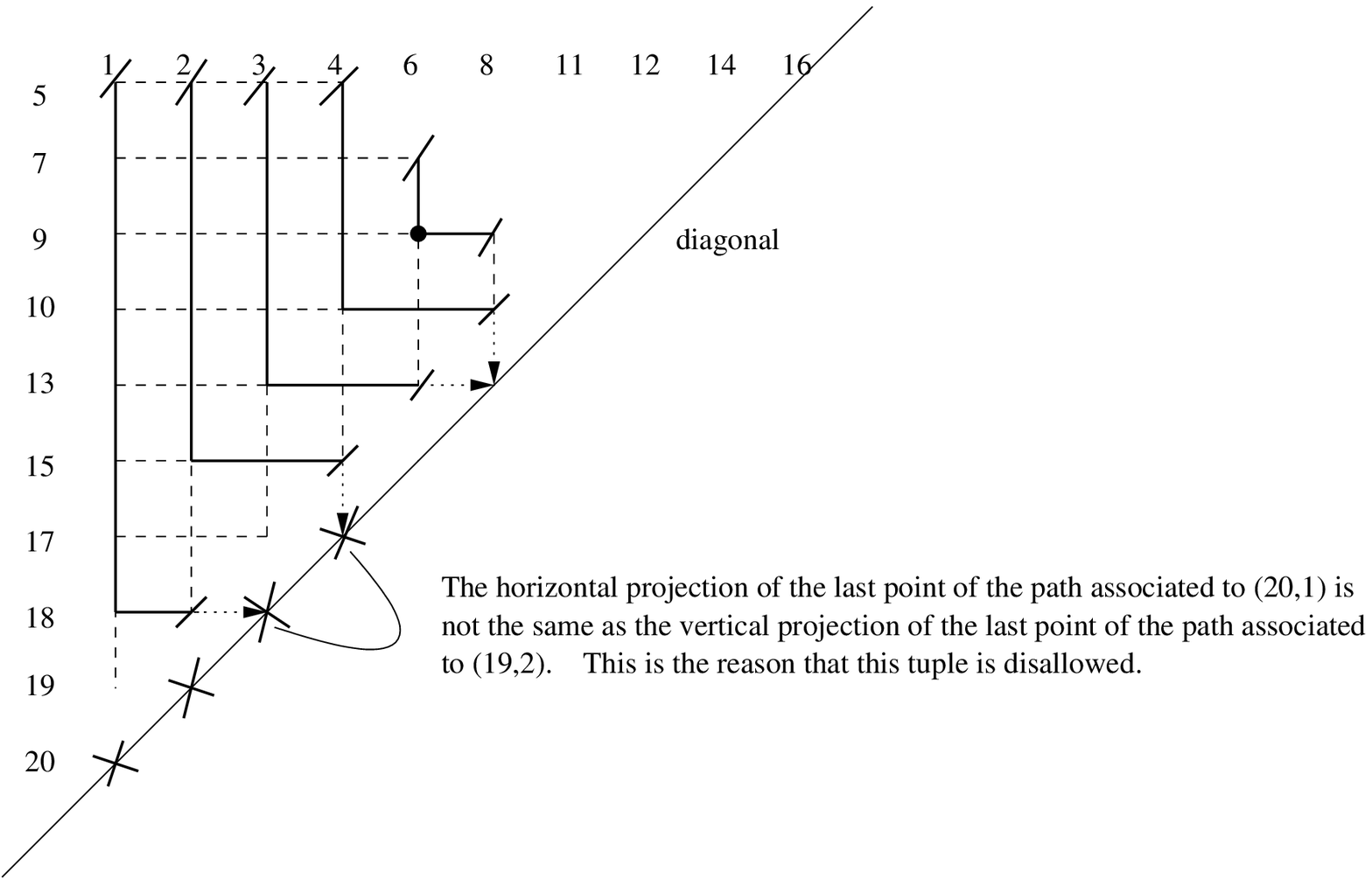}
\end{center}
\caption{\label{f.disallowedpath}
A disallowed tuple of lattice paths (see Example~\ref{x.disallowedpath})}
\end{figure}
\mysubsection{Justification for the interpretation}\mylabel{ss.justify}
We now justify the interpretation in the previous 
subsection of the multiplicity.
Corollary~\ref{c.t.main} says that the multiplicity is the number of monomials
in $\rootsv$ of maximal cardinality that are square-free and \ortho-dominated
by $w$.  Any such monomial contains $\rootsv\setminus\posv$,  for,  by the
definition of \ortho-domination,  adding or removing elements 
of $\rootsv\setminus\posv$ to or from a monomial does not alter the status 
of its \ortho-domination.      One could therefore equally well consider 
the number of monomials in $\posv$  of maximal cardinality that are 
square-free and \ortho-dominated by $w$.     We now establish a bijection 
between the set $\monomvw$\index{monw@$\monomvw$} of such monomials and the set $\pathsvw$ of 
non-intersecting lattice paths as in~\S\ref{ss.dandi}.

Each element $\path$ 
of $\pathsvw$ can be thought of, in the obvious way,
as a monomial in $\posv$.   
We will continue to denote the corresponding monomial
by~$\path$.
It is clear that the monomial $\path$ is square-free and
that all such monomials $\path$ have the same cardinality (in particular, that
if $\path_1\subseteq\path_2$ for two such monomials then $\path_1=\path_2$).
In order to establish the bijection it therefore suffices to prove 
the following proposition.
\begin{proposition}\mylabel{p.justify}
\begin{enumerate}
\item
$w$ is the element of $\id$ obtained on application of $\opi$ to 
the monomial~$\path$ (in particular (see Proposition~\ref{p.4.1.kr}),
the monomial $\path$ is \ortho-dominated by~$w$).
\item Given a monomial $\mont$ of $\posv$ that is square-free 
and \ortho-dominated by $w$,  there exists $\path$ such
that $\mont\subseteq\path$.
\end{enumerate}
\end{proposition}
\begin{myproof}
(1)~Write $\path=\left(\path_\beta\right)_{\beta\in\mon_w\up}$. 
From the description of the map~$\opi$ in~\S\ref{s.opi}, it follows that
it suffices to show that $\path_k^\pr$ (in the notation of~\S\ref{s.opi})
is the union~$\cup\path_\beta$ where~$\beta$ runs over all elements 
of depth~$k$ in~$\mon_w\up$.    In other words, it suffices to 
show that the \odepth\ in $\path$ of any element of $\path_\beta$ equals the
depth in~$\mon_w$ of $\beta$.    To prove this, we observe the following
(these assertions are easily seen to be true thinking in terms of pictures):
for fixed $\beta\in\mon_w\up$ and $\alpha\in\path_\beta$,
\begin{enumerate}
\item[(A)] For $\beta'$ in $\mon_w\up$ such that $\beta'>\beta$,  there
exists $\alpha'\in\path_{\beta'}$ such that $\alpha'>\alpha$.
\item[(B)] If $\alpha'>\alpha$ for some $\alpha'$ in $\path_{\beta'}$ for some
$\beta'$ in $\mon_w\up$,   then $\beta'>\beta$.   If, furthermore,  $\beta$ and
$\beta'$ are diagonal, their depths in $\mon_w$ are $1$ apart, and the
depth in $\mon_w$ of $\beta$ is even,  then the following is not
possible: $p_h(\alpha')$ belongs to $\andposv$ and $p_h(\alpha')>\alpha$.
\end{enumerate}
From (A) it is immediate that the \odepth~$e$ in $\path$ 
of an element $\alpha$ of
$\path_\beta$ is not less than the depth~$d$ in $\mon_w$ of $\beta$. 
We now show, by induction on~$e$, that $e\leq d$.   For $e=1$
there is nothing to show.    Suppose that $e\geq 2$.
Let~$C$ be a \vchain\ in~$\path$ having tail $\alpha$ and the good 
property of Proposition~\ref{p.goodCexists},
$\alpha'$ the immediate predecessor in~$C$ of $\alpha$, 
$e'$ the \ortho-depth of $\alpha'$ in $\path$,
$\beta'$ the element of $\mon_w\up$ such that $\alpha'\in\path_{\beta'}$,
and $d'$ the depth in~$\mon_w$ of $\beta'$.   
From Corollary~\ref{c.p.odepth}~(3) it follows that $e'\leq e-1$, so
we may apply induction.   From~(B) it follows that $d'\leq d-1$, so that,
by induction, $e'\leq d-1$.    If $e'\leq d-2$,  then we are done by
Lemma~\ref{l.odepths}~(\ref{l.i.odepths}).  So suppose that $e'=d'=d-1$.
If $d$ is odd, then the conclusion $e\leq d$ follows 
from~(\ref{l.i.odepths}) and~(\ref{l.i.hone}) of the same lemma.
In case $d$ is even, then it follows from condition~(B) and 
(\ref{l.i.odepths}) of the same lemma.

(2)   Let $\mont$ be a square-free monomial in $\posv$ that is
\ortho-dominated by~$w$.   To construct $\path$ such that $\mont\subseteq\path$,
we construct the ``components'' $\path_\beta$.    
As in~\S\ref{s.ophi},   let $\pbeta$ denote the piece of $\mont$ 
corresponding to $\beta\in\mon_w$. 
From every point belonging to $\pbeta\up$ and also from $\bstart$
carve out the South-West quadrant;    if $\beta$ is not diagonal,  then
do this also from $\bfinish$.   The boundary of the carved out
portion (intersected with $\posv$) gives a lattice path starting from $\bstart$.   In case
$\beta$ is not diagonal,   the path ends in $\bfinish$.  In this case
as well as in the case when $\beta$ is diagonal and of even depth in $\mon_w$,
we take $\path_\beta$ to be this lattice path.   In case $\beta$ is diagonal
and of odd depth in~$\mon_w$ we do the carving out from one more point
before taking $\path_\beta$ to be the boundary of the carved out region,
namely from the point that is one step away from the diagonal and whose
horizontal projection is the vertical projection of the end point of
$\path_\gamma$ where $\gamma$ is the diagonal element of $\mon_w$ of depth~$1$
more than $\beta$.   We need to justify why carving out from the
extra point is still valid,  and we do this now by 
applying Lemma~\ref{l.forphi.3}.

Let us first choose notation that is consistent with that of that lemma.
Let $\beta$ and $\gamma$ be diagonal elements in $\monw$ of 
depths $d$ and $d+1$.   Assume that~$d$ is odd.   Let the pieces of
$\mont$ corresponding to $\beta$ and $\gamma$,  when their elements
are arranged in increasing order of row and column indices,  look like this:
\[ \ldots,~(r_1,a^*),~(a,r_1^*),~\ldots;\quad\quad
 \ldots,~(r_2,b^*),~(b,r_2^*),~\ldots  \]  
It is easy to see that 
the conditions on the numbers in the above display that provide
the requisite justification are: $r_1\leq b$ 
and $a^*<b^*$ (if $\pbeta$ is empty then the justification is easy).    
To prove that $a^*<b^*$, observe that 
the diagonal elements in $\pbeta^*$ and $\piece_\gamma^*$ are respectively
$(a,a^*)$ and $(b,b^*)$,  and apply Lemma~\ref{l.forphi.2}~(2).
That $r_1\leq b$ now follows from Lemma~\ref{l.forphi.3}~(1). This finishes
the justification.

It suffices to prove the following claim: 
the lattice paths $\path_\beta$ as $\beta$ varies are non-intersecting.
Suppose that $\path_\beta$ and $\path_{\beta'}$ intersect for $\beta\neq\beta'$. Let $\alpha$ be a point of intersection.
Clearly $\beta$ dominates all elements of $\path_\beta$ 
and in particular $\alpha$;
for the same reason $\beta'$ also dominates $\alpha$.
By the distinguishedness of $\mon_w$,  we may assume without loss of generality
that $\beta'>\beta$.    It is easy to see graphically that if $\gamma$
in $\mon_w$ is such that $\beta'>\gamma>\beta$ then $\path_\gamma$ intersects
either $\path_{\beta'}$ or $\path_\beta$:  consider the open portion of $\posv$ 
``caught between'' the segment 
of $\path_{\beta'}$ from $\beta'\start$ to $\alpha$ and
the segment of $\path_\beta$ from $\bstart$ to $\alpha$;   the starting
point $\gamma\start$ of $\path_\gamma$ lives in this region but its ending
point does not (points strictly to the Northwest of $\alpha$ can neither be
of the form $\gamma\finish$ for $\gamma$ not on the diagonal nor can they
be one step away from the diagonal);  so $\path_\gamma$ must intersect one
of the two lattice path segments.     We may therefore assume that 
the depths of $\beta'$ and $\beta$ differ by~$1$.

We now apply Lemma~\ref{l.4.21.kr}.
From the construction of $\path_\beta$ it readily follows that 
$\alpha$ satisfies the hypotheses (a), (b), and (c) of that lemma.  
By the conclusion of Lemma~\ref{l.4.21.kr},
there exists $\alpha'\in\pbpstar\up$ such that $\alpha'>\alpha$.
On the other hand, it follows from the construction of $\pbpstar$ from
$\pbetap$, and from the construction of $\path_{\beta'}$ that two elements one from $\pbpstar$ and another from $\path_{\beta'}$
are not comparable.    This is a contradiction to the comparability of
$\alpha'$ and $\alpha$.
\end{myproof}


	\newcommand\citenumfont[1]{\textbf{#1}}

\bibliographystyle{bibsty-final-no-issn-isbn}
\addcontentsline{toc}{section}{References}
\ifthenelse{\equal{\finalized}{no}}{
\bibliography{abbrev,references}
}{
}
\vfill\eject
\addcontentsline{toc}{section}{Index of definitions and notation}
\begin{theindex}

  \item $>$, relation on $\posv$\dotfill 11
  \item $\leq$, partial order on $I(d,2d)$\dotfill 10

  \indexspace

  \item $\affinev$, affine patch $q_v\neq 0$ of $\miso$\dotfill 14
  \item $\alpha\down$, for $\alpha\in\andposv$\dotfill 25
  \item $\alpha^\hash$ for $\alpha$ in $\andposv$\dotfill 22
  \item $\alpha\up$, for $\alpha\in\andposv$\dotfill 25
  \item anti-domination\dotfill 18

  \indexspace

  \item $B$, a specific Borel subgroup\dotfill 8
  \item $\bfinish$, for $\beta\in\mon_w\up$\dotfill 65
  \item $\bstart$, for $\beta\in\mon_w\up$\dotfill 65
  \item block
    \subitem in the sense of ~\cite{kr}\dotfill 47
    \subitem of a monomial $\mon$ in $\posv$\dotfill 51

  \indexspace

  \item comparability, of elements of $\roots$\dotfill 29
  \item connected components of a \vchain\dotfill 23
  \item connectedness of two succcessive elements in a $v$-chain\dotfill 
		22
  \item critical element (of a $v$-chain)\dotfill 24

  \indexspace

  \item $d$, integral part of $n/2$ (unfortunately also used otherwise)\dotfill 
		7, 9
  \item degree, of a monomial\dotfill 10
  \item degree, of a standard monomial\dotfill 16
  \item $\delta_j$, for $j$ odd\dotfill 44
  \item depth (of an element $\alpha$ in a monomial $\mon$ in~$\andposv$) = $\depthin{\mon}{\alpha}$\dotfill 
		46
  \item depth (of a monomial $\mon$ in $\andposv$)\dotfill 46
  \item diagonal, $\diagv $\dotfill 10, 11
  \item distinguished (a subset of $\andposv$)\dotfill 21
  \item domination (among elements of $\roots$)\dotfill 29
  \item domination map\dotfill 18

  \indexspace

  \item $e_1,\ldots,e_n$, a specific basis of $V$\dotfill 8
  \item $\pointe^v$, $T$-fixed point\dotfill 10

  \indexspace

  \item $\langle\ ,\ \rangle$, bilinear form on $V$\dotfill 7
  \item $f_\theta:=q_\theta/q_v$\dotfill 14

  \indexspace

  \item head, of a $v$-chain\dotfill 11
  \item horizontal projection $p_h(\alpha)$\dotfill 22

  \indexspace

  \item $\id$\dotfill 10
  \item $I(d,2d)$\dotfill 10
  \item $i\even$, for an integer $i$\dotfill 29
  \item $i\odd$, for an integer $i$\dotfill 29
  \item intersection (of a monomial in a set with a subset)\dotfill 11
  \item isotropic subspace\dotfill 7
  \item $\isubn$\dotfill 8
  \item $\isubn'$\dotfill 8

  \indexspace

  \item $\field$, base field, ($\textup{characteristic}\neq2$)\dotfill 
		7, 15
  \item $k^*(:=n+1-k)$\dotfill 7, 10

  \indexspace

  \item $L$, line bundle\dotfill 13
  \item $\path_\beta$, for $\beta\in\mon_w\up$\dotfill 66
  \item lattice path, from $\bstart$ to $\bfinish$, denoted $\path_\beta$\dotfill 
		65
  \item legs of $\alpha$, for $\alpha\in\posv$\dotfill 22
  \item legs, intertwining of\dotfill 22
  \item length, of a $v$-chain\dotfill 11

  \indexspace

  \item $\miso$, orthogonal Grassmannian\dotfill 7, 8
  \item $\miso'$\dotfill 7
  \item monomial\dotfill 10
  \item $\monomvw$\dotfill 69
  \item multiplicity, of $X(w)$ at $\pointe^v$\dotfill 65
  \item $\textup{multiset}:=\textup{monomial}$\dotfill 10

  \indexspace

  \item $\andposv $\dotfill 10
  \item $n:=\dim V$, (even from \S\ref{ss.basicnotn} on)\dotfill 7, 9

  \indexspace

  \item ${\rm O}(V)$\dotfill 7
  \item $\ortho$-depth\dotfill 27
  \item $\posv $\dotfill 10
  \item $\opbstar$\dotfill 44
  \item $\ophi$\dotfill 17, 42
  \item $\opi$\dotfill 17, 34, 35
  \item $\rootsv $\dotfill 10
  \item $\ortho$-domination\dotfill 11
  \item orthogonal Grassmannian ($\miso$)\dotfill 7

  \indexspace

  \item $\pathsvw$\dotfill 66
  \item $\pbeta$\dotfill 43
  \item $\pbstar$\dotfill 43
  \item Pfaffian $q_\theta$\dotfill 14
  \item $p_h(\alpha)$, horizontal projection\dotfill 22
  \item piece of $\mont$ (see also caution)\dotfill 43
  \item $p_\theta$, Pl\"ucker coordinate\dotfill 13
  \item $p_v(\alpha)$, vertical projection\dotfill 22

  \indexspace

  \item $\myrho_{C,\alpha}$, for $\alpha$ in a \vchain~$C$\dotfill 24
  \item $q_\theta$, Pfaffian\dotfill 14

  \indexspace

  \item $\androotsv$\dotfill 10

  \indexspace

  \item $\mon$, fixed monomial in $\posv$ in~\S\ref{s.opi},~\S\ref{ss.l.ortho}\dotfill 
		\ignore{51}
  \item $\sv$, set of monomials in $\rootsv$\dotfill 18
    \subitem modifications\dotfill \see{Notation~\ref{n.explain}}{18}
  \item $\mon_C$, where $C$ is a $v$-chain\dotfill 23
  \item $\mon_{C,\alpha}$, for $\alpha$ in a \vchain~$C$\dotfill 24
  \item Schubert varieties\dotfill 7, 8
  \item $\mon\down$, for a monomial $\mon$\dotfill 25
  \item $\mon^\hash$, for monomial $\mon$ in $\andposv$ or $\androotsv$\dotfill 
		22
  \item $\sigma_k$\dotfill 34
  \item $\monjjpone$, for $\mon$ in $\posv$, $j$ odd\dotfill 34
  \item $\monjjponeext$, for $\mon$ in $\posv$, $j$ odd\dotfill 38
  \item $\monjpr$, for $\mon$ in $\posv$, $j$ odd\dotfill 34
  \item $\mon_k$, for monomial $\mon$ in $\andposv$\dotfill 46
  \item $\mon_k$, for monomial $\mon$ in $\posv$\dotfill 34
  \item $\mon_k\ext$, for monomial $\mon$ in $\posv$\dotfill 39
  \item $\mon_{k,k+1}$, for monomial $\mon$ in $\andposv$\dotfill 49
  \item $\monkpr$, for monomial $\mon$ in $\posv$\dotfill 34
  \item $\smv$\dotfill 18
    \subitem modifications\dotfill \see{Notation~\ref{n.explain}}{18}
  \item $\smvadv$\dotfill 18
  \item $\smvw$\dotfill 15
  \item ${\rm SO}(V)$\dotfill 7
  \item $\mon'$, for monomial $\mon$ in $\posv$\dotfill 35
  \item standard monomial\dotfill 14
    \subitem $v$-compatible\dotfill 15
    \subitem $w$-dominated\dotfill 14
  \item $\mon\up$, for a monomial $\mon$\dotfill 25
  \item $\monsupjjpone$, for monomial $\mon$ in $\posv$\dotfill 32
  \item $\mon^k$, for monomial $\mon$ in $\andposv$\dotfill 46
  \item $S^w(v)(m)$\dotfill 11
  \item $\monw$, $w$ in $\idd$ or $I(d,n)$\dotfill 21, 48
  \item $\monwsubjjpone$, $w$ in~$\id$, $j$ odd\dotfill 42
  \item $\mon_w^j$, $w$ in~$\id$, $j$ odd\dotfill 42
  \item symmetric (monomial of $\andposv$)\dotfill 22

  \indexspace

  \item $T$, a specific maximal torus\dotfill 8
  \item $\tv$, set of monomials in $\posv$\dotfill 18
    \subitem modifications\dotfill \see{Notation~\ref{n.explain}}{18}
  \item tail, of a $v$-chain\dotfill 11
  \item $\mont$, fixed monomial in~$\posv$ in~\S\ref{s.ophi}\dotfill 
		\ignore{42}
  \item truly orthogonal at $j$ ($j$ odd)\dotfill 34
  \item $\montwjjpone$\dotfill 42
  \item $\montwjjponehashstar$\dotfill 43
  \item $\montwjjponestar$\dotfill 44
  \item $\montwstar(:=\ophi(w,\mont))$\dotfill 44
  \item type (V, H, S), of an element in a $v$-chain\dotfill 23--24

  \indexspace

  \item $\uv$, set of monomials in $\rootsv\setminus\posv$\dotfill 18
    \subitem modifications\dotfill \see{Notation~\ref{n.explain}}{18}
  \item $u^*$, for $u\in\id$\dotfill 19

  \indexspace

  \item $V$, vector space of dimension $n$\dotfill 7
  \item $v$, fixed element of $I(d)$\dotfill 10
  \item $v$-chain\dotfill 11
  \item $v$-degree\dotfill 16
  \item vertical projection $p_v(\alpha)$\dotfill 22

  \indexspace

  \item $w(C)$ (or $w_C$), where $C$ is $v$-chain\dotfill 11
  \item $w^\hash$, for $w$ an element of $\idd$\dotfill 21
  \item $w(k)$\dotfill 34
  \item $(w,\mon')(:=\opi(\mon))$, for $\mon$ in $\posv$\dotfill 34, 35
  \item $w^*$, for $w$ an element of $\idd$\dotfill 21
  \item $\wsubjjpone$, $j$ odd\dotfill 42
  \item $\wsupj$, $j$ odd\dotfill 42

  \indexspace

  \item $X(w)$, Schubert variety\dotfill 10
  \item $x_k$, $x^k$, $x_{k,k+1}$, for $x\in I(d,n)$\dotfill 48--49
  \item $X_{r,c}$, variable\dotfill 16

  \indexspace

  \item $Y(w)(:=X(w)\cap\affinev)$\dotfill 15

\end{theindex}

\end{document}